\documentclass[article,10pt]{amsart} 
  
\usepackage[colorlinks = true, linkcolor=blue,
            urlcolor=red,
            citecolor=olive]{hyperref}
\usepackage[utf8]{inputenc}
\usepackage[T1]{fontenc}
\usepackage{amsmath}
\usepackage{graphicx}
\usepackage{dynkin-diagrams}
\usepackage{mathtools}      
\usepackage{booktabs}
\usepackage{mathrsfs}
\usepackage{tikz}
\usepackage{tikz-cd}
\usepackage{tikz-3dplot}

\usepackage{amssymb}
\usepackage{amsfonts}
\usepackage{enumitem}

\usepackage{array}

\usepackage{multicol}
\usepackage{subfiles} 

\usepackage[colorinlistoftodos]{todonotes}

\hypersetup{bookmarksdepth  =  3} 

\newtheorem{thm}{Theorem}[section]
\newtheorem{mainthm}{Theorem}

\newtheorem{lem}[thm]{Lemma}
\newtheorem{prop}[thm]{Proposition}
\newtheorem{cor}[thm]{Corollary}
\newtheorem{assumption}[thm]{Assumption}

\newtheorem{defn}[thm]{Definition}

\newtheorem{exmp}[thm]{Example}
\newtheorem{fact}[thm]{Fact}

\newtheorem{rem}[thm]{Remark}

\newtheorem{que}[thm]{Question}

\newcommand{\Z}{\mathbb{Z}}
\renewcommand{\H}{\mathbb{H}}
\newcommand{\R}{\mathbb{R}}
\newcommand{\C}{\mathbb{C}}

\newcommand{\g}{\mathfrak{g}}
\newcommand{\cc}{\mathfrak{c}}

\newcommand{\Id}{\text{Id}}

\newcommand{\sslash}{\mathbin{/\mkern-6mu/}}

\DeclareMathOperator{\Aut}{Aut}
\DeclareMathOperator{\Ad}{Ad}
\DeclareMathOperator{\ad}{ad}

\DeclareMathOperator{\End}{End}

\DeclareMathOperator{\Nest}{Nest}
\DeclareMathOperator{\tr}{tr}
\DeclareMathOperator{\Ker}{Ker}

\DeclareMathOperator{\SL}{SL}
\DeclareMathOperator{\GL}{GL}
\DeclareMathOperator{\Hom}{Hom}

\DeclareMathOperator{\SO}{SO}
\DeclareMathOperator{\SU}{SU}
\DeclareMathOperator{\Sp}{Sp}

\DeclareMathOperator{\rk}{rank}
\DeclareMathOperator{\Gr}{Gr}

\makeatletter
\DeclareRobustCommand*{\bfseries}{%
  \not@math@alphabet\bfseries\mathbf
  \fontseries\bfdefault\selectfont
  \boldmath
}
\makeatother

\newcommand{\ie}{i.e. }

\newcommand{\X}{\mathbb{X}}
\DeclareMathAlphabet{\mathpzc}{OT1}{pzc}{m}{it}
\newcommand{\harm}{\mathpzc{h}}

\newcommand{\calE}{\mathcal{E}}

\title{Anosov representations from cyclic Higgs bundles}

\copyrightinfo{2024}{Samuel Bronstein, Colin Davalo, Qiongling Li}

\author{Samuel Bronstein}
\address{Max-Planck-Institut für Mathematik in den Naturwissenschaften, Inselstrasse 22, 04103 Leipzig, Germany}
\email{bronstein@mis.mpg.de}

\author{Colin Davalo}
\address{Dipartimento di Matematica ``Giuseppe Peano'', Università degli Studi di Torino, Via Carlo Alberto 10, 10123 Torino, Italy}
\email{colinmarcthierry.davalo@unito.it}

\author{Qiongling Li}
\address{Chern Institute of Mathematics and LPMC, Nankai University, Tianjin 300071, China}
\email{qiongling.li@nankai.edu.cn}

\date{\today}
\begin{document}

\begin{abstract}
For every complex simple Lie group $G$ (not of type $E_8$) we construct the minimal Slodowy slice, which is a half-dimensional slice in the character variety of a closed surface group into $G$, containing the composition of Fuchsian representations with the minimal representation of $\mathfrak{sl}_2$ in $G$ and invariant under the mapping class group action. This slice is defined via cyclic Higgs bundles, and we show that when $G$ is not of type $E_8$, the corresponding representations are Anosov, so in particular discrete and faithful, and dominate some Fuchsian representation.

For each of these Higgs bundles, we consider the equivariant harmonic map to the symmetric space, and we show it is an embedding and a quasi-isometric embedding with regard to the induced metric. From the cyclicity of these Higgs bundles we are also able to construct an auxiliary harmonic map to a pseudoriemannian symmetric space $G/G^\sigma$ which is also an embedding.
Our techniques also apply to harmonic Higgs bundles over the unit disk satisfying a certain completeness condition. These techniques also apply to two other families of cyclic Higgs bundles.
\end{abstract}
\maketitle

\tableofcontents
\section{Introduction}

The celebrated non-Abelian Hodge correspondence, developed by \cite{hitchin1987self, Simpson88, donaldson1987twisted, corlette1988flat}, defines a bijection between geometric objects, which are representations of surface groups into a complex Lie group $G$, and algebraic objects, which are $G$-Higgs bundles over a Riemann surface. This correspondence is realized in one direction by constructing equivariant harmonic maps and on the other direction by resolving a PDE called Hitchin's equation. In particular, this correspondence does not translate every property of a Higgs bundle well to its associated representation. 
\begin{que}
    What analytic properties on a Higgs bundle imply that its associated holonomy representation is discrete and faithful?
\end{que}

 Given a complex Lie group $G$, and a Riemann surface structure $X$ on an orientable closed surface $S$ of genus $g\geq 2$, Hitchin constructed a half-dimensional slice in the character variety $\chi(\pi_1S,G)=\Hom(\pi_1S,G)\sslash G$, the Hitchin slice, that defines a connected component of the character variety $\chi(\pi_1S,G^\R)=\Hom(\pi_1S,G^\R)\sslash G^\R$ that is now called the Hitchin component \cite{Hit92}. When $G=\SL(n,\C)$, the associated representations are all discrete and faithful  \cite{Lab06,FG06}, and satisfy an additional strong property introduced by Labourie on the occasion: the Anosov property with respect to the Borel subgroup. This property that he introduced was generalized for general Parabolic subgroups, and simpler characterizations of this property were introduced by \cite{GW12,GGKW2017,KLP16,KLP17,BPS2019,KP22, T25}.

This construction of a slice by Hitchin was generalized by Collier, Sanders \cite{CS21} and Bradlow, Collier, Garcia-Prada, Gothen and Oliveira \cite{BCGGO24}, who defined given a representation of $\mathfrak{sl}_2$ into $\mathfrak{g}=\text{Lie}(G)$, also called a $\mathfrak{sl}_2$-triple, a slice of the space of $G$-Higgs bundles that they call the \emph{Slodowy slice}.
They show that the Slodowy slice of a principal $\mathfrak{sl}_2$-triple coincides with the Hitchin locus. They define a notion of \emph{magical} $\mathfrak{sl}_2$-triple for which the image of the Slodowy slice in the character variety $\chi(\pi_1S,G)$ has the property of being included in the character variety of $\chi(\pi_1S,G^\R)$ for some real form $G^\R$ of $G$.
The image of this slice is then a union of connected components of $\chi(\pi_1S,G^\R)$.
It is known that those components consist entirely of discrete and faithful representations. For some other $\mathfrak{sl}_2$-triples called \emph{odd magical}, these slices form a union of connected components of the character variety \cite{Hsiao2026}. 

However, in the case of magical $\mathfrak{sl}_2$-triples, the representations in these slices are Anosov since these representations are $\Theta$-positive for some set of roots $\Theta$ in $G^\R$~\cite{BCGGO24,BGLPW24}, a notion
introduced by Guichard--Wienhard~\cite{GW18}. Positivity implies that the representation is Anosov and hence discrete and faithful~\cite{GLW26}.
This phenomenon can only occur for connected components of character varieties consisting only of discrete and faithful representations, also called Higher Teichmüller spaces~\cite{W18}. Among $\Theta$-positive representations one has in particular Hitchin representations, maximal representations who were shown to be Anosov in \cite{BIW10, BILW05} and studied in terms of Higgs bundles in \cite{CTT19}, and $\Theta$-positive representations in $\SO_o(p,q)$ which were studied form the point of view of Higgs bundles in \cite{Collier20, ABCGGO19} and geometrically by Beyrer-Pozzetti \cite{BP24}.

\medskip

Recently, new results have emerged for Higgs bundles in Slodowy slices that are also cyclic Higgs bundles. Filip~\cite{Fil22} exhibited a family of cyclic Higgs bundles of Hodge type in $\SO(5,\C)$, whose holonomy is in $\SO(2,3)$ and satisfies the $P_2$-Anosov property. This result was extended by Zhang~\cite{Zha25} to a larger family, containing cyclic Higgs bundles that are not fixed by the $\C^\ast$-action. Other results have been obtained for cyclic Higgs bundles defining representations in $\SL(3,\R)$ \cite{BD25} and $G_2'$ \cite{DE26}. Note that \cite{DE26} uses the same technique as this paper, but on top of that leverages them to obtain control on the Anosov limit maps, and obtain consequences on the associated geometric structures. We also use a few facts about nesting cones from \cite{DE26} in our proof.

\medskip

All these are examples of cyclic Higgs bundles in Slodowy slices with discrete holonomy.

In this paper, we introduce and study, for each complex simple Lie group $G$ (except of type $E_8$), the \emph{restricted} Slodowy slice for the minimal $\mathfrak{sl}_2$-triple, that is an immersed slice of Higgs bundles that is half dimensional at smooth points. We prove that the corresponding representations are Anosov for some parabolic subgroup of $G$. We also prove that the length spectrum with respect to some weight for a representation in this slice always dominates that of some Fuchsian representation. We also apply some of these techniques to two other slices of cyclic Higgs bundles.

Our techniques are flexible enough so that they allow one to obtain results about families of harmonic Higgs bundles over the unit disk, without requiring the equivariance by a group action. Under the right hypotheses, we show that the associated minimal surface in the symmetric space is an embedding and a quasi-isometric embedding, that is $\Theta$-undistorted where $\Theta$ is a set of two roots for $G$ of type $A_n$ and a single root for other types.

\subsection{Minimal \texorpdfstring{$\mathfrak{sl}_2$}{sl2}-triple and Slodowy slice.}

We now present our main results which concern the restricted Slodowy slice of representations of a closed surface group.

Let $G$ be a simple complex Lie group with Lie algebra $\mathfrak{g}$, with finite center, of type different from $E_8$, and let $S$ be a closed oriented surface of genus $g\geq 2$.

\medskip

There exists a unique \emph{minimal} $\mathfrak{sl}_2$-triple $(e,h,f)$ in $\mathfrak{g}$ characterized by the fact that $e\in \mathfrak{g}_{-\alpha}$ where $\alpha\in \Sigma$ is the highest root of $\mathfrak{g}$. In this case the kernel of $\mathrm{ad}_e$ decomposes as $R\oplus \mathfrak{c} \oplus \langle f\rangle$ where $\mathfrak{c}$ is the Lie algebra of the centralizer $C$ of $\langle e,h,f\rangle$ in $G$. Let $T\subset G$ be the one parameter subgroup generated by $h$.

\medskip 

We consider the following family of Higgs bundles which forms the \emph{restricted Slodowy slice} associated to the minimal $\mathfrak{sl}_2$-triple, following \cite{CS21}, similar to \cite{BCGGO24}. 

Given any Riemann surface $X$ let us choose a principal $T$-bundle over $X$ such that $\mathcal{E}_T[\langle e\rangle]\simeq \mathcal{K}$\footnote{Such choices of a bundle $\mathcal{E}_T$ differ only by the choice of a $T\cap C$-bundle. Any change of the bundle $\mathcal{E}_T$ can be obtained by changing the bundle $\mathcal{E}_C$.} and a pair $(\mathcal{E}_C, \mathbf{r})$ where $\mathcal{E}_C$ is a holomorphic $C$-bundle over $X$ and $\mathbf{r}\in H^0(\mathcal{E}_T\star \mathcal{E}_C[R]\otimes \mathcal{K})$.

\begin{defn}
   A Higgs bundle over $X$ in the \emph{restricted Slodowy slice} associated to the minimal $\mathfrak{sl}_2$-triple is a Higgs bundle of the form:
$$ (\mathcal{E},\varphi)=(\mathcal{E}_T\star \mathcal{E}_C,\mathbf{e}+\mathbf{r}).$$
\end{defn}

\begin{rem}
    We call this slice \emph{restricted} as we do not allow the section $\mathbf{r}$ to have any component in $\langle f\rangle \oplus \mathfrak{c}$. 
\end{rem}

We obtain the following result in Corollary \ref{cor:Anosov closed surface groups}.

\begin{mainthm}
\label{mainthm:Anosov}
Suppose that $G$ is a simple complex group with finite center not of type $E_8$. Let $X$ be a Riemann surface structure over $S$ and $(\mathcal{E},\varphi)$ be a polystable Higgs bundle over $X$ in the restricted Slodowy slice associated to the minimal $\mathfrak{sl}_2$-triple. 

The holonomy $\rho:\pi_1S\to G$ of this Higgs bundle is a representation that is Anosov with respect to some parabolic subgroup, and in particular is discrete and faithful.
\end{mainthm}

The parabolic subgroup for which these representations are Anosov is the one associated to the set of roots $\Theta$ written in the table of Figure \ref{table:Slodowy slices}. For $G=\SL(n+1,\C)$ this parabolic subgroup is $P_{1,n}$, the stabilizer of a pair $(\ell, U)$ for some line $\ell$ and some hyperplane $U$ such that $\ell\subset U$.

In order to prove this result we construct nested cones in $\mathbb{CP}^n$ and measure their nestedness. This technique is a simplification of what was done in \cite{BD25}. The same techniques  are pushed further for representation into the split real Lie group $G_2'$ in the parallel work \cite{DE26}, to obtain information on the position of the Anosov limit set and show the existence of various connected components of the space of Anosov representations. 

\medskip

We also obtain the following domination result, also in Corollary \ref{cor:Anosov closed surface groups}.
Let $\lambda:G\to \mathfrak{a}^+$ be the Jordan projection.

\begin{mainthm}
\label{mainthm:domination}

Under the hypothesis of Theorem \ref{mainthm:Anosov}, let $\rho:\pi_1S\to G$ be the holonomy of $(\mathcal{E},\phi)$ and let $\rho_0:\pi_1S\to G$ be the holonomy of the associated Fuchsian Higgs bundle $(\mathcal{E}_T,\mathbf{e})$ over the same Riemann surface $X$. Then for some fundamental weight $\omega_V=\omega_V(G)\in \mathfrak{a}^*$, for all $\gamma\in \pi_1S$:
$$\omega_V(\lambda (\rho(\gamma)))\geq \omega_V(\lambda (\rho_0(\gamma))).$$
 Moreover if $\mathbf{r}\neq 0$, one has for some $\epsilon>0$, for all $\gamma\in \pi_1S$:
 $$\omega_V(\lambda (\rho(\gamma))\geq (1+\epsilon)\omega_V(\lambda (\rho_0(\gamma)).$$
\end{mainthm}

When $G=\SL(n+1,\C)$, the weight $\omega_V$ is the first fundamental weight:
$$\omega(\lambda_1,\cdots, \lambda_{n+1})=\lambda_1.$$
This result is motivated by the very recent work of Yao-Zhang \cite{YZ26} who proved a version of this result for representations in $\SL(3,\R)$ with $\mathbf{r}$ small. This result is also similar to the domination results for maximal representations in $\SO_o(2,n)$ \cite{CTT19} and for Hitchin representations in $\SL(3,\R)$ \cite{Tholozan_2017}.

\medskip

The slice of the character variety of Higgs bundles of the previous form for a fixed Riemann surface structure is the image of moduli space $\mathcal{M}_\mathbf{e}(TC,\langle e\rangle\oplus R)$ of $(TC,\langle e\rangle \oplus R)$-Higgs pairs satisfying that the Higgs field projection to $\mathcal E[\langle e\rangle]$ is nowhere vanishing. This moduli space is the upper extremal Toledo locus in the
space of cyclic Higgs pairs associated with the quaternionic
grading studied systematically by Garc\'ia-Prada and Gonz\'alez \cite{GG} (also \cite{GP24});
see Section \ref{sec:Toledo}. There they studied the Cayley correspondence for the lower extremal locus. Here, we give the the Cayley correspondence
for the upper extremal locus and study its properties in the following, see Proposition \ref{prop:Moduli} and Theorem \ref{thm:cayley}:

\begin{mainthm}
    \label{mainthm:Moduli}
\begin{itemize}
    \item[(i)]  $\mathcal{M}_{\mathbf e}(TC, \langle e\rangle \oplus R))$ is isomorphic to $\mathcal{M}^{\mathcal{K}^{3/2}}(C, R)$, the moduli space of polystable $\mathcal{K}^{3/2}$-twisted $(C,R)$-pairs.
    \item[(ii)] $\mathcal{M}_{\mathbf e}(TC, \langle e\rangle \oplus R))$ is a complex quasi-projective variety of pure dimension $(g-1)(\dim \mathfrak{g}-3)$. Each stable point has a neighborhood homeomorphic to $\mathbb{H}^1\sslash \Aut(\mathcal{E},\phi)$, where $\Aut(\mathcal{E},\phi)$ is finite.
\end{itemize}
\end{mainthm}

Note that $\mathbb{H}^1$ is the first hypercohomology which is a vector bundle of dimension $(g-1)(\dim \mathfrak{g}-3)$ for every stable point. For types $A$ and $D_3$, this moduli space is disconnected and admits topological invariants $d$ indexed by $[-(g-1),g-1]$; see Proposition \ref{prop:controlDegree}.

When we vary the Riemann surface structure, we can consider the joint moduli space $\mathcal{J}\mathcal{M}_{\mathbf e}(TC, \langle e\rangle \oplus R))$. Using the Non-abelian Hodge correspondence, we have a holonomy map $\mathrm{hol}$ from $\mathcal{J}\mathcal{M}_{\mathbf e}(TC, \langle e\rangle \oplus R))$ to the character variety. We obtain the following, see Theorem \ref{thm:JointModuliHalf} and Theorem \ref{thm:JointModuliproper}. 

\begin{mainthm}
\label{mainthm:ModuliSlice}
    The set of representations $\rho:\pi_1S\to G$ associated to Higgs bundles over $X$ in the restricted Slodowy slices for varying Riemann surface structure $X$ over $S$ is a closed subset in the character variety and it is the image of the map $\mathrm{hol}$ that is a half dimensional immersion at stable points into the character variety $\Hom(\pi_1S,G)\sslash G$. 
    
    Moreover, image of stable points are infinitesimally irreducible representations in $G$, and are irreducible when $G=\SL(n+1,\C)$. 
\end{mainthm}

Here a  representation $\rho$ is called infinitesimally irreducible if $\text{Lie}(Z_G(\rho))=0$, where $Z_G(\rho)$ is the centralizer of $\rho$ in $G$. It is irreducible if $Z_G(\rho)=0$.

\subsection{Higgs bundles over the disk}
 Let us fix $X=\mathbb{D}$ to be the unit disk in $\C$. A \emph{harmonic Higgs bundle} in the restricted Slodowy slice for the minimal $\mathfrak{sl}_2$-triple is a triple $(\mathcal{E},\varphi,H)$ where $(\mathcal{E},\varphi)=(\mathcal{E}_T\star \mathcal{E}_C,\mathbf{e}+\mathbf{r})$ as above and $H$ is a harmonic metric on $\mathcal{E}$, \ie a reduction of the structure group to the maximal compact $K$ of $G$ solving Hitchin's equation. We say it is \emph{parallelisable} if $H$ is compatible with the holomorphic reduction of $\mathcal{E}$ as a $TC$-principal bundle, \ie  $H$ is a reduction of structure group to the maximal compact $TC\cap K$ of $TC$.

\begin{rem}
    In the case when $X$ is compact, we prove in fact that any harmonic metric is necessarily compatible with the $TC$-structure, which comes from the cyclicity of our bundles and the special form of Hitchin's equation.
\end{rem}

 Given such a harmonic bundle, we define its associated \emph{projected metric} $g_\pi$ as the metric on $\mathbb{D}$ defined so that for $v\in T\mathbb{D}$, $\lVert v\rVert_{g_\pi}=\lVert\mathbf{e}(v^{(1,0)})\rVert_H$. We will first assume that this metric is complete. In this case a maximum principle shows that its curvature is non-positive. We obtain the following in Theorem \ref{thm:Unsdistorted On the Unit Disk}.

 \begin{mainthm}
\label{mainthm;AnosovDisk}
 Let $(\mathcal{E},\varphi,H)$ be a parallelisable harmonic bundle over $\mathbb{D}$ in the restricted Slodowy slice of the minimal $\mathfrak{sl}_2$-triple, such that the projected metric on $\mathbb{D}$ is complete and the eigen one-forms of $\varphi$ are bounded with respect to $g_{\text{hyp}}$. Then the corresponding minimal immersion $\harm:\mathbb{D}\to \X$ into the symmetric space induced by the harmonic metric $H$ and the flat connection $\nabla=\nabla^H+\varphi+\varphi^{*H}$ is $\Theta$-undistorted for some set of simple roots $\Theta\subset \Sigma$, meaning that for all $\alpha \in \Theta$ there exist $A,B>0$ such that for all $x,y\in \mathbb{D}$:
$$\alpha(\mu(\harm(x),\harm(y)))\geq A\cdot d_\pi(x,y)-B.$$
 \end{mainthm}
Here $\mu:\X\times\X\to\mathfrak{a}^+$ denotes the Cartan projection (vector-valued distance) and $d_\pi$ the distance associated to the projected metric $g_\pi$ on $\mathbb{D}$.

 \begin{rem}
     This result allows to define a universal version of the Slodowy slice, similarly to the universal Teichmüller space. Note that the condition that the eigen one-forms of $\varphi$ are bounded is equivalent to an upper bound $-\epsilon<0$ on the curvature of the projected metric $g_\pi$, see Theorem \ref{thm:NormBound3-cyclic}.
 \end{rem}

 The set $\Theta$ is described in Figure \ref{table:Slodowy slices}. Note that we obtain explicit lower bounds, with $B$ depending only on $G$ and $A$ depending only on the maximum of the curvature of the projected metric.

 \medskip

 With only the assumption that the projected metric is complete, if $X=\C$ or $\mathbb{D}$ we also prove that the harmonic map is quasi-isometrically embedded with uniform constants in Theorem \ref{thm:Undistorted weight on the disk}. 

  \begin{mainthm}
\label{mainthm:dominationDisk}
 Let $(\mathcal{E},\varphi,H)$ be a parallelisable harmonic bundle over $X=\mathbb{D}$ or $\C$ in the restricted Slodowy slice of the minimal $\mathfrak{sl}_2$-triple, such that the projected metric on $X$ is complete. Let $\harm:X\to \X$ be the minimal immersion into the symmetric space induced by the harmonic metric $H$ and the flat connection $\nabla=\nabla^H+\varphi+\varphi^{*H}$.
For some fundamental weight $\omega_V \in \mathfrak{a}^*$ depending on $G$, for all $x,y\in X$:
$$\omega_V(\mu(\harm(x),\harm(y)))\geq d_\pi(x,y).$$
 \end{mainthm}
 
Note that the induced metric of the corresponding map in the symmetric space is controlled from above by $\sqrt{3}$ times the projected metric by Corollary \ref{cor:Maximum principle}, so we have a quasi-isometric embedding.

\subsection{Other slices of cyclic bundles}

Our techniques also allow us to prove the Anosov property for two other families of Higgs bundles. We first consider $\SO(2,n+1)$ Higgs bundles that are 4-cyclic of the following shape, where $\mathcal{V}$ is rank $n-2$ orthogonal bundle, and $\mathcal{I}^2=\mathcal{O}$.
\begin{equation}\label{eq:IntroSliceSp4}
\begin{tikzcd}
&\mathcal{V}\arrow[dl, "\gamma^{\dagger}"']&\\
\mathcal{IK}\arrow[r, "1"]&\mathcal{I}\arrow[r, "1"] & \mathcal{IK}^{-1}\arrow[lu, "\gamma"']
\end{tikzcd}
\end{equation}

For closed surface groups these correspond to maximal representations into $\SO(2,n+1)$. When the underlying Riemann surface $X$ is a closed surface, the corresponding representation is already known to be Anosov, see \cite{CTT19}. In this case our argument gives a proof using more directly the shape of the Higgs bundles, rather than general results about maximal representations \cite{BILW05}.

Our method applied to Higgs bundles over non-compact Riemann surfaces allow us to prove the following results on quasi periodic spacelike surfaces into the pseudo-Riemannian hyperbolic space $\mathbb{H}^{2,n}$, introduced by Labourie--Toulisse \cite{LT}. In Theorem \ref{thm:CompleteMaximalSurface} we prove this analog of Theorem 
\ref{thm:Unsdistorted On the Unit Disk}.

\begin{mainthm}
\label{mainthm:CompleteMaximalSurface}
    Let $M\subset \mathbb{H}^{2,n}$ be a complete maximal spacelike surface whose induced metric has uniformly negative sectional curvature. Let $\mathcal{G}:M\to \X$  the corresponding Gauss map into the Riemannian symmetric space of $\SO_o(2,n+1)$. This map is $\lbrace \alpha_1\rbrace$-undistorted, meaning that there exist $A,B>0$ such that for all $x,y\in M$:
    $$\alpha_1(\mu(\mathcal{G}(x),\mathcal{G}(y))\geq Ad_\pi(x,y)-B.$$
\end{mainthm}

Here $\alpha_1$ in the root associated to the stabilizer of an isotropic line in $\R^{2,n+1}$ and $d_\pi$ is the induced metric on $M$, which can also be viewed the projected metric associated to the $4$-cyclic Higgs bundle in the Slodowy slice.

\medskip

We finally consider the following slice of nilpotent $\SL(n+4,\R)$-Higgs bundles. Let $(\mathcal V,Q_{\mathcal V})$ be an orthogonal bundle of rank $n$ and let $\mathcal L$ be a holomorphic line bundle. Let $Q_0$ be the natural orthogonal structure on $\mathcal {LK}\oplus \mathcal L\oplus \mathcal L^{-1}\oplus \mathcal L^{-1}\mathcal K^{-1}$. Let $\mu\in H^0(X,\mathcal L^{-2}\mathcal K)$ and $\nu\in H^0(X,\Hom(\mathcal V,\mathcal L\mathcal K)\otimes \mathcal K)$, and let $\nu^{\dagger}$ be the adjoint of $\nu$ with respect to $Q_{\mathcal V}$ and $Q_0$. Assume that $\nu\nu^{\dagger}=0$. Consider the following Higgs bundle:
\begin{equation}\label{eq:IntroSliceSL5}
\begin{tikzcd}
\mathcal{V}\arrow[r, "\nu"]&\mathcal{KL}\arrow[r, "1"] & \mathcal{L}\arrow[r, "\mu"] &\mathcal{L}^{-1}\arrow[r, "1"]&\mathcal{L}^{-1}\mathcal K^{-1}\arrow[bend right,llll, "\nu^{\dagger}"]&
\end{tikzcd}
\end{equation}

We prove the following analog of Theorem \ref{mainthm:Anosov} in Theorem \ref{thm:SL5}.
\begin{mainthm}
\label{mainthm:SL5}
The representations $\rho:\pi_1S \to \SL(n+4,\R)$ associated to polystable Higgs bundles of the form \eqref{eq:IntroSliceSL5} over a Riemann surface structure $X$ on $S$ are Anosov with respect to the stabilizer of a plane in $\R^{n+4}$.
\end{mainthm}

\subsection{Our techniques.}

The Higgs bundles that we consider have in common that they belong to a Slodowy slice of some $\mathfrak{sl}_2$-triple, and satisfy an additional technical property : they lie in a \emph{parallelisable} Slodowy slice\footnote{Except for our slice in $\SL(n+4,\R)$ that belongs to a "partially parallelisable" Slodowy slice.}. This implies a strong restriction on the harmonic metric, namely for compact Riemann surfaces it implies that the metric is compatible with the $TC$-structure on the bundle. This is related to some cyclicity of the underlying Higgs bundle.

\medskip

This compatibility of the metric allows us to lift the equivariant harmonic map to the symmetric space $G/K$ to a map into the space $G/(TC\cap K)$, which can geometrically be viewed as a space of pointed copies of $\mathbb{H}^2$ in $\X$. This lift that we call $Y$ is Yano-holomorphic for the correct Yano structure on $G/(TC\cap K)$, which geometrically implies that it describes a parallel distribution of planes along the associated harmonic map.

\begin{figure}[htbp!]

\begin{center}

\begin{tikzcd}[column sep=small, row sep=small,font=\small]
T^1\widetilde{X} \arrow[rr, "\mathscr{V}"]\arrow[d] 
&& \begin{array}{c} G/{(C\cap K)} \\ = \left\lbrace\substack{\strut (\mathcal{P},x,v)\text{ pointed } \mathbb{H}^2 \text{ in } \mathbb{X} \\ \strut \text{with a tangent vector $v\in T^1_x\mathcal{P}$}}\right\rbrace \end{array}\arrow[r]\arrow[d] 
& \left\{
\substack{\strut\begin{array}{c}
\substack{\strut\left\{\text{Cones in }\mathbb{P}(V)\right\} } \\
\substack{\strut\left\{\text{Filled hyperboloids in $V$}\right\}}
\end{array}}
\right\rbrace \\
\widetilde{X} \arrow[rr,"Y"]\arrow[rdd,"\harm"] \arrow[rrrdd,"u" swap,pos=.3] 
&& \begin{array}{c} G/{(TC\cap K)} \\= \left\lbrace (\mathcal{P},x) \text{ pointed } \mathbb{H}^2 \text{ in } \mathbb{X} \right\rbrace \end{array}\arrow[ddl, dotted] \arrow[ddr, dotted]
&& \\
\\
&G/K=\mathbb{X}&&  G/G^\sigma &&
\end{tikzcd}
\end{center}
  
    \caption{Illustration of all the equivariant objects constructed from a Higgs bundle in the restricted Slodowy slice for the minimal $\mathfrak{sl}_2$-triple.}
    \label{fig:Diagram}
\end{figure}

The proof of the regularity of the harmonic equivariant map into the symmetric space $\X$, Theorem \ref{mainthm:Anosov}, comes from a construction for every unit tangent vector $v\in T^1\widetilde{X}$ of a cone in $\mathbb{P}(V)$ for some linear representation $V$ of $G$, so that these cones are nested when $v$ varies along the geodesic flow of the projected metric on $\widetilde{X}$. For the domination, Theorem \ref{mainthm:domination}, we construct for every unit tangent vector $v\in T^1\widetilde{X}$ a related filled hyperboloid in $V$, that are also nested along the geodesic flow. In practice we prove directly the versions of these theorems over the unit disk, Theorems~\ref{mainthm;AnosovDisk},\ref{mainthm:dominationDisk}.

\medskip

The idea of constructing equivariant objects associated to cyclic Higgs bundles has been studied for Higgs bundles in $\SO_o(2,n+1)$ \cite{CTT19}, in $G_2'$ \cite{CT23}, in $\SO(p,q)$ with $A$-surfaces~\cite{Nie24},  isotropic $\mathbf{P}$-alternating surfaces \cite{RungiTamburelli2025}, alternating surfaces \cite{LZ1, LZ26}, and branched affine spheres~\cite{FR26}. The same construction of parallel distribution of planes can be found in~\cite{BD25}  for representations in $\SL(3,\R)$ and~\cite{DE26b,DE26} for representations into $G_2'$. Different techniques have been used to prove that cyclic Higgs bundle define Anosov representations into $\SU(1,n)$~\cite{VirgilioPhD}, $\SO(2,n+1)$ \cite{Fil22,Zha25} or at least quasi-isometric embeddings in $\text{PU}(2,n+1)$ \cite{LZ26}.

In order to develop an analogy with the theory of maximal spacelike surfaces into $\mathbb{H}^{2,n}$ \cite{CTT19,LTW24,LT}, and of affine spheres \cite{Labourie2007FlatProjective,Loftin2001AffineSpheres, BH13, BH14, DW15}, we define for minimal $\mathfrak{sl}_2$-triples a real form $G^\sigma$ of $G$ such that the map $Y$ descends to a harmonic immersion $u:\widetilde{X} \to G/G^\sigma$, whose induced metric is exactly the projected metric.

\subsection{Organization of the Paper}

In Section \ref{sec:ParallelizableSlodowySlices} we introduce the notions of parallelisable Slodowy slices in Lie algebras, illustrated with examples in the classical Lie algebras. We then discuss the corresponding restricted Slodowy slice of Higgs bundles, with examples in the classical cases. We then explain how when $G$ is not of type $E_8$, the Higgs bundles in the minimal restricted Slodowy slice define $3$-cyclic bundles in $\SL(V)$ for some linear representation.

In Section \ref{sec:Analysis} we do the necessary analysis to obtain a maximum principle, and deduce from it some explicit nestedness of cones in projective space and filled hyperboloids. We consider first $3$-cyclic bundles and then consider the two other families we study. We then reinterpret the objects we constructed in terms of Yano-holomorphic curves.

In Section \ref{sec:Dynamics} we show how nestedness of cones and of filled hyperboloids gives control on the Cartan and Jordan projections. We then deduce from it the Anosov properties, the relevant domination properties and the regularity of minimal surfaces in symmetric spaces.

In Section \ref{sec:ModuliSpace} we study the structure of the Moduli space of Higgs bundles and the character variety for closed surface groups in order to describe the slice obtained as the union of all restricted Slodowy slices for varying Riemann surface structures $X$ on $S$.

\subsection*{Acknowledgements}
We would like to thank Brian Collier in particular for promoting the idea of studying non-magical Slodowy slices, and for discussing with us. We would also like to thank Enya Hsiao and Oscar Garcia-Prada, for helpful discussions on Slodowy slices and quaternionic triples. 
C. Davalo was funded by the ERC Consolidator Grant 101124349 of Andrea Seppi. Q. Li was supported by the National Key R\&D
Program of China No. 2022YFA1006600, the Fundamental Research Funds for the Central Universities, Nankai Zhide foundation, and the Alexander von Humboldt Foundation. This project has received funding from the European Research Council (ERC) under the European Union's Horizon 2020 research and innovation programme (grant agreement No 101018839)

The authors acknowledge support of the Institut Henri Poincaré (UAR 839 CNRS-Sorbonne Université), and LabEx CARMIN (ANR-10-LABX-59-01). This project initiated during the IHP thematic semester in 2024.

Views and opinions expressed are however those of the authors only and do not
necessarily reflect those of the European Union or the European Research Council Executive Agency.
Neither the European Union nor the granting authority can be held responsible for them. 

\subsection*{AI usage}
The classification result that we added as Appendix \ref{appendix:Classification} was obtained by AI tools. We also used AI tools to help fix gaps in the proofs of Proposition~\ref{prop:StableEquiv}, Proposition~\ref{prop:properFInite}, Proposition~\ref{prop:Genreal fact Singular values}, Proposition \ref{prop:Gap singular values from nestedness}, and find references for Appendix \ref{sec:AppendixLieAlgebras}. The idea of considering filled hyperboloids and not just projective cones to control the representations was found by the authors, during a discussion with AI tools on the paper \cite{YZ26}. In these cases the writing and verification of the results were done by the authors.

Otherwise AI tools were used solely to edit the paper: correct typos, help typesetting and proof-reading. The results and ideas were obtained by the authors themselves.

\section{Parallelisable Slodowy slices}
\label{sec:ParallelizableSlodowySlices}
In this section we introduce the notion of parallelisable Slodowy slice, with the necessary background in Lie algebra, and on harmonic Higgs bundles.

\subsection{Lie algebraic Slodowy slices}\label{sec:Slodowy}
 First let us define the Slodowy slice associated to a $\mathfrak{sl}_2$-triple $(e,h,f)$ in a complex simple Lie algebra. Let $G$ be a complex simple Lie group with Lie algebra $\mathfrak{g}$, with finite center.

 \begin{defn}
     A \emph{$\mathfrak{sl}_2$-triple} is a triple $(e,h,f)$ in a Lie algebra $\mathfrak{g}$ such that:
     $$[h,e]=-2e,\; [h,f]=2f,\;[f,e]=h.$$
 \end{defn}

 The model $\mathfrak{sl}_2$-triple is the following basis of $\mathfrak{sl}_2$:
 $$e=\begin{pmatrix}
     0 & 0  \\
     1 & 0 \\
 \end{pmatrix},\; h=\begin{pmatrix}
     1 & 0  \\
     0 & -1 \\
 \end{pmatrix},\; f=\begin{pmatrix}
     0 & 1  \\
     0 & 0 \\
 \end{pmatrix}.$$

We will denote by $\mathfrak{s}=\langle e,h,f\rangle\subset \mathfrak{g}$ the span of the triple, $S$ the connected Lie group it generates $\mathfrak{c}= \Ker(\mathrm{ad}_f)\cap \Ker(\mathrm{ad}_e)$ the centralizer of $\mathfrak{s}$ in $\mathfrak{g}$ and $C\subset G$ the centralizer of $S$.
 
\begin{defn}
The Slodowy slice associated to the $\mathfrak{sl}_2$-triple $(e,h,f)$ is the affine subspace:
$$e+\Ker(\mathrm{ad}_f)\subset \mathfrak{g}.$$

We define the \emph{restricted Slodowy slice} associated to the $\mathfrak{sl}_2$-triple $(e,h,f)$ as:
$$e+R.$$

Here $R$ is the orthogonal in $\Ker(\mathrm{ad}_f)$ of $\langle f\rangle \oplus \mathfrak{c}$ for the Killing form.
\end{defn}

The $\mathfrak{sl}_2$-triples we will consider have the following remarkable property. Let $\theta$ be any Cartan involution $\theta$ of $\mathfrak{g}$ preserving $\mathfrak{c}\oplus\langle h\rangle$ such that $\theta(e)=-f$ and $\theta(h)=-h$.

\begin{defn}
We say that a $\mathfrak{sl}_2$-triple $(e,h,f)$ in $\mathfrak g$ is \emph{parallelisable} if the subspace~$R$ satisfies
$$[R,\theta(R)]\subset \langle h \rangle \oplus \mathfrak{c}\,.$$
\end{defn}

The name \emph{parallelisable} comes from the fact that for such triples, the Higgs bundles in the Slodowy slice will define a parallel distribution of planes in the symmetric space. We classify these $\mathfrak{sl}_2$-triples in Appendix \ref{appendix:Classification}.

\begin{rem}
Note that $\theta(R)$ does not depend on $\theta$ since $R$ is $\mathfrak{c}\oplus \langle h\rangle$-invariant. Note also that $\theta(R)$ is the orthogonal in $\Ker(\mathrm{ad}_e)$ of $\langle e\rangle \oplus \mathfrak{c}$.\footnote{For the hermitian pairing defined via the Killing form and the Cartan involution.}
\end{rem}

The main example of parallelisable $\mathfrak{sl}_2$-triple is the \emph{minimal} $\mathfrak{sl}_2$-triple. Recall that the space of $\mathfrak{sl}_2$-triples up to conjugation is in bijection with the set of nilpotent elements $e\in \mathfrak{g}$ up to conjugation. The space of nilpotent elements is equipped with an order such that $e\leq e'$ if $\mathrm{ad}(G)\cdot e\subset \overline{\mathrm{ad}(G)\cdot e'}$. The maximum nilpotent orbits corresponds to \emph{principal} $\mathfrak{sl}_2$-triples. The minimal nilpotent orbit corresponds to what we call the minimal $\mathfrak{sl}_2$-triples \cite{Jut08}. 

The minimal $\mathfrak{sl}_2$-triple is the $\mathfrak{sl}_2$-triple associated to long roots (for instance the longest root). Let us fix a Cartan subalgebra $\mathfrak{a}$ of $\mathfrak{g}$ and let $\Sigma$ be the associated root system. We consider the associated root space decomposition of $\mathfrak{g}$:
$$\mathfrak{g}=\mathfrak{a}\oplus\bigoplus_{\alpha\in \Sigma}\mathfrak{g}_\alpha.$$

For every root $\alpha\in \Sigma$ one can pick any non-zero element $e\in \mathfrak{g}_\alpha$ which is necessarily nilpotent. This element can be extended to an $\mathfrak{sl}_2$-triple $\lbrace e,h,f\rbrace$ so that $h\in \mathfrak{a}$ is a multiple of the co-root of $\alpha$ and $f\in \mathfrak{g}_{-\alpha}$. 

The conjugacy class of $e$ depends only on the Weyl group orbit of $\alpha$. For simple Lie algebras, two roots are in the same Weyl group orbit if and only if they have the same length, and there are at most two possible lengths. Hence one can define the \emph{long root $\mathfrak{sl}_2$-triple}, or minimal triple, constructed so that $e\in \mathfrak{g}_\alpha$ for $\alpha$ one of the long roots. If there are two lengths of roots one can additionally define the \emph{short root $\mathfrak{sl}_2$-triple} constructed so that $e\in \mathfrak{g}_\alpha$ for $\alpha$ one of the short roots. 

\medskip

Let us discuss properties of minimal $\mathfrak{sl}_2$-triples. Let $\alpha\in \Sigma$ be a long root and let us fix nonzero elements $e\in \mathfrak{g}_{-\alpha}$ and $f\in \mathfrak{g}_{\alpha}$, normalized so that $h=[f,e]\in \mathfrak{a}$ is the coroot $\alpha^\vee$ of $\alpha$ (so that $[h,e]=-2e$). For all root $\beta\in \Sigma$, the root space $\mathfrak{g}_\beta$ is an eigenline of $\mathrm{ad}_h$ with eigenvalue $2\frac{\langle\alpha,\beta\rangle}{\langle \alpha,\alpha\rangle}$, which is an integer. Since $\alpha$ has norm greater than or equal to the norm of all other roots this integer is in $\lbrace 2,1,0,-1,-2\rbrace$, and is equal to $\pm 2$ only if $\beta=\pm \alpha$.
Hence we can decompose $\mathfrak{g}$ into the corresponding weight spaces of $\mathrm{ad}_h$. This decomposition is often called a quaternionic $\Z$-grading:
$$\mathfrak{g}=\mathfrak{g}_{-2}\oplus \mathfrak{g}_{-1}\oplus \mathfrak{g}_0\oplus \mathfrak{g}_{1}\oplus \mathfrak{g}_{2}.$$

In section \ref{sec:Slodowyhiggs}, we will construct cyclic Higgs bundles related to this splitting, studied in \cite{GPR}.
We have the following:
\begin{prop}
The $h$-weight space decomposition of $\mathfrak{g}$ for a minimal $\mathfrak{sl}_2$-triple satisfies:
\begin{itemize}
    \item[(i)]$\mathfrak{g}_{2}=\mathfrak{g}_{\alpha}$, $\mathfrak{g}_{-2}=\mathfrak{g}_{-\alpha}$,
    \item[(ii)]$\mathfrak{g}_0=\mathfrak{c}\oplus \langle  h \rangle$,
    \item[(iii)]$\mathfrak{g}_{1}=R$ is the restricted Slodowy slice.
\end{itemize}
\end{prop}
Hence this quaternionic $\Z$-grading can be rewritten, for any Cartan involution $\theta$ with $\theta(h)=-h$, $\theta(e)=-f$ as:
$$\mathfrak{g}=\langle e\rangle\oplus \theta(R)\oplus (\langle h\rangle \oplus \mathfrak{c})\oplus R \oplus \langle f\rangle.$$

\begin{proof}
Part $(i)$ is a consequence of the fact that $\alpha$ is a long root. For part (ii) note that $\mathfrak{c}\subset \mathfrak{g}_0$.The space $\mathfrak{g}_0$ is equal to the sum of some root spaces and $\mathfrak{a}$. For every non-zero root $\beta\in \Sigma$ with $2\frac{\langle\alpha,\beta\rangle}{\langle \alpha,\alpha\rangle}=0$, necessarily $\beta\pm \alpha \notin \Sigma$ so $\mathfrak{g}_\beta\subset \mathfrak{c}$. If $v\in \mathfrak{a}$ is orthogonal to $h$ one can also check that $[e,v]=[f,v]=0$ so $v\in \mathfrak{c}$. This concludes part (ii).

Finally for part (iii), note that $\mathfrak{g}_1$ decomposes as a sum of root spaces associated to $\beta\in \Sigma$ such that $2\frac{\langle\alpha,\beta\rangle}{\langle \alpha,\alpha\rangle}=1$. Then the symmetry $\sigma_\alpha$ associated to $\alpha$ sends $\beta$ to $\beta -\alpha\in \Sigma$. Hence $\Ker(\mathrm{ad}_f)$ decomposes exactly as $\mathfrak{c}\oplus \mathfrak{g}_{1}\oplus \langle f\rangle$. This decomposition is orthogonal for the Killing form, hence $R=\mathfrak{g}_{1}$.
\end{proof}

\begin{prop}
Minimal $\mathfrak{sl}_2$-triples are parallelisable.
\end{prop}

\begin{proof}
    Let us fix the Cartan involution on $\mathfrak{s}$ such that $\theta(h)=-h$ and $\theta(e)=-f$ and let us extend it to a Cartan involution of $\mathfrak{g}$. This Cartan involution preserves $\langle h\rangle \oplus \mathfrak{c}$, and we have $\theta(R)=\theta(\mathfrak{g}_{1})=\mathfrak{g}_{-1}$. Then note that $[\mathfrak{g}_1,\mathfrak{g}_{-1}]\subset\mathfrak{g}_0=\langle h\rangle \oplus \mathfrak{c}$, and hence the triple is parallelisable.
\end{proof}

In order to study in general the minimal $\mathfrak{sl}_2$-triples, we will use a linear representation $V$ of the Lie algebra $\mathfrak{g}$.

\begin{prop}
\label{prop:RepresentationV}
    Let $\mathfrak{g}$ be a simple complex Lie algebra, which is not of type $E_8$, and consider the minimal $\mathfrak{sl}_2$-triple $\lbrace e,h,f\rbrace$. There exists a linear representation $V$ of $\mathfrak{g}$ such that the action of $(e,h,f)$ are given by the following block matrices $(E,D,F)$, in a vector space decomposition $V=V_{1}\oplus V_{0}\oplus V_{-1}$ with an isomorphism $I:V_{1}\to V_{-1}$:
    \[
E= \begin{pmatrix}
0 & 0 & 0 \\
0 & 0 & 0 \\
I & 0 & 0 \\
\end{pmatrix} , \; 
D= \begin{pmatrix}
\Id & 0 & 0 \\
0 & 0 & 0 \\
0 & 0 & -\Id \\
\end{pmatrix}, \;
F= \begin{pmatrix}
0 & 0 & I^{-1} \\
0 & 0 & 0 \\
0 & 0 & 0 \\
\end{pmatrix}. 
\]  
\end{prop}

If we are given such a representation, for any $r\in R$ the matrix of the corresponding action on $V$ takes the following form:
$$\begin{pmatrix}
0 & R_1 & 0 \\
0 & 0 & R_2 \\
0 & 0 & 0 \\
\end{pmatrix}. $$

We will write $\omega_V$ the highest weight of the representation $V$, which is a fundamental weight.

\begin{proof}
To find such a representation it suffices to find a representation $V$ whose associated weights $\Omega$ satisfy $\omega(h)\in \lbrace 1,0,-1\rbrace$ for all $\omega\in \Omega$. This defines the splitting with the correct action of $D$. Finally the adjoint action of $e$ identifies the weight space associated to $\omega\in \Omega$ with $\omega(h)=1$ with $\sigma_\alpha(\omega)=\omega-\alpha\in \Omega$ satisfying $(\omega-\alpha)(h)=-1$, so it defines an isomorphism between $V_{-1}$ and $V_{1}$.

The proof of the existence of such a representation, and the exact formula for $\omega_V$ is given for classical cases in Section \ref{sec:ListLie Algebra}. For the exceptional groups different from $E_8$, one checks case by case that in the smallest irreducible linear representation the weights of $h$ are $\lbrace 1,0,-1\rbrace$, which is done in Appendix \ref{sec:AppendixLieAlgebras}. We also determine there the weight $\omega_V$ in each case. 
\end{proof}

Given a minimal $\mathfrak{sl}_2$-triple and a Cartan involution $\theta$ such that $\theta(e)=-f$ and $\theta(h)=-h$ , let us describe a particular real form $\mathfrak{g}^\sigma$ of $\mathfrak{g}$.

\medskip

Let $\mathfrak{g}_\text{odd}=\mathfrak{g}_1\oplus \mathfrak{g}_{-1}$ and $\mathfrak{g}_\text{even}=\mathfrak{g}_2\oplus \mathfrak{g}_0\oplus \mathfrak{g}_{-2}$. Let $\mathfrak{k}\oplus \mathfrak{p}=\mathfrak{g}$ be the eigenspaces of $\theta$ associated to the eigenvalue $1$ and $-1$ respectively.
This real form is defined as the following Lie subalgebra:
\begin{equation}  \label{eq:defnGsigma}\mathfrak{g}^\sigma=(\mathfrak{g}_\text{even}\cap \mathfrak{k})\oplus (\mathfrak{g}_{\text{odd}}\cap \mathfrak{p}).\end{equation}

The smallest group in $G$ admitting with this Lie algebra containing $\mathcal{L}:=(TC)\cap K$ is called $G^\sigma$, which is a subgroup of $G$ containing $\mathcal{L}:=(TC)\cap K$, where $K$ is the maximal compact of $G$. 

\begin{rem}
    In what follows we don't use a lot this homogeneous space $G/G^\sigma$. It serves mostly to complete the picture in Figure \ref{fig:Diagram}, and to provide an analogy with $\mathbb{H}^{2,n}$ in the case of maximal representations into $\SO(2,n+1)$, see \cite{CTT19}. Note that this real form exists for all parallelisable $\mathfrak{sl}_2$-triples and plays a role in their classification in Appendix \ref{appendix:Classification}.
\end{rem}

\begin{rem}
    We will write $TC$ for the product of $C$ and of the one-parameter subgroup $T\subset G$ generated by $h$, whose Lie algebra is $\langle h\rangle \oplus \mathfrak{c}$. Note that depending on the Lie group $G$ chosen, $TC$ might not be exactly $T\times C$, as $T\cap C$ can be a non-trivial discrete subgroup of $T$.
\end{rem}

\begin{rem}
    Writing any nilpotent in $\rm{SL}(n+1,\C)$ in its normal Jordan form, it is easy to check that
    the only parallelisable $\mathfrak{sl}_2$-triples in $\rm{SL}(n+1,\C)$ are the minimal ones, the principal $\mathfrak{sl}_2$-triple in $\rm{SL}(3,\mathbb \C)$ and the $\mathfrak{sl}_2$-triple corresponding to the partition $2+2$ in $\rm{SL}(4,\C)$. Refer to Appendix~\ref{appendix:Classification} for a characterization and classification of parallelisable triples.
\end{rem}

\subsection{Concrete description of the Slodowy slices.}
\label{sec:ListLie Algebra}

In this section we go through the list of all minimal $\mathfrak{sl}_2$-triples, i.e. minimal $\mathfrak{sl}_2$-triples , as well as a few other $\mathfrak{sl}_2$-triples, including all the short root $\mathfrak{sl}_2$-triples. The following tables summarize the results we obtain.

\begin{figure}[h!]
\setlength\extrarowheight{4pt}
\scalebox{0.6}{\begin{tabular}{|c|c|c|c|c|c|c|}
    \hline
     Lie group & Centralizer & Slodowy slice & & Regularity of $\harm$ & Dimensions of & Real form\\
     $G$& $C$& $R$ &$\dim(R)$& $\Theta$ & $V_1,V_0,V_{-1}$ & $\mathfrak{g}^\sigma$\\ 
    \hline
     $\SL(n+1,\C)$, ($A_n$) & 
     locally $\GL(n-1,\C)$ & $\C^{n-1}\oplus (\C^{n-1})^*$& $2n-2$& \dynkin A{*o.o*}& $(1,n-1,1)$& $\mathfrak{su}(2,n-1)$\\ 
    \hline
    $\SO(2n+1,\C)$, ($B_n$) & 
    $\SL(2,\C)\times \SO(2n-3,\C)$& $\Hom(\C^2,\C^{2n-3})$ & $4n-6$& \dynkin B{o*o.oo}& $(2,2n-3,2)$ & $\mathfrak{so}(4,2n-3)$\\ 
     \hline
    $\Sp(2n,\C)$, ($C_n$)& 
    $\Sp(2n-2,\C)$ & $\C^{2n-2}$ & $2n-2$ & \dynkin C{*oo.oo}&$(1,2n-2,1)$ & $\mathfrak{sp}(1,n-1)$\\ 
    \hline
    $\SO(2n,\C)$, ($D_n$) & 
    $\SL(2,\C)\times \SO(2n-4,\C)$& $\Hom(\C^2,\C^{2n-4})$ & $4n-8$& \dynkin D{o*.oooo}& $(2,2n-4,2)$& $\mathfrak{so}(4,2n-4)$\\ 
    \hline
        $G_2$ &  $A_1$ (locally $\SL(2,\C)$) & $S^3\C^2 $& $4$& 
    \dynkin G{o*}& $(2,3,2)$ & split real form\\ 
      \hline
            $F_4$ & 
    $C_3$ (locally $\Sp(6,\C)$)& $\ker(\bigwedge^3\C^6\to\C^6) $ & $14$ & \dynkin F{*ooo}& $(6,14,6)$ & quaternionic\\ 
    \hline
    $E_6$ & 
    $A_5$ (locally $\SL(6,\C)$)&  $\bigwedge ^3\C^6$ & $20$ & \dynkin E{o*oooo}& $(6,15,6)$ & quaternionic\\ 
    \hline
    $E_7$ & 
    $D_6$ (locally $\SO(12,\C)$) &  Half-spin representation & $32$ & \dynkin E{*oooooo}& $(12,32,12)$ & quaternionic\\ 
    \hline
    $E_8$ & 
    $E_7$ & Fundamental representation & $56$ & \dynkin E{ooooooo*} & *& quaternionic\\ 
     \hline
\end{tabular}}
\caption{Information associated to the minimal $\mathfrak{sl}_2$-triples.}
\label{table:Slodowy slices}
\end{figure}

In this first table we write for each type of Lie algebra the type of the centralizer $C$ of the minimal $\mathfrak{sl}_2$-triple. The Slodowy slice is a vector space on which $C$ acts, and we describe which representation it is. We then illustrate the set of simple roots $\alpha$ such that $\alpha(h)> 0$ (always one or two roots). Finally we describe the dimensions of the splitting of the representation $V$ adapted to the minimal $\mathfrak{sl}_2$-triple, except for the type $E_8$. See \cite{GrossWallach1996} for the definition of quaternionic real forms, and see Appendix \ref{sec:AppendixLieAlgebras} for the explaination of this table for the exceptional cases.

\begin{figure}[h!]
\setlength\extrarowheight{4pt}
\scalebox{.7}{\begin{tabular}{|c|c|c|c|}

    \hline
    Lie group & $\mathfrak{sl}_2$-triple & Is parallelisable ? & Comments\\
    \hline
     $\SL(3,\C)$ &  principal $\mathfrak{sl}_2$-triple & Yes & This triple is magical:\\
     &&&Hitchin representations in $\SL(3,\R)$. \\ 
    \hline
    $\SO(2n+1,\C)$ &  short root $\mathfrak{sl}_2$-triple & Yes & This triple is magical:\\
     &&&Maximal representations in $\SO(2,2n-1)$. \\ 
       \hline
    $\SO(2n,\C)$ & Some $\mathfrak{sl}_2$-triple & Yes & This triple is magical:\\
     &&&Maximal representations in $\SO(2,2n-2)$. \\ 
       \hline
    $G_2$ &  short root $\mathfrak{sl}_2$-triple & Yes & Representations studied in \cite{CT23, DE26b,DE26} \\ 
    &&&The Slodowy slice contains Hitchin representations in $G_2'$. \\ 
    \hline
    $\SL(n+4,\C)$ & Some $\mathfrak{sl}_2$-triple & No & A partial slice defining parallel distributions of planes. \\ 
    &&&These representations have holonomy in $\SL(n+4,\R)$. \\ 
     \hline
    $\Sp(2n,\C)$ &  short root $\mathfrak{sl}_2$-triple & No & \\
    \hline
    Type $F_4$ &  short root $\mathfrak{sl}_2$-triple & No & \\
    \hline
\end{tabular}}
\caption{Other $\mathfrak{sl}_2$-triples discussed below and in Appendix \ref{sec:AppendixLieAlgebras}.}
\label{table:Slodowy slicesbis}
\end{figure}

In this second table we describe the other $\mathfrak{sl}_2$-triples that we discuss below and in Appendix \ref{subsec:AppendixG2}. In particular some of them are magical in the sense of \cite{BCGGO24}, and the corresponding Slodowy slices of Higgs bundles parametrize Higher Teichmüller space.

\subsubsection*{\texorpdfstring{$A_n$}{An}, root \texorpdfstring{$\mathfrak{sl}_2$}{sl2} :}

We fix a family $e_i$ in $\mathfrak{a}$ such that $(e_i,e_j)=\delta_{i,j}-\frac{1}{n+1}$ and $\sum_{i=1}^{n+1}e_i=0$. The root system of type $A_n$ can be written as: 
$$\Sigma=\lbrace \alpha_{i,j}=e_i-e_j\mid 1\leq i\neq j\leq n+1\rbrace.$$

Let $\alpha=e_1-e_{n+1}$ be the highest root, and let us pick a non-zero element $e\in \mathfrak{g}_{-\alpha}$ so that we can complete the triple into $f\in \mathfrak{g}_{\alpha}$ and $h=e_1-e_{n+1}=\alpha^\vee$.

The centralizer $\mathfrak{c}$ decomposes as the sum of all root spaces associated to roots $e_i-e_j$ for $2\leq i\neq j\leq n$, as well as the kernel of $\alpha$ in $\mathfrak{a}$. This Lie algebra is isomorphic to $\mathfrak{sl}_{n-1}\oplus \C$.

\medskip

The restricted Slodowy slice $R$ is the sum of root spaces associated to roots $e_1-e_i$, $e_i-e_{n+1}$ for $2\leq i\leq n$. One can check that, for $\theta$ an adapted Cartan involution, $[R, \theta(R)]=\mathfrak{c}\oplus \langle h\rangle$. This triple is parallelisable. By identifying the weights, one can see that the representation of the centralizer $\mathfrak{sl}_{n-1}\oplus \C$ defined by $R$ is $\C^{n-1}\oplus (\C^{n-1})^*$.

\medskip

In block matrix form, these take the following form:

\[
e= \begin{pmatrix}
0 & 0 & 0 \\
0 & 0 & 0 \\
1 & 0 & 0 \\
\end{pmatrix} 
\]

\[
\mathfrak{c}=\left\lbrace \begin{pmatrix}
-\ell & 0 & 0 \\
0 & U & 0 \\
0 & 0 & -\ell \\
\end{pmatrix} \middle\vert U\in \mathfrak{gl}_{n-1}(\C),\, \text{tr}(U)=2\ell \right\rbrace
\]

\[
R_e=\left\lbrace \begin{pmatrix}
0 & \beta & 0 \\
0 & 0 & \gamma \\
0 & 0 & 0 \\
\end{pmatrix} \middle\vert \beta\in \left(\C^{n-1}\right)^*,\, \gamma\in \C^{n-1} \right\rbrace
\]

The weight $\omega_V$ is simply the first fundamental weight $e_1$.

\subsubsection*{\texorpdfstring{$A_n$}{An}, principal \texorpdfstring{$\mathfrak{sl}_2$}{sl2} :}
This triple is magical, but it is also parallelisable, which is easy to check since $R_e$ is a line. However it has the property that every Higgs bundle in the associated restricted Slodowy slice is also inside the restricted Slodowy slice for the root $\mathfrak{sl}_2$, so it is not so much of interest for us.

\subsubsection*{\texorpdfstring{$A_{n+3}$}{An+3}, other \texorpdfstring{$\mathfrak{sl}_2$}{sl2} :}
Let us consider a $\mathfrak{sl}_2$-triple in $\SL(n+4,\C)$ for $n\geq 1$ whose element $e$ can be written in a decomposition of $\C^{n+4}$ into blocks of size $(1,1,n,1,1)$ as follows. We denote by $(v_1,v_2,W,v_4,v_5)$ the corresponding decomposition, so that $e$ maps $v_1$ to $v_2$ and $v_4$ to $v_5$, i.e.\ $e$ has Jordan type $(2,2,1^n)$:
\[e=\begin{pmatrix}
0&0&0&0&0\\
1&0&0&0&0\\
0&0&0&0&0\\
0&0&0&0&0\\
0&0&0&1&0
\end{pmatrix}.\]

An element of $\mathfrak{g}_0$ commuting with $e$ acts by the same matrix on $\langle v_1,v_4\rangle$ and on $\langle v_2,v_5\rangle$, identified via $e$. Hence:
\[
\mathfrak{c}=\left\lbrace \begin{pmatrix}
a_{11} & 0 & 0 & a_{12} & 0\\
0 & a_{11} & 0 & 0 & a_{12}\\
0 & 0 & U & 0 & 0\\
a_{21} & 0 & 0 & a_{22} & 0\\
0 & a_{21} & 0 & 0 & a_{22}
\end{pmatrix} \middle\vert ((a_{ij}), U)\in\mathfrak{gl}_2(\C) \times\mathfrak{gl}_{n}(\C) \right\rbrace \cap \mathfrak{sl}_{n+4}(\C).
\]
The centralizer of the triple is the group $C$ of pairs $(A,U)\in\GL(2,\C)\times\GL(n,\C)$ with $\det(A)^2\det(U)=1$, acting as above. Its Lie algebra is isomorphic to $\mathfrak{sl}_2(\C)\oplus\mathfrak{sl}_n(\C)\oplus\C$, so that $C$ is locally isomorphic to $\C^*\times\SL(2,\C)\times \SL(n,\C)$.

\medskip

The restricted Slodowy slice $R$ therefore has dimension $4n+3$. In block matrix form it is equal to:
\[
R=\left\lbrace \begin{pmatrix}
0 & x & a & 0 & y\\
0 & 0 & 0 & 0 & 0\\
0 & c & 0 & 0 & d\\
0 & z & b & 0 & -x\\
0 & 0 & 0 & 0 & 0
\end{pmatrix} \middle\vert x,y,z\in \C;\,a,b\in (\C^n)^*;\, c,d\in \C^n \right\rbrace.
\]

This Slodowy slice decomposes with respects to the weights of $h$ as $R=R_1\oplus R_2$, where the parameters $a,b,c,d$ parametrize $R_1$ and the parameters $x,y,z$ parametrize $R_2$ 

\medskip

We fix the Cartan involution $\theta(X)=-\overline{X}^t$, which preserves $\mathfrak{c}\oplus\langle h\rangle$ and satisfies $\theta(h)=-h$, $\theta(e)=-f$. Since $[R_2,\theta(R_1)]\subset\mathfrak{g}_1$ which is disjoint from $\mathfrak{g}_0=\mathfrak{c}\oplus \langle h\rangle$, this triple is \emph{not} parallelisable. For instance, $[E_{15},\theta(E_{13})]=[E_{15},-E_{31}]=E_{35}\in\mathfrak{g}_1$, which corresponds to the pair $(y,a)$. 

\medskip

However if let us fix a subalgebra $\mathfrak{so}(n+4,\C)\subset \mathfrak{sl}(n+4,\C)$ such that $e,f,h$ are orthogonal to $\mathfrak{so}(n+4,\C)$. We look at the following:
\[
P=\left\lbrace \begin{pmatrix}
0 & 0 & a & 0 & 0\\
0 & 0 & 0 & 0 & 0\\
0 & 0 & 0 & 0 & a^t\\
0 & z & 0 & 0 & 0\\
0 & 0 & 0 & 0 & 0
\end{pmatrix} \middle\vert z\in\C,\, a\in (\C^n)^* \right\rbrace \subset R.
\]

This sub-slice is in some sense  parallelisable, as one has $[P,\theta(P)]\subset (\mathfrak{c}\oplus \mathfrak{h})\cap \mathfrak{so}(n+4,\C)$ for every Cartan involution $\theta$ compatible with $(\mathfrak{c}\oplus \langle h\rangle)\cap \mathfrak{so}(n+4,\C)$. 

\subsubsection*{\texorpdfstring{$B_n$}{Bn}, long root \texorpdfstring{$\mathfrak{sl}_2$}{sl2} :}
Recall that the root system of type $B_n$ can be written as: 
$$\Sigma=\lbrace \alpha_{\pm i,\pm j}=\pm e^*_i\pm e_j^*\mid 1\leq i\neq j\leq n\rbrace\cup \lbrace \beta^{\pm}_{i}=\pm e_i^*\mid 1\leq i\leq n \rbrace.$$

Let $\alpha=e^*_1+e_2^*$, and let us pick a non-zero element $e\in \mathfrak{g}_{-\alpha}$ so that we can complete the triple into $f\in \mathfrak{g}_{\alpha}$ and $h=e_1+e_2$.

\begin{rem}
    Here we can choose any long root that we want, and we obtain the same $\mathfrak{sl}_2$-triple. The matrices form will be the same in the end, up to conjugating by a Weyl group element.
\end{rem}

The centralizer $\mathfrak{c}$ decomposes as the sum of all root spaces associated to roots $\pm e^*_i\pm e^*_j$ for $3\leq i\neq j\leq n$ and roots $\pm e^*_i$ for $3\leq i\leq n$ and the roots $e_1^*-e_2^*$ and $-e_1^*+e_2^*$, as well as the kernel of $\alpha$ in $\mathfrak{a}$.
 This Lie algebra is isomorphic to $\mathfrak{so}_{2n-3}\oplus \mathfrak{sl}_2$.

\medskip

The restricted Slodowy slice $R_e$ is the sum of root spaces associated to roots $e^*_1\pm e^*_i$, $e^*_2 \pm e^*_i$ for $3\leq i\leq n$, as well as $e_1^*$ and $e_2^*$. One can check that, for $\theta$ an adapted Cartan involution, $[R_e, \theta(R_e)]\subset\mathfrak{c}\oplus \langle h\rangle$. This triple is hence parallelisable. The representation of the centralizer $\mathfrak{so}_{2n-3}\oplus \mathfrak{sl}_2$ defined by $R$ is $\Hom(\C^2,\C^{2n-3})$.

In block matrix form, these take the following form, in blocks of size $(2,2n-3,2)$ respectively, where $J=\begin{pmatrix}0&1\\-1&0\end{pmatrix}$:
\[
e=\begin{pmatrix}
0 & 0 & 0\\
0 & 0 & 0 \\
J & 0 & 0  
\end{pmatrix}
\]
\[
\mathfrak{c}=\left\lbrace \begin{pmatrix}
A & 0 & 0\\
0 & U & 0 \\
0 & 0 & -A^t  
\end{pmatrix} \middle\vert U\in \mathfrak{so}_{2n-3}(\C),\, A\in \mathfrak{sl}_2\right\rbrace
\]
\[
R_e=\left\lbrace \begin{pmatrix}
0 & -\beta^t & 0 \\
0 & 0 & \beta \\
0 & 0 & 0 \\
\end{pmatrix} \middle\vert \beta\in \Hom(\C^2,\C^{2n-3})\right\rbrace
\]

Here we choose a basis adapted with the Cartan decomposition, so the quadratic form preserved is:
\[
q=\begin{pmatrix}
0 & 0 & I\\
0 & I & 0 \\
I & 0 & 0 \end{pmatrix}
\]

The weight $\omega_V$ is the weight $e_1$.

\subsubsection*{\texorpdfstring{$B_n$}{Bn}, short root \texorpdfstring{$\mathfrak{sl}_2$}{sl2} (which is also magical):}

Let $\alpha=e^*_n$, and let us pick a non-zero element $e\in \mathfrak{g}_{-\alpha}$ so that we can complete the triple into $f\in \mathfrak{g}_{\alpha}$ and $h=\alpha^\vee=2e_n$.  This triple is in fact the magical triple associated to maximal representations into $\SO(2, 2n-1)$.

\medskip

The centralizer $\mathfrak{c}$ decomposes as the sum of all root spaces associated to roots $\pm e^*_i\pm e^*_j$ for $1\leq i\neq j< n$, as well as the kernel of $\alpha$ in $\mathfrak{a}$. 
 This Lie algebra is isomorphic to $\mathfrak{so}_{2n-2}$.

\medskip

The restricted Slodowy slice $R_e$ is the sum of root spaces associated to roots $e^*_n\pm e^*_i$, for $1\leq i< n$. One can check that, for $\theta$ an adapted Cartan involution, $[R_e, \theta(R_e)]\subset\mathfrak{c}\oplus \langle h\rangle$. This triple is hence parallelisable.

In block matrix form, these take the following form, in blocks of size $(n-1,1,1,1,n-1)$:
\[
e=\begin{pmatrix}
0 & 0 & 0& 0& 0 \\
0 & 0 & 0& 0& 0 \\
0 & 1 & 0& 0& 0 \\
0 & 0 & 1& 0& 0 \\
0 & 0 & 0& 0& 0 
\end{pmatrix}
\]
\[
\mathfrak{c}=\left\lbrace \begin{pmatrix}
A & 0 & 0& 0& D \\
0 & 0 & 0& 0& 0 \\
0 & 0 & 0& 0& 0 \\
0 & 0 & 0& 0& 0 \\
C & 0 & 0& 0& B 
\end{pmatrix} \middle\vert \begin{pmatrix}
A & D \\
C & B 
\end{pmatrix}\in \mathfrak{so}_{2n-2}(\C)\right\rbrace
\]

\[
R_e=\left\lbrace \begin{pmatrix}
0 & \mu^t & 0& \nu^t& 0 \\
0 & 0 & 0& 0& \nu\\
0 & 0 & 0& 0& 0 \\
0 & 0 & 0& 0& \mu \\
0 & 0 & 0& 0& 0 
\end{pmatrix} \middle\vert \mu,\nu \in \C^{n-1} 
\right\rbrace
\]

Here we choose a basis adapted with the Cartan decomposition, so the quadratic form preserved is also $q$ as before.

\subsubsection*{\texorpdfstring{$C_n$}{Cn}, long root \texorpdfstring{$\mathfrak{sl}_2$}{sl2} :}

Recall that the root system of type $C_n$ can be written as: 
$$\Sigma=\lbrace \alpha_{\pm i,\pm j}=\pm e^*_i\pm e_j^*\mid 1\leq i\neq j\leq n\rbrace\cup \lbrace \beta^{\pm}_{i}=\pm 2e_i^*\mid 1\leq i\leq n \rbrace.$$

Let $\alpha=2e^*_1$ be the highest root, and let us pick a non-zero element $e\in \mathfrak{g}_{-\alpha}$ so that we can complete the triple into $f\in \mathfrak{g}_{\alpha}$ and $h=\alpha^\vee=e_1$.

\medskip

The centralizer $\mathfrak{c}$ decomposes as the sum of all root spaces associated to roots $\pm e^*_i\pm e^*_j$ for $1< i\neq j\leq n$ and $\pm 2e^*_i$ for $1<i\leq n$, as well as the kernel of $\alpha$ in $\mathfrak{a}$. 
 This Lie algebra is isomorphic to $\mathfrak{sp}_{2n-2}$.

\medskip

The restricted Slodowy slice $R_e$ is the sum of root spaces associated to roots $e^*_1\pm e^*_i$, for $1< i\leq n$. One can check that, for $\theta$ an adapted Cartan involution, $[R_e, \theta(R_e)]\subset\mathfrak{c}\oplus \langle h\rangle$. This triple is hence parallelisable. The representation of the centralizer $\mathfrak{sp}_{2n-2}$ associated to $R$ is $\C^{2n-2}$.

\medskip

In block matrix form, these take the following form, in blocks of size $(1,2n-2,1)$:
\[
e=\begin{pmatrix}
0 & 0 & 0 \\
0 & 0 & 0 \\
1 & 0 & 0 
\end{pmatrix}
\]
\[
\mathfrak{c}=\left\lbrace \begin{pmatrix}
0 & 0 & 0 \\
0 & A & 0 \\
0 & 0 & 0 
\end{pmatrix} \middle\vert A\in \mathfrak{sp}_{2n-2}(\C)\right\rbrace
\]

\[
R_e=\left\lbrace \begin{pmatrix}
0 & -\beta^\perp & 0 \\
0 & 0 & \beta \\
0 & 0 & 0 
\end{pmatrix} \middle\vert \beta\in \C^{2n-2} 
\right\rbrace
\]
Here $\beta^\perp\in(\C^{2n-2})^*$ denotes the linear form associated to $\beta$ by the symplectic form of $\C^{2n-2}$.

\medskip

The weight $\omega_V$ is the first fundamental weight $e_1$.

\subsubsection*{\texorpdfstring{$C_n$}{Cn}, short root \texorpdfstring{$\mathfrak{sl}_2$}{sl2} (not parallelisable):}

Before starting, we will assume here that $n\geq 3$, else, $B_2=C_2$.

Let $\alpha=e^*_1+e_2^*$, and let us pick a non-zero element $e\in \mathfrak{g}_{-\alpha}$ so that we can complete the triple into $f\in \mathfrak{g}_{\alpha}$ and $h=\alpha^\vee=e_1+e_2$.

The centralizer $\mathfrak{c}$ decomposes as the sum of all root spaces associated to roots $\pm e^*_i\pm e^*_j$ for $3\leq i\neq j\leq n$ and roots $\pm 2e^*_i$ for $3\leq i\leq n$, as well as the kernel of $\alpha$ in $\mathfrak{a}$.
 This Lie algebra is isomorphic to $\mathfrak{sp}_{2n-4}\oplus \C$.

\medskip

The restricted Slodowy slice $R$ is the sum of root spaces associated to roots $e^*_1\pm e^*_i$, $e^*_2 \pm e^*_i$ for $3\leq i\leq n$ as well as $2e_1^*$ and $2e_2^*$.

One can check that, for $\theta$ an adapted Cartan involution we have $[R_e, \theta(R_e)]\not\subset\mathfrak{c}\oplus \langle h\rangle$, so this triple is \emph{not} parallelisable ! The issue is that $[\mathfrak{g}_{2e^*_1},\theta(\mathfrak{g}_{e^*_1-e^*_3})]=\mathfrak{g}_{2e^*_1-(e^*_1-e^*_3)}=\mathfrak{g}_{e^*_1+e^*_3}$ and this root space is not in the centralizer (it is actually in $R_e$). We can however decompose $R_e$ into two slices $P_1$ and $P_2$ where $P_1$ is the sum of root spaces $e^*_1\pm e^*_i$, $e^*_2 \pm e^*_i$, whose elements are denoted $\beta, \gamma$, and $P_2$ is the sum of root spaces $2e_1^*$ and $2e_2^*$, whose elements are denoted $\mu,\nu$.

In block matrix form, these take the following form, in blocks of size $(1,1,2n-4,1,1)$ respectively:
\[
e=\begin{pmatrix}
0 & 0 & 0& 0& 0 \\
0 & 0 & 0& 0& 0 \\
0 & 0 & 0& 0& 0 \\
1 & 0 & 0& 0& 0 \\
0 & 1 & 0& 0& 0 
\end{pmatrix}
\]
\[
\mathfrak{c}=\left\lbrace \begin{pmatrix}
0 & 0 & 0& 0& 0 \\
0 & \lambda & 0& 0& 0 \\
0 & 0 & A& 0& 0 \\
0 & 0 & 0& -\lambda& 0 \\
0 & 0 & 0& 0& 0   
\end{pmatrix} \middle\vert A\in \mathfrak{sp}_{2n-4}(\C),\, \lambda\in \C\right\rbrace
\]
\[
R_e=\left\lbrace \begin{pmatrix}
0 & 0 & \beta& 0& \mu \\
0 & 0 & 0& \nu& 0 \\
0 & 0 & 0& 0& \gamma \\
0 & 0 & 0& 0& 0 \\
0 & 0 & 0& 0& 0 
\end{pmatrix}   \middle\vert \beta\in (\C^{2n-4})^*, \gamma\in \C^{2n-4}, \mu,\nu\in \C\right\rbrace
\]

\begin{rem}
Here if we want to have a parallel distribution of planes we need either $\mu,\nu=0$, or $\beta,\gamma=0$.
\end{rem}

\subsubsection*{\texorpdfstring{$D_n$}{Dn}, root \texorpdfstring{$\mathfrak{sl}_2$}{sl2}:}

Recall that the root system of type $D_n$ can be written as: 
$$\Sigma=\lbrace \alpha_{\pm i,\pm j}=\pm e^*_i\pm e_j^*\mid 1\leq i\neq j\leq n\rbrace.$$

Let $\alpha=e^*_1+e_2^*$, and let us pick a non-zero element $e\in \mathfrak{g}_{-\alpha}$ so that we can complete the triple into $f\in \mathfrak{g}_{\alpha}$ and $h=e_1+e_2$.

\medskip

The centralizer $\mathfrak{c}$ decomposes as the sum of all root spaces associated to roots $\pm e^*_i\pm e^*_j$ for $3\leq i\neq j\leq n$ and roots $e^*_1-e^*_2$ and $e^*_2-e^*_1$, as well as the kernel of $\alpha$ in $\mathfrak{a}$.
This Lie algebra is isomorphic to $\mathfrak{so}_{2n-4}\oplus \mathfrak{sl}_2(\C)$.

\medskip

The restricted Slodowy slice $R_e$ is the sum of root spaces associated to roots $e^*_1\pm e^*_i$ and $e^*_2 \pm e^*_i$ for $3\leq i\leq n$. One can check that, for $\theta$ an adapted Cartan involution we have $[R_e, \theta(R_e)]\subset\mathfrak{c}\oplus \langle h\rangle$, so this triple is  parallelisable. The representation of the centralizer $\mathfrak{so}_{2n-4}\oplus \mathfrak{sl}_2$ defined by $R$ is $\Hom(\C^2,\C^{2n-4})$.
\medskip

In block matrix form, these take the following form, in blocks of size $(2,2n-4,2)$ respectively, with $J$ as in the case $B_n$:
\[e=\begin{pmatrix}
0 & 0 & 0\\
0 & 0 & 0 \\
J & 0 & 0  
\end{pmatrix}
\]
\[
\mathfrak{c}=\left\lbrace \begin{pmatrix}
A & 0 & 0\\
0 & U & 0 \\
0 & 0 & -A^t  
\end{pmatrix} \middle\vert U\in \mathfrak{so}_{2n-4}(\C),\, A\in \mathfrak{sl}_2\right\rbrace
\]
\[
R_e=\left\lbrace \begin{pmatrix}
0 & -\beta^t & 0 \\
0 & 0 & \beta \\
0 & 0 & 0 \\
\end{pmatrix} \middle\vert \beta\in \Hom(\C^2,\C^{2n-4})\right\rbrace
\]

The weight $\omega_V$ is the first fundamental weight $e_1$. 

\subsubsection*{\texorpdfstring{$D_n$}{Dn}, magical \texorpdfstring{$\mathfrak{sl}_2$}{sl2} associated to \texorpdfstring{$\SO(2,1)\subset \SO(2,2n-2)$}{SO(2,1) in SP(2,2n-2)}:}

This magical triple is not associated to a root, unlike the $B_n$ case. However the description of its Slodowy slice is completely analogous, and it is parallelisable in a very similar fashion.

\medskip

We consider the exceptional Lie algebras in Appendix \ref{sec:AppendixLieAlgebras}.

\subsection{Higgs bundles, Slodowy slices and parallelisable \texorpdfstring{$\mathfrak{sl}_2$}{sl2}-triples}
\label{sec:Slodowyhiggs}

Let $X$ be any Riemann surface, possibly non-compact. Let $G$ be a complex Lie group with Lie algebra $\mathfrak{g}$, that we assume has finite center. Following the construction in \cite{CS21, BCGGO24}, we define the Slodowy slice of Higgs bundles associated to a $\mathfrak{sl}_2$-triple. Let us fix a triple $(e,f,h)$ and let $S\subset G$ be the connected group with Lie algebra $\langle e,h,f\rangle$. Let $T\subset G$ be the one parameter subgroup generated by $h$. Let $C\subset G$ be the centralizer of the $\mathfrak{sl}_2$-triple.

\begin{defn}
A $G$-Higgs bundle is a pair $(\mathcal E, \varphi)$ of a $G$-principal holomorphic bundle $\mathcal{E}$ and a holomorphic section $\varphi\in H^0(\mathcal{E}[\mathfrak{g}]\otimes \mathcal{K})$. 
\end{defn}

Here we write $\mathcal{E}[V]$ for any linear representation $V$ of $G$ the associated vector bundle over $X$.

\begin{defn}
A $G$-Higgs bundle in the \emph{Slodowy slice} associated to $\lbrace e,h,f\rbrace$ is a pair $(\mathcal{E}_T\star \mathcal{E}_C,\mathbf{e}+\mathbf{r})$ where :

\begin{itemize}
    \item[-] $\mathcal{E}_T$ is a $T$-principal holomorphic bundle over $X$ with a never vanishing section $\mathbf{e}\in H^0(\mathcal{K}\otimes \mathcal{E}_T[\langle e\rangle])$,
    \item[-] Let $\mathcal{E}_C$ be a principal holomorphic $C$-bundle over $X$,
    \item[-] $\mathcal{E}_T\star \mathcal{E}_C$ is the associated $G$-principal bundle over $X$,
    \item[-] $\mathbf{r}$ is a holomorphic section of $\mathcal{K}\otimes\mathcal{E}_T\star \mathcal{E}_C[\Ker(\mathrm{ad}_f)]$.
\end{itemize}
We define the corresponding \emph{restricted} Slodowy slice as the space of bundles in the Slodowy slice where $\mathbf{r}$ is a section of $\mathcal{K}\otimes\mathcal{E}_T\star \mathcal{E}_C[R]$.
\end{defn}

\begin{rem}
    Note that when $X$ is a compact Riemann surface, taking $\mathbf{r}=0$ amounts, if the bundle is still polystable, to considering a representation $\rho:\pi_1X\to G$ of the form $\rho(\gamma)=\chi(\gamma)\times \rho_0(\gamma)$ where $\rho_0:\pi_1(X)\to \SL(2,\R)\to G$ is a lift of the uniformization of $X$, and $\chi:\pi_1(X)\to C\cap K$ is a representation that is directly associated to the principal bundle $\mathcal{E}_C$.
\end{rem}

Here $\mathcal{E}_T\star \mathcal{E}_C$ is the principal $G$-bundle induced by the natural map $T\times C\to TC\subset G$.
Note that $T$ preserves $\langle e\rangle\subset \mathfrak{g}$ and $TC$ preserve $R\subset \mathfrak{g}$ so the reductions $\mathcal{E}_T[\langle e\rangle]$ and $\mathcal{E}_T\star \mathcal{E}_C[R]$ are well defined.

\medskip

Let us now discuss Higgs bundles equipped with harmonic metrics. We consider here possibly non-compact surfaces and we don't discuss the existence of such metrics before Section \ref{sec:ModuliSpace}.

\begin{defn}
    A \emph{harmonic $G$-Higgs bundle} is a triple $(\mathcal{E},\varphi, H)$ where $(\mathcal E, \varphi)$ is a $G$-Higgs bundle and $H$ is a reduction of structure group of $\mathcal{E}$ to $K$, such that $D=\nabla^H+\varphi+\varphi^{*H}$ is a flat connection, where $\nabla^H$ is the Chern connection on $\mathcal{E}$ associated to $H$. 
\end{defn}

The flatness of the connection $D$ is equivalent to Hitchin's equation :
$$ F_{\nabla^H}+\varphi\wedge\varphi^{*H}=0.$$

The data of a harmonic bundle $(\mathcal{E},\varphi,H)$ defines a flat connection on $\mathcal{E}$, which defines a representation $\rho:\pi_1X\to G$ up to conjugation. Hence together with the restriction of structure group $H$, this defines a map $\harm:\widetilde{X}\to G/K$ that is harmonic and $\rho$-equivariant, that is well defined up to the action of $G$.

\medskip

Now let us focus on the case when the Higgs bundle lies in a parallelisable restricted Slodowy slice, or in a partial parallelisable Slodowy slice associated to some $P\subset \Ker(\mathrm{ad}_f)$ which is $TC$-invariant.

Let $(\mathcal{E}_T\star \mathcal{E}_C, \mathbf{e+r})$ be a $G$-Higgs bundle in the associated partial Slodowy slice. Let $K$ be a maximal compact of $G$ such that $K\cap TC$ is a maximal compact of $TC$.

\begin{defn}
\label{defn:Parallelizable Harmonic bundle}
A \emph{parallelisable harmonic bundle} over a Riemann surface $X$ in the restricted Slodowy slice associated to the triple $\lbrace e,h,f\rbrace$ is a harmonic bundle of the form $(\mathcal{E}_T\star\mathcal{E}_C,\mathbf{e+r},H)$ where the reduction $H$ to the maximal compact $K$ of $G$ is compatible with the reduction to $TC$, meaning that $H$ defines with the holomorphic reduction to $TC$ a reduction of structure group to $TC\cap K$.

A \emph{parallelisable harmonic bundle} in the partial Slodowy slice associated to $P\subset \Ker(\mathrm{ad}_f)$ for the triple $\lbrace e,h,f\rbrace$ is defined similarly, with $\mathbf{r}\in H^0(\mathcal{E}_T\star\mathcal{E}_C[P]\otimes \mathcal{K})$.
\end{defn}

In general one cannot expect harmonic metrics to always be compatible with $TC$. However the condition that $P$ is a parallelisable partial slice implies that the image of $\varphi\wedge \varphi^{*H}$ lies in the Lie subalgebra $\langle h\rangle\oplus \mathfrak{c}$ for every $H$ compatible with $TC$. Hence one can look for solutions $H_T\times H_C$ compatible with $TC$. 

The following proposition motivates the definition of parallelisable $\mathfrak{sl}_2$-triples. 

\begin{prop}\label{prop:SolutionHitchin equation compacible0}Let $X$ be a closed hyperbolic surface and $(\mathcal{E}_{TC},\mathbf{e}+\mathbf{r})$ be a polystable Higgs bundle in the corresponding restricted Slodowy slice, where the triple $\lbrace e,h,f\rbrace$ is parallelisable. There exists a harmonic metric $H$ on this bundle is compatible with $TC$.

When it is stable, the harmonic metric is unique.
\end{prop}

Proposition~\ref{prop:SolutionHitchin equation compacible0}
is included in Corollary~\ref{cor:SolutionHitchin equation compacible},
which we will prove later in the paper.

\subsection{Examples of Slodowy slices of Higgs bundles}
We discuss in the classical cases how the Lie algebra study from Section \ref{sec:ListLie Algebra} translates into Higgs bundles.

\subsubsection*{The case \texorpdfstring{$A_n$}{An} (\texorpdfstring{$\SL_{n+1}(\C)$}{SL(n+1,C)}), root \texorpdfstring{$\mathfrak{sl}_2$}{sl2}:}
The centralizer of the triple is the group $C$ of matrices of the following form, with blocks of size $(1,n-1,1)$:
$$\begin{pmatrix}
\lambda & 0 & 0 \\
0 & M & 0 \\
0 & 0 & \lambda  \\
\end{pmatrix},\qquad \lambda\in\C^*,\ M\in\GL(n-1,\C),\ \lambda^2\det(M)=1 .$$
The projection $(\lambda,M)\mapsto M$ realizes $C$ as a double cover of $\GL(n-1,\C)$, and its Lie algebra is $\mathfrak{gl}_{n-1}(\C)$.

Note that the one-parameter subgroup $T$ generated by $h$ acts on $\C^{n+1}$ as follows:
$$s\mapsto\begin{pmatrix}
s & 0 & 0 \\
0 & \Id & 0 \\
0 & 0 & s^{-1}  \\
\end{pmatrix},\quad s\in\C^* .$$

\begin{rem}
   In particular $T$ acts with weight $s^{-2}$ on $e$, so that a $T$-bundle $\mathcal{E}_T$ with $\mathcal{E}_T[\langle e\rangle]\simeq\mathcal{K}^{-1}$ is the frame bundle of a square root $\mathcal{K}^{1/2}$ of $\mathcal{K}$.
\end{rem}

Hence to construct the slice we fix:
\begin{itemize}
    \item a holomorphic rank $n-1$ vector bundle $\mathcal{U}$,
    \item a line bundle $\mathcal{L}$ such that $\mathcal{L}^{2}\otimes\det(\mathcal{U})\simeq\mathcal{O}$,
\end{itemize}
which together define a $C$-principal bundle $\mathcal{E}_C$.

The vector bundle associated to $\mathcal{E}_T\star \mathcal{E}_C$ and $\C^{n+1}$ writes as $\mathcal{L}\mathcal{K}^{1/2}\oplus \mathcal{U}\oplus \mathcal{L}\mathcal{K}^{-1/2}$.
The sections to construct the Higgs bundle we pick are:
\begin{itemize}
    \item $\beta \in H^0(\mathcal{U}\mathcal{L}^{-1}\mathcal{K}^{3/2})$,
    \item $\gamma \in H^0(\mathcal{U}^*\mathcal{L}\mathcal{K}^{3/2})$.
\end{itemize}
\begin{equation}
\begin{tikzcd}
\mathcal{L}\mathcal{K}^{-1/2} \arrow[r, "\beta"] & \mathcal{U} \arrow[r, "\gamma"] &\mathcal{L}\mathcal{K}^{1/2}\arrow[bend left,ll, "1"]
\end{tikzcd}.
\end{equation}

The bundles considered in \cite{BD25} are particular examples of such bundles in the case $n=2$, with $\beta=\gamma$ and $\mathcal{L}=\mathcal{O}$.

\subsubsection*{The case  \texorpdfstring{$B_n$}{Bn} ( \texorpdfstring{$\SO_{2n+1}(\C)$}{SO(2n+1,C)}), long root):} Let $n\geq 4$. The centralizer of the triple is $C=\SL(2,\C)\times \SO(2n-3,\C)$, which acts on $\C^{2n+1}$ as follows, where the matrix has blocks of size $(2,2n-3,2)$:
$$(A,M)\mapsto\begin{pmatrix}
A & 0 & 0 \\
0 & M & 0 \\
0 & 0 & (A^t)^{-1} \\
\end{pmatrix} .$$

Note that $T$ acts on $\C^{2n+1}$ as follows:
$$s\mapsto\begin{pmatrix}
s\,\Id & 0 & 0 \\
0 & \Id & 0 \\
0 & 0 & s^{-1}\Id  \\
\end{pmatrix},\quad s\in\C^* .$$

Hence to construct the slice we fix:
\begin{itemize}
    \item A $\SO_{2n-3}(\C)$-principal bundle, hence a holomorphic rank $2n-3$ vector bundle $\mathcal{U}$ with a holomorphic trivialization of $q_\mathcal{U}=S^2\mathcal{U}$,
    \item A $\SL(2,\C)$-principal bundle, \ie a rank $2$ bundle $\mathcal{V}$ with trivial determinant.
\end{itemize}

The bundle $\mathcal{E}_T\star \mathcal{E}_C$ writes as $\mathcal{V}^*\mathcal{K}^{-1/2}\oplus \mathcal{U}\oplus \mathcal{V}\mathcal{K}^{1/2}$.
The sections to construct the Higgs bundle we pick are:
\begin{itemize}
    \item $\beta \in H^0(\mathcal{V}\mathcal{U}\mathcal{K}^{3/2})$.
\end{itemize}
\begin{equation}
\begin{tikzcd}
\mathcal{V}^*\mathcal{K}^{-1/2} \arrow[r, "\beta"] & \mathcal{U} \arrow[r, "-\beta^t"] &\mathcal{V}\mathcal{K}^{1/2}\arrow[bend left,ll, "J"]
\end{tikzcd}.
\end{equation}
The symmetric bilinear pairing on this bundle  has the following block matrix:
$$\begin{pmatrix}
0 & 0 & I \\
0 & I & 0 \\
I & 0 & 0  \\
\end{pmatrix} .$$

Here $J$ is induced by the identification $\mathcal{V}\simeq \mathcal{V}^*$ induced by the symplectic form on $\mathcal{V}$, preserved by $\SL(2,\C)$. In other words $J$ is the following matrix:

$$\begin{pmatrix}
0 & 1 \\
-1 & 0 
\end{pmatrix} .$$

\begin{rem}
    Note that it makes sense to write $J$ in this Higgs bundle for the following reason: for every element $A\in \SL(2,\C)$, $A^tJA=J$. In other words, the map $J$ is the composition of the symplectic form of $\mathcal{V}$ and the bilinear pairing of the whole bundle.
\end{rem}

One can check that for this pairing, the part $J$ of the Higgs field is anti-adjoint and hence is in the Lie algebra $\mathfrak{so}(2n+1)$.
Here $\beta^t$ is obtained as the adjoint for the bilinear form to $\beta$, so that $\beta-\beta^t$ is anti-adjoint and hence in the Lie algebra $\mathfrak{so}(2n+1)$. 
In matrix form, $\beta^t$ is just the usual transpose of $\beta$.

\medskip
The bundles considered by Filip \cite{Fil22} belong to this restricted minimal Slodowy slice, in the case when $n=2$, and $\mathcal{V}$ splits as a sum of two lines bundles. The bundles considered by Zhang \cite{Zha25} do not belong to the restricted Slodowy slice, as the Higgs field contains a component in the centralized of the $\mathfrak{sl}_2$-triple.

\subsubsection*{The case \texorpdfstring{$B_n$}{Bn} (\texorpdfstring{$\SO_{2n+1}(\C)$}{SO(2n+1,C)}), short roots):} 

These are the bundles associated to maximal representations in $\SO(2,2n-1)$. 

\subsubsection*{The case  \texorpdfstring{$C_n$}{Cn} ( \texorpdfstring{$\Sp_{2n}(\C)$}{Sp(2n,C)}), long root:} 

The centralizer is $\Sp(2n-2,\C)\times \lbrace \pm \rbrace$ To construct the slice we fix:
\begin{itemize}
    \item A $\Sp_{2n-2}(\C)$-principal bundle, hence a holomorphic rank $2n-2$ vector bundle $\mathcal{W}$ with a holomorphic trivialization of $\omega_\mathcal{W}=\wedge^2\mathcal{W}$
    \item A square root $\mathcal{I}$ of the trivial line bundle $\mathcal{O}$.
\end{itemize}

The bundle $\mathcal{E}_T\star \mathcal{E}_C$ writes as $\mathcal{I}\mathcal{K}^{-1/2}\oplus \mathcal{W}\oplus \mathcal{I}\mathcal{K}^{1/2}$.
The sections to construct the Higgs bundle we pick are:
\begin{itemize}
    \item $\beta \in H^0(\mathcal{W}\mathcal{I}\mathcal{K}^{3/2})$.
\end{itemize}
\begin{equation}
\begin{tikzcd}
\mathcal{I}\mathcal{K}^{-1/2} \arrow[r, "\beta"] & \mathcal{W} \arrow[r, "-\beta^\perp"] &\mathcal{I}\mathcal{K}^{1/2}\arrow[bend left,ll, "1"]
\end{tikzcd}.
\end{equation}

\subsubsection*{The case \texorpdfstring{$C_n$}{Cn} (\texorpdfstring{$\Sp_{2n}(\C)$}{Sp(2n,C)}), short root:} 

Actually this case is not parallelisable, so we are only partially interested.

\subsubsection*{The case \texorpdfstring{$D_n$}{Dn} (\texorpdfstring{$\SO_{2n}(\C)$}{SO(2n,C)}):} 

Very similar to $B_n$, to construct the slice we fix:
\begin{itemize}
    \item A $\SO_{2n-4}(\C)$-principal bundle, hence a holomorphic rank $2n-4$ vector bundle $\mathcal{U}$ with a holomorphic trivialization of $q_\mathcal{U}=S^2\mathcal{U}$,
    \item A $\SL(2,\C)$-principal bundle, \ie a rank $2$ bundle $\mathcal{V}$ with trivial determinant.
\end{itemize}

The bundle $\mathcal{E}_T\star \mathcal{E}_C$ writes as $\mathcal{V}^*\mathcal{K}^{-1/2}\oplus \mathcal{U}\oplus \mathcal{V}\mathcal{K}^{1/2}$.
The sections to construct the Higgs bundle we pick are:
\begin{itemize}
    \item $\beta \in H^0(\mathcal{V}\mathcal{U}\mathcal{K}^{3/2})$.
\end{itemize}
\begin{equation}
\begin{tikzcd}
\mathcal{V}^*\mathcal{K}^{-1/2} \arrow[r, "\beta"] & \mathcal{U} \arrow[r, "-\beta^t"] &\mathcal{V}\mathcal{K}^{1/2}\arrow[bend left,ll, "J"]
\end{tikzcd}.
\end{equation}

\subsection{Parallel distribution of planes associated to an  \texorpdfstring{$\mathfrak{sl}_2$}{sl2}-triple.}
\label{subsec:parralel distribution of planes}

Let $\mathbb{X}=G/K$ be the symmetric space of non-compact type associated with the semi-simple Lie group $G$. Fix a parallelisable harmonic bundle, meaning a  $G$-Higgs bundle $(\mathcal{E}_T\star \mathcal{E}_C,\mathbf{e}+\mathbf{r})$ over $X$ in the Slodowy slice together with a harmonic metric $H$ compatible with the $TC$-reduction defines a corresponding harmonic map $\harm:\widetilde{X}\to \mathbb{X}$, up to the action of $G$. The additional structure on the Higgs bundle allows us a \emph{parallel distribution of planes}.

We begin by assuming that $X$ is simply connected, which we can do up to lifting the harmonic Higgs bundle to its universal cover.

Let $\nabla=D_H+\varphi-\theta_H(\varphi)$ be the flat connection constructed by the harmonic metric. This allows to trivialize the principal bundle $\mathcal{E}=\mathcal{E}_T\star \mathcal{E}_C$ as $X \times G$. This identification is well defined up to acting by an element $g\in G$. The harmonic metric then defines a map $\harm:X \to \mathbb{X}=G/K$ that is harmonic. The section $\varphi$ defines a section of $T^*X\otimes\mathfrak{g}$, so that the Maurer--Cartan form identifies $\varphi-\theta_H(\varphi):TX \to \mathfrak{g}$ with $\mathrm{d}\harm:TX\to T\mathbb{X}$. 

\medskip

The Higgs field decomposes as $\varphi=\mathbf{e}+\mathbf{r}$. We define $\mathscr{V}:TX\to T\mathbb{X}$ as the map corresponding via the Maurer--Cartan form to $e-\theta_H(e)$. Since $H$ is a metric compatible with $TC$, one can build a $\mathfrak{sl}_2$-triple $\theta_H(e)$.

\begin{defn}
    The \emph{projected} metric is the Riemannian metric $g_\pi$ on $X$ that is the pullback by $\mathscr{V}$ of the Riemannian metric on $T\mathbb{X}$.
    \end{defn}
    
\begin{rem}
The projected metric depends on the metric fixed on the symmetric space $\mathbb{X}$.
\end{rem}

\begin{prop}
\label{prop:Parallelism}
    For all $x\in X$, the exponential $\mathcal{P}_x$ of $\mathscr{V}(T_xX)$ is a totally geodesic copy of $\mathbb{H}^2$ inside $\mathbb{X}$ associated to the $\mathfrak{sl}_2$-triple $\lbrace e,h,f\rbrace$. 

This distribution of planes along $\harm:X\to \mathbb{X}$ is \emph{parallel} in the sense that for every local section $v$ of $TX$, if $\nabla^\mathbb{X}$ is the Levi-Civita connection on $T\mathbb{X}$ and  $\nabla^{g_\pi}$ the Levi-Civita connection of the projected metric $g_\pi$ on $X$:
$$\nabla^\mathbb{X} \mathscr{V}(v)=\mathscr{V}(\nabla^{g_\pi} v).$$
\end{prop}

\begin{proof}
Since the harmonic metric $H$ is compatible with the $TC$-reduction, one can assume that at a point $x\in X$, $\theta_H(e)=-f$ up to identifying the fiber at $x$ with $\mathfrak{g}$ in a way that is compatible with the $TC$ reduction of the bundle. Hence for all $v\in T_xX$, $\mathscr{V}(v)$ is identified with the plane spanned by $e+f$ and $ie-if$, which is tangent to the totally geodesic $\mathbb{H}^2$ associated with the copy of $\mathfrak{sl}_2(\R)$ spanned by $e+f$, $ie-if$ and $ih$, which is conjugated to the one spanned by $e$, $h$ and $f$.

\medskip

For the parallelism property, note that the Levi-Civita connection on $T\mathbb{X}$ corresponds via the Maurer--Cartan form to the Chern connection $D_H$. We therefore compute the quantity corresponding to $\nabla^\mathbb{X} \mathscr{V}(v)$ via the Maurer--Cartan form which is :
$$D_H (\mathbf{e}-\theta_H(\mathbf{e}))= D_H (\mathbf{e})-\theta_H(D_H (\mathbf{e})).$$

Since $H$ is compatible with $TC$, which preserves the complex line $\langle e \rangle\subset \mathfrak{g}$, the Chern connection $D_H$ restricts to the Chern connection of the induced metric on $\langle e\rangle$. This holomorphic line bundle over $X$ is naturally identified via the section $\mathbf{e}$ with $TX$, and the restriction of $H$ restricts to the projected metric. Hence one obtains that $D_H (\mathbf{e})(v)=\mathbf{e}(\nabla^{g_\pi}v)$. This concludes the proof.

\end{proof}

The most important consequence of the parallelism of this distribution is the following. Recall that to every oriented curve in $\mathbb{X}$ one can associate a unique isometry $G$ of $\mathbb{X}$ realizing the parallel transport along this curve.

\begin{cor}
\label{cor:ParallelTransport Distribution planes}
    Let $\gamma:\R\to X$ be a geodesic in $X$ for the projected metric $g_\pi$. Let $g_t\in G$ be the realization of the parallel transport between $\gamma(0)$ and $\gamma(t)$ along $\gamma$, meaning that $g'_t$ corresponds to $\varphi-\theta_H(\varphi)$. Then $\mathscr{V}(\gamma'(t))=g_t\cdot \mathscr{V}(\gamma'(0))$ for all $t>0$.
\end{cor}

\subsection{ \texorpdfstring{$3$}{3}-cyclic Higgs bundles}
\label{subsec:3-cyclic}
Let $\mathcal W$ be a holomorphic rank $n$ vector bundle, $\mathcal V$ be a holomorphic rank $m$ vector bundle. Let $\alpha\in H^0(X,\Hom(\mathcal{W},\mathcal{V})\otimes\mathcal{K}), \beta\in H^0(X,\Hom(\mathcal V,\mathcal{WK})\otimes\mathcal{K}).$ Let $(\mathcal E, \varphi)$ be the following Higgs bundle: 
\begin{equation}
\begin{tikzcd}
\mathcal{WK} \arrow[rr, "I_n"] && \mathcal{W} \arrow[ld, "\alpha"] \\&\mathcal{V}\arrow[lu, "\beta"]
\end{tikzcd}.
\end{equation}

\begin{assumption}\label{ass:3cyclic}

\begin{enumerate}
\item The harmonic metric $H=H|_{\mathcal{WK}}\oplus H|_{\mathcal W}\oplus H|_{\mathcal V}.$ And there exists a Hermitian metric $g_\pi$ on $\mathcal K^{-1}$ such that: \[H|_{\mathcal{WK}}= H|_{\mathcal W}\otimes g_\pi^{-1}.\]
\item The metric $g_\pi$ is a complete metric on $X$. 
\end{enumerate}

We call the metric $g_\pi$ the \emph{projected metric} of the $3$-cyclic bundle.
\end{assumption}

Such $3$-cyclic bundles arise naturally when considering harmonic bundles. Let $G$ be a simple complex group whose type is not $E_8$ with finite center. Consider the linear representation $V$ of $\mathfrak{g}$ from Proposition \ref{prop:RepresentationV}, listed in Figure \ref{table:Slodowy slices}. Let us assume that this representation defines a representation $G\to \SL(V)$. 

\begin{rem}
\label{rem:Finite Cover Representation}
If the Lie algebra representation does not extend to $G$, there is a finite central extension $G_0$ of $G$ such that the Lie algebra representation extends to $G_0\to \SL(V)$ and one can lift any $G$-harmonic bundle on $X$ to a $G_0$-harmonic bundle on $\hat{X}$, up to eventually taking a finite cover $\hat{X}$ of $X$. Our main results will then still hold.
\end{rem}

The following proposition says that bundles in the restricted Slodowy slices are $3$-cyclic, when viewed through the right representation.

\begin{prop}
\label{prop:Restricted Slodowy Implies 3 cyclic}
Let $G$ be a simple complex group whose type is not $E_8$ with finite center with the linear representation $G\to\SL(V)$ described above. Let $(\mathcal{E},\varphi,H)$ be a parallelisable harmonic $G$-Higgs bundle in the restricted Slodowy slice whose projected metric is complete on $X$.

The associated $\SL(V)$-Higgs bundle is a $3$-cyclic bundle satisfying Assumption \ref{ass:3cyclic}, for which $n,m=\dim(V_1),\dim(V_0)$ where $V=V_{-1}\oplus V_0\oplus V_1$ described in Proposition \ref{prop:RepresentationV}. Moreover up to a positive scalar the projected metric $g_\pi$ on $X$ for the bundle $(\mathcal{E},\varphi,H)$ in the restricted Slodowy slice coincides with the projected metric of the $3$-cyclic bundle.
\end{prop}

\begin{proof}
Since the structure of the $G$-principal bundle $\mathcal{E}$ reduces to $TC$, which preserves the decomposition $V=V_{-1}\oplus V_0\oplus V_1 $, the corresponding vector bundle $\mathcal{E}[V]$ decomposes into $\mathcal{E}[V_{-1}]=\mathcal{W}$, $\mathcal{E}[V_0]=\mathcal{V}$, $\mathcal{E}[V_{1}]=\mathcal{W}'$. Applying Proposition \ref{prop:RepresentationV}, if we write $\varphi=\mathbf{e}+\mathbf{r}$ the never-vanishing section $\mathbf{e}$ defines an isomorphism $\mathcal{W}'\simeq \mathcal{W}\mathcal{K}$ and the section $\mathbf{r}\in H^0(\mathcal{E}[R]\otimes \mathcal{K})$ defines sections $\alpha:\mathcal{W}\to \mathcal{V}\mathcal{K}$ and $\beta:\mathcal{V}\to \mathcal{W}\mathcal{K}^2$. 

\medskip

The reduction $H$ is compatible with the $TC$-structure, which preserves the block decomposition $V=V_1\oplus V_0\oplus V_{-1}$. This implies that the corresponding metric on $\mathcal{W}\mathcal{K}\oplus \mathcal{V}\oplus \mathcal{W}$ is block diagonal. The Cartan involution $\theta$ on $\mathcal{E}[\mathfrak g]$ corresponding to the reduction $H$ satisfies $[\mathbf{e},\theta(\mathbf{e})]\in \langle h\rangle$. This implies that $[I_n,I_n^{*H}]$ acts by a scalar multiple on $\mathcal{W}$ and $\mathcal{W}\mathcal{K}$. This implies that $I_n$ is conformal for the metrics $H_{\mid \mathcal{W}\mathcal{K}}$ and $H_{\mid\mathcal{W}}$, which implies the statement of Assumption \ref{ass:3cyclic}. 

The projected metric for the harmonic bundle in the Slodowy slice for $v\in TX$ is $\lVert \mathbf{e}(v)\rVert$ where the norm is taken with respect to the harmonic metric, which is equal up to a fixed constant to $\lVert I_n(v)\rVert$. In particular the completeness of both metrics are equivalent.
\end{proof}

The following is the other example of $3$-cyclic Higgs bundles satisfying Assumption \ref{ass:3cyclic} that we will consider. Tehse bundles belong to the restricted minimal Slodowy slice for $\SO(n,\C)$.

\begin{exmp}
Let $(\mathcal V,Q_{\mathcal V}), (\mathcal W,Q_{\mathcal W})$ be two orthogonal bundles. Let $\mathcal L$ be a holomorphic line bundle. Let $\mu\in H^0(X,\mathcal L^{-1}\mathcal {WK})$, $\gamma\in H^0(X,\mathcal V\mathcal L\mathcal K^2).$ Let $Q_0$ be the natural orthogonal structure on $\mathcal {LK}\oplus \mathcal L\oplus \mathcal L^{-1}\oplus \mathcal L^{-1}\mathcal K^{-1}$. Let $\gamma^{\dagger}, \mu^{\dagger}$ be the adjoint of $\gamma, \mu$ with respect to $Q_0$ and $Q_{\mathcal V}, Q_{\mathcal W}$ respectively. Consider the following Higgs bundle:
\begin{equation}
\begin{tikzcd}
&&\mathcal{V}\arrow[dll, "\gamma^{\dagger}"']&&\\\mathcal{KL}\arrow[r, "1"'] & \mathcal{L}\arrow[r, "\mu"'] &\mathcal{W}\arrow[r, "\mu^{\dagger}"']&\mathcal{L}^{-1}\arrow[r, "1"']&\mathcal{L}^{-1}\mathcal K^{-1}\arrow[llu, "\gamma"']&
\end{tikzcd}.
\end{equation}
Assume the harmonic metric $H=H_{\mathcal V}\oplus H_1g_\pi^{-1}\oplus H_1\oplus H_{\mathcal W}\oplus H_1^{-1}\oplus H_1^{-1}g_\pi,$ where $H_{\mathcal V}, H_{\mathcal W}$ are Hermitian metrics compatible with $Q_{\mathcal V}, Q_{\mathcal W}$ respectively, $H_1$ is a Hermitian metric on $\mathcal L$, and $g_\pi$ is a Hermitian metric on $\mathcal K^{-1}$. It is an $\SO(\mathrm{rk}(W)+2,\mathrm{rk}(V)+2)$-harmonic bundle.

Assume the metric $g_\pi$ is complete.

When $\mathcal W=\mathcal O, (\mathcal V,Q_{\mathcal V})=(\mathcal N\oplus \mathcal N^{-1},\begin{pmatrix}0&1\\1&0\end{pmatrix}),$ it is an $\SO(3,4)$-harmonic bundle. Let $\gamma=(\gamma_1,\gamma_2)$ with respect to the decomposition $\mathcal V=\mathcal N\oplus \mathcal N^{-1}$. If in addition $\gamma_1=\mu$, it is a $G_2'$-harmonic bundle (a $\alpha$-bundle according to the terminology of \cite{DE26}). 

\medskip

These bundles are $3$-cyclic if we group together $\mathcal{KL}\oplus \mathcal{L}^{-1}$, $ \mathcal{L}^{-1}\mathcal{K}^{-1}\oplus \mathcal{L}$ and $\mathcal{V}\oplus \mathcal{W}$.
\end{exmp}

\section{Analysis for Slodowy slices}

\label{sec:Analysis}
In this section we consider analytic aspects of parallel distributions of planes. We consider parallelisable harmonic Higgs bundles over a (possibly non-compact) Riemann surface $X$ in $3$ families:
\begin{itemize}
    \item[-] the minimal restricted Slodowy slice, our main focus,
    \item[-] maximal $\SO(2,n+1)$-bundles,
    \item[-] a family of nilpotent $\SL(n+4,\R)$-Higgs bundles.
\end{itemize}

We first use a maximum principle to control the harmonic metric in each case, which geometrically allows to bound the angle between the tangent planes of the minimal surface in the symmetric space and the planes in the parallel distribution.

We then construct for each family filled hyperboloids and cones in projective space associated to each tangent vector to the surface. We then use the bounds given by the maximum principle to show that these cones and hyperboloids are nested, and estimate their nestedness.

Finally we discuss the differential geometric properties of the map $Y:\widetilde{X}\to G/(TC\cap K)$ defining the parallel distribution of planes, and discuss the properness of the associated map in $G/G^\sigma$.

\subsection{Maximum principle}
\subsubsection{Preliminaries}
Consider a Hermitian metric $g$ on $\mathcal K^{-1}$ and write it $g=\tilde gdz\otimes d\bar z$ on $X$. Let $\omega=\sqrt{-1}\tilde gdz\wedge d\bar z$ be the associated $2$-form. The associated Riemannian metric is $\Re(g)=\frac{1}{2}(g+\bar g)$. Let $\Lambda_g$ be the contraction operator on $2$-forms such that for any $2$
-form, $\alpha=\Lambda_g(\alpha)\cdot \omega$. The Laplacian with respect to $g$, the Gaussian curvature $K_g$ of $\Re(g)$
 are 
 \[\Delta_g=\frac{1}{\tilde g}\partial_z\partial_{\bar z}, \quad K_g=-\frac{2}{\tilde g}\partial_z\partial_{\bar z}\log \tilde g. \]
The curvature of the Chern connection on $(\mathcal K^{-1}, g)$ is $F(g)=\bar\partial\partial\log\tilde g$. Thus \[K_g=2i\Lambda_gF(g).\]
For a section $s$ of $\End(E)$, denote the square norm with respect to $H$, $|s|_H^2=\tr(ss^*).$ For a section $\sigma$ of $\End(E)\otimes K^{\otimes p}\otimes \bar K^{\otimes q}$, locally $\sigma=sdz^{\otimes p}\otimes d\bar z^{\otimes q}$, denote by $\|\sigma\|^2$ the square norm with respect to $H$ and $g$, explicitly, \[\|\sigma\|^2=\tr(ss^*)\cdot \tilde g^{-p-q}.\]
Denote by $\|\sigma\|_{g}^2$ if we want to specify the role of $g$. For $(p,q)=(1,0)$, denote the $(1,1)$-tensor $|\sigma|_H^2=\tr(ss^*)dz\otimes d\bar z.$

We first recall the following formula, which is a slight modification of \cite[Lemma 2.4]{Moc16}, where it was called Simpson's main estimate.
\begin{lem}(\cite[Lemma 3.3]{DaiLiBounded})\label{lem:DaiLiSection}Let $(\mathcal E, \varphi, H)$ be a harmonic bundle. Locally, $\varphi=fdz$ for a local holomorphic section f of $\End(\mathcal E)$. 
For a holomorphic section $s$ of $\End(\mathcal E)$, we have
\[\partial_z\partial_{\bar z}\log |s|_H^2\geq \frac{|[s, f^{*_H}]|_H^2-|[s,f]|_H^2}{|s|_H^2}.\]
\end{lem}
As an immediate corollary, we obtain
\begin{lem}\label{lem:EquationSection}
For a holomorphic section $s$ of $\End(\mathcal E)\otimes \mathcal K$, we have 
\[\Delta_g\log \|s\|^2\geq \frac{\|[s, \varphi^{*_H}]\|^2-\|[s,\varphi]\|^2}{\|s\|^2}+\frac{1}{2}K_g.\]
\end{lem}

\begin{lem}(\cite[Lemma 2.5]{DaiLiBounded})\label{lem:NegativeCurvature}
    Suppose $(X,g)$ is a complete surface with curvature $K_g$ bounded from below. Suppose there exists positive constants $c$ and $d$ such that 
    \[\Delta_gK_g\geq cK_g(K_g+d).\] Then $K_g<0$ or $K_g\equiv 0$.
Moreover, if $X$ is hyperbolic, let $g_{hyp}$ be the unique complete conformal hyperbolic metric. If $g\leq Cg_{hyp}$ for some positive constant $C$, then there exists a constant $\delta>0$ only depending on $C$ such that $K_g\leq -\delta.$
\end{lem}
\begin{rem}
In Lemma \ref{lem:NegativeCurvature}, the dependence on constant $C$ is implicit in \cite[Section 2]{DaiLiBounded}, although not explicitly stated.
\end{rem}
\begin{lem}(\cite[Proposition 3.14]{LiMochizukiComplete})\label{lem:BoundedSpectralImpliesBoundedNorm}
Let $(\mathcal E, \varphi, H)$ be a harmonic bundle. Assume that $\sup_{1\leq k\leq n, X}{|\tr(\varphi^k)|_{g_{hyp}}}$ is bounded. Then there exists a constant $C$ only depending on $\sup_{1\leq k\leq n, X}{|\tr(\varphi^k)|_{g_{hyp}}}$, such that $\|\varphi\|_{g_{hyp}}^2\leq C$.
\end{lem}

As a corollary of Lemma \ref{lem:NegativeCurvature} and Lemma \ref{lem:BoundedSpectralImpliesBoundedNorm}, we obtain the following:
\begin{lem}\label{lem:NegativeCurvature1}
    Suppose $(X,g)$ is a complete surface with curvature $K_g$ bounded from below. Suppose there exists positive constants $c$ and $d$ such that 
    \[\Delta_gK_g\geq cK_g(K_g+d).\] Then $K_g<0$ or $K_g\equiv 0$.
Moreover, the following are equivalent:
\begin{itemize}
    \item $X$ is equipped with a  complete conformal hyperbolic metric $g_{hyp}$\\ and $\sup_{1\leq k\leq n, X}{|\tr(\varphi^k)|_{g_{hyp}}}$ is bounded;
    \item There exists a constant $\epsilon>0$ such that $K_g\leq -\epsilon.$
\end{itemize}
Furthermore, the constants $\sup_{X,  1\leq k\leq n}|\tr(\varphi^{k})|_{g_{hyp}}$ and $\epsilon$ depend only on each other. 
\end{lem}

\subsection{Maximum principle for \texorpdfstring{$3$}{3}-cyclic bundles}

We now consider $3$-cyclic Higgs bundles, described in Section \ref{subsec:3-cyclic}. Recall that they have the following form:
\begin{equation}
\begin{tikzcd}
\mathcal{WK} \arrow[rr, "I_n"] && \mathcal{W} \arrow[ld, "\alpha"] \\&\mathcal{V}\arrow[lu, "\beta"]
\end{tikzcd}.
\end{equation}

Assumption \ref{ass:3cyclic} requires that the harmonic metric decomposes as $H=H|_{\mathcal{WK}}\oplus H|_{\mathcal W}\oplus H|_{\mathcal V}$ and there exists a complete metric $g_\pi$ on $\mathcal K^{-1}$ called the \emph{projected metric} such that: \[H|_{\mathcal{WK}}= H|_{\mathcal W}\otimes g_\pi^{-1}.\]

\begin{rem}
\label{rem:ProjectedMetricNormalization}
What we previously defined the projected metric depended on a choice of a normalization of the metric on the symmetric space. There is such a choice such that the projected metric is exactly this one.
\end{rem}

We obtain the following maximum principle:

\begin{thm}\label{thm:NormBound3-cyclic}
Under Assumption \ref{ass:3cyclic}, we have:
\begin{enumerate}
    \item $i\Lambda_{g_\pi}(\beta\wedge\beta^*-\alpha^*\wedge\alpha)=(2+\frac{1}{2}K_{g_\pi})\cdot I_n$,
    \item $K_{g_\pi}=\frac{2}{n}(\|\beta\|^2+\|\alpha\|^2)-4\leq 0$. If the equality holds at one point, it holds everywhere.
    \item if $X$ is hyperbolic, $g_\pi\geq g_\text{hyp}$ where $g_\text{hyp}$ is the hyperbolic metric of curvature $-4$. The inequality is strict everywhere if $(\alpha,\beta)\neq (0,0)$.
\end{enumerate}
Moreover, the following are equivalent:
\begin{itemize}
    \item $X$ is equipped with a  complete conformal hyperbolic metric $g_{hyp}$\\ and $\sup_{X,  3k\leq 2n+m}|\tr(\varphi^{3k})|_{g_{hyp}}$ is bounded;
    \item There exists a constant $\epsilon>0$ such that $K_{g_\pi}\leq -\epsilon.$
\end{itemize}
Furthermore, the constants $\sup_{X,  3k\leq 2n+m}|\tr(\varphi^{3k})|_{g_{hyp}}$ and $\epsilon$ depend only on each other. 

If in addition, $\beta\alpha=0$, then $K_{g_\pi}\leq -1.$

If $\beta=0$ or $\alpha=0$, then $K_{g_\pi}\leq -2.$
\end{thm}

\begin{proof}
The Hitchin equation gives :
\[-i\Lambda_{g_\pi} (F_H|_{\mathcal{WK}})=i\Lambda_{g_\pi}(\beta\wedge\beta^*+I_n^*\wedge I_n)=i\Lambda_{g_\pi}(\beta\wedge\beta^*)- I_n\]
\[-i\Lambda_{g_\pi}(F_H|_{\mathcal W})=i\Lambda_{g_\pi}(I_n\wedge I_n^*+\alpha^*\wedge\alpha)= I_n+i\Lambda_{g_\pi}(\alpha^*\wedge\alpha).\]
Thus :
\[i\Lambda_{g_\pi} F_H|_{\mathcal{WK}}-i\Lambda_{g_\pi} F_H|_{\mathcal W}=i\Lambda_{g_\pi}(\beta\wedge\beta^*-\alpha^*\wedge\alpha)-2I_n.\]
By Assumption \ref{ass:3cyclic}, we have $-i\Lambda_{g_\pi} F_H|_{\mathcal{WK}}+i\Lambda_{g_\pi} F_H|_{\mathcal W}=i \Lambda_{g_\pi} F(g_\pi)\cdot I_n.$ 
We thus have part 1. The curvature 
\[K_{g_\pi}=2i\Lambda_{g_\pi}F(g_\pi)=\frac{2}{n}\tr(\Lambda_{g_\pi}(\beta\wedge\beta^*-\alpha^*\wedge\alpha)-2I_n)=\frac{2}{n}(\|\beta\|^2+\|\alpha\|^2)-4.\]
Let $\delta$ be the section of $\End(\mathcal E)\otimes \mathcal K$ given by $\alpha$ and $\beta$. Thus $\|\delta\|^2=\|\alpha\|^2+\|\beta\|^2.$ From Lemma \ref{lem:EquationSection}, we have
\begin{eqnarray*}
\Delta_{g_\pi}\log\|\delta\|^2&\geq&\frac{\|[\delta,\varphi^*]\|^2-\|[\delta,\varphi]\|^2}{\|\delta\|^2}+\frac{1}{2}K_{g_\pi}\\&\geq&\frac{2\|\alpha^*\alpha\|^2+2\|\beta^*\beta\|^2-2\|\beta\alpha\|^2}{\|\alpha\|^2+\|\beta\|^2}-1+\frac{1}{2}K_{g_\pi}\\
&\geq&\frac{\|\alpha\alpha^*\|^2+\|\beta^*\beta\|^2}{\|\alpha\|^2+\|\beta\|^2}-1+\frac{1}{2}K_{g_\pi}\\
&\geq&\frac{1}{2n}(\|\alpha\|^2+\|\beta\|^2)-1+\frac{1}{2}K_{g_\pi}\\
&=&\frac{3}{2n}(\|\alpha\|^2+\|\beta\|^2)-3.
\end{eqnarray*}
In the third inequality, we make use of the Cauchy-Schwarz inequality $\|AB\|^2=\tr(B^*A^*AB)=\tr(A^*ABB^*)\leq \|A^*A\|\cdot\|BB^*\|$ for two $n\times n$ matrices $A$ and $B$. In the fourth inequality, we make use of the fact that $|AA^*|^2\geq \frac{1}{n}|A|^4$ for an $n\times n$ matrix $A$.

If in addition, $\beta\alpha=0$, we have 
\begin{eqnarray*}\Delta_{g_\pi}\log(\|\alpha\|^2+\|\beta\|^2)
&\geq&\frac{2\|\alpha\alpha^*\|^2+2\|\beta^*\beta\|^2}{\|\alpha\|^2+\|\beta\|^2}-1+\frac{1}{2}K_{g_\pi}\\
&\geq&\frac{1}{n}(\|\alpha\|^2+\|\beta\|^2)-1+\frac{1}{2}K_{g_\pi}\\
&=&\frac{2}{n}(\|\alpha\|^2+\|\beta\|^2)-3.
\end{eqnarray*}

If $\beta$=0, then 
we have 
\begin{eqnarray*}\Delta_{g_\pi}\log(\|\alpha\|^2)
&\geq&\frac{2\|\alpha\alpha^*\|^2}{\|\alpha\|^2}-1+\frac{1}{2}K_{g_\pi}\geq\frac{2}{n}(\|\alpha\|^2)-1+\frac{1}{2}K_{g_\pi}\\
&=&\frac{3}{n}(\|\alpha\|^2)-3.
\end{eqnarray*}

If $\alpha=0$, then 
we have 
\begin{eqnarray*}\Delta_{g_\pi}\log(\|\beta\|^2)
&\geq&\frac{2\|\beta^*\beta\|^2}{\|\beta\|^2}-1+\frac{1}{2}K_{g_\pi}\geq\frac{2}{n}(\|\beta\|^2)-1+\frac{1}{2}K_{g_\pi}\\
&=&\frac{3}{n}(\|\beta\|^2)-3.
\end{eqnarray*}

By assumption $g_\pi$ is complete and $K_{g_\pi}\geq -4$, hence we can apply the Cheng-Yau maximum principle and obtain $\|\alpha\|^2+\|\beta\|^2\leq 2n$. Moreover, if $\beta\alpha=0$, $\|\alpha\|^2+\|\beta\|^2\leq \frac{3n}{2};$  if $\beta=0$, $\|\alpha\|^2\leq n;$ if $\alpha=0$, $\|\beta\|^2\leq n;$

\medskip

If the projected metric $g_\pi=e^{2u}g_\text{hyp}$ is written in terms of the hyperbolic metric of curvature $-4$, using the curvature formula for both $g_\pi$ and $g_{\text{hyp}}$, the function $u$ satisfies the following equation:
$$ -4\Delta_{g_{\pi}}u=K_{g_{\pi}}-e^{-2u}K_{g_{\text{hyp}}}\geq 4(e^{-2u}-1).$$
Since $g_\pi$ is complete with $K_{g_\pi}\geq -4$, once again we can apply the Cheng-Yau maximum principle and obtain $-u\leq 0$. If $u=0$ at some point, the strong maximum principle implies that $u=0$ everywhere, and hence $K_{g_\pi}=-4$ everywhere, so $(\alpha,\beta)=(0,0)$.

 \medskip

We continue to calculate $\Delta_{g_\pi}K_{g_\pi}$ in order to apply Lemma \ref{lem:NegativeCurvature}. 
\begin{eqnarray*}\Delta_{g_\pi}K_{g_\pi}&=&\frac{2}{n}\Delta_{g_\pi}(\|\alpha\|^2+\|\beta\|^2)\geq\frac{2}{n}(\|\alpha\|^2+\|\beta\|^2)\cdot \Delta_{g_\pi}\log(\|\alpha\|^2+\|\beta\|^2)\\
&\geq& \frac{2}{n}(\|\alpha\|^2+\|\beta\|^2)\cdot( \frac{3}{2n}(\|\alpha\|^2+\|\beta\|^2)-3)\\
&=&\frac{3}{4}(K_{g_\pi}+4)\cdot K_{g_\pi}
\end{eqnarray*}
By the continuity, the inequality extends to the whole surface $X$. By Lemma \ref{lem:NegativeCurvature}, $K_{g_\pi}<0$ or $K_{g_\pi}\equiv 0$. 

By the $3$-cyclic property, $\tr(\varphi^k)=0$ unless $3$ divides $k$. Applying Lemma \ref{lem:NegativeCurvature1}, we obtain the result.
\end{proof}

In particular we obtain the following result, the second part coming from Proposition \ref{prop:Restricted Slodowy Implies 3 cyclic}.

\begin{cor}
\label{cor:Maximum principle}
For a $3$-cyclic Higgs bundle satisfying Assumption \ref{ass:3cyclic} we have :
$$\lVert \alpha\rVert^2+\lVert \beta\rVert ^2\leq 2\lVert I_n\rVert^2.$$ 

In other words for a harmonic $G$-Higgs bundle $(\mathcal{E},\varphi,H)$ in the restricted Slodowy slice of the minimal $\mathfrak{sl}_2$-triple with complete projected metric,  if we write $\varphi=\mathbf{e}+\mathbf{r}$:
$$\lVert \mathbf{r}\rVert^2\leq 2\lVert \mathbf{e}\rVert^2.$$ 
\end{cor}
\subsection{Nestedness of cones}
\label{subsec:Nestedness}
Given a unit vector for the metric $g_\pi$, $v\in T^1X$, in some orthonormal basis for $H$ we obtain :
\[\Psi(v)=\varphi(v)+\varphi^{*H}(v)=\begin{pmatrix}0&I_n&b\\I_n&0&\bar a^T\\\bar b^T&a&0\end{pmatrix}.\]

In this expression $a=\alpha(v)$ and $b=\beta(v)$. Let us write the Fuchsian part of $\Psi$, which corresponds to $\mathbf{e}+\mathbf{e}^*$ for Higgs bundles in the Slodowy slice of the minimal $\mathfrak{sl}_2$-triple:
$$\Psi_0(v)=\begin{pmatrix}0&I_n&0\\I_n&0& 0\\0&0&0\end{pmatrix}.$$
We consider the following change of basis diagonalizing the Fuchsian part $\Psi_0$ of $\Psi$ the Higgs field:
\[P=\begin{pmatrix}\tfrac{1}{\sqrt{2}}I_n&\tfrac{1}{\sqrt{2}}I_n&0\\-\tfrac{1}{\sqrt{2}}I_n&\tfrac{1}{\sqrt{2}}I_n&0\\
0&0&I_m\end{pmatrix}.\]
\[\hat\Psi(v):=P\Psi(v)P^{-1}=
\begin{pmatrix}
I_n&0&\tfrac{b+\bar a^T}{\sqrt{2}}\\
0&-I_n&\tfrac{-b+\bar a^T}{\sqrt{2}}\\
\tfrac{\bar b^T+a}{\sqrt{2}}&\tfrac{-\bar b^T+a}{\sqrt{2}}&0\end{pmatrix},\quad \hat\Psi_0(v):=P\Psi_0(v)P^{-1}.\]

We work in this new basis, in which every vector $\vec{X}$ decomposes as $X=\vec{X_1}+\vec{X_2}+\vec{X_3}$. Using this basis depending on the unit vector $v\in TX$, we define the following hermitian form on $\mathcal{W}\oplus \mathcal{V}\oplus \mathcal{W}\mathcal{K}$:
\begin{equation}
\label{eq:Definition qv}
q_v(\vec X, \vec X')=\Re(\vec X_1\cdot \overline{\vec X_1'}-\vec X_2\cdot \overline{\vec X_2'}-\vec X_3\cdot \overline{\vec X_3'}).
\end{equation}

This hermitian form $q_v$ does not depend on any choice of basis, as it is characterized by the hermitian metric $H$ and the eigenspaces of the endomorphism $\Psi_0(v)$.
It will allow us to define cones $C_v\subset \mathbb{P}(\C^{2n+m})$ and filled hyperboloids $\mathcal{H}_v\subset \C^{2n+m}$ as follows:
$$C_v=\{[\vec X]|q_v(\vec X, \vec X)>0\},$$
$$\mathcal{H}_v=\{\vec X|q_v(\vec X, \vec X)\geq 1\}.$$

\begin{rem}
This construction of cones and filled hyperboloids for Slodowy slices associated to the minimal $\mathfrak{sl}_2$ can be also viewed as follows: given $v\in T^1X$ we consider $\mathscr{V}(v)\in T^1\X$. Through the linear representation $V$ of $G$, $\X$ can be viewed as a space of hermitian metrics on $V$, and $\mathscr{V}(v)$ as a self adjoint endomorphism of $V$, which is exactly $\Psi_0(v)$. The sets $C_v\subset \mathbb{P}(V)$ and $\mathcal{H}_v\subset V$ are determined by this endomorphism associated to $\mathscr{V}(v)$.
\end{rem}

These cones and these filled hyperboloids, when viewed in the same ambient space by using the flat connection associated to the harmonic bundle, are nested along the geodesic flow. This will be a consequence of the parallelism and of the following property:

\begin{prop}\label{prop:3cycliclowerbound}
Under Assumption \ref{ass:3cyclic}, for all $\vec X\in \partial C_v$ meaning $q_v(\vec X,\vec X)=0$:
\begin{equation}\label{eq:q_v Cone}
    q_v(\hat\Psi(v) \vec X,\vec X)\geq\frac{3-\sqrt{9+2K_{g_\pi}}}{4}\lVert\vec X\rVert^2\geq -\frac{1}{12}K_{g_\pi} \lVert\vec X\rVert^2.
\end{equation}

Under Assumption \ref{ass:3cyclic}, for all $\vec X\in \partial \mathcal{H}_v$ meaning $q_v(\vec X,\vec X)=1$:
\begin{equation}\label{eq:q_v Hyperboloid}
q_v(\hat\Psi(v) \vec X,\vec X) \geq  1. 
\end{equation}
Moreover, if $K_{g_\pi}<0$ equality in \eqref{eq:q_v Hyperboloid} can only occur if $\hat\Psi_0(\vec X)=\vec X$.
\end{prop}

In the terminology that we introduce in Section \ref{sec:Dynamics}, the endomorphism $\Psi(v)\in\mathfrak{sl}_{2n+m}(\C)$ is at least $\frac{1}{12}|K_{g_\pi}|$-nesting the cone $C_v$ and $1$-nesting the filled hyperboloid $\mathcal{H}_v$.

\begin{rem}
    Recall that here $K_{g_\pi}=2\frac{\lVert \mathbf{r}\rVert^2}{\lVert \mathbf{e}\rVert^2}-4$, since $\lVert \mathbf{e}\rVert^2 =n$ and $\lVert \mathbf{r}\rVert^2=\lVert \alpha\rVert^2+\lVert \beta\rVert^2$.
\end{rem}

\begin{proof}

Note that when computing $q(\hat\Psi(v)\vec X,\vec X)$ two terms cancel, as the $(1,3)$ and $(3,1)$ component of $\hat\Psi(v)$ in the basis given by $P$ form a matrix in $\mathfrak{su}(q_v)$:
\begin{eqnarray*}
(*)=q_v(\hat\Psi(v) \vec X,\vec X)=|\vec X_1|^2+|\vec X_2|^2-\sqrt{2}\Re(\overline{\vec X_1}^T\cdot (-b+\bar a^T)\cdot \vec X_3).\end{eqnarray*}
Let us first consider the case when $\vec X\in \partial C_v:$
\begin{eqnarray*}
(*)&=&(2|\vec X_2|^2+|\vec X_3|^2)-\sqrt{2}\Re(\overline{\vec X_2}^T\cdot (-b+\bar a^T)\cdot \vec X_3)\\
&\geq&(2|\vec X_2|^2+|\vec X_3|^2)-2\frac{\lVert \mathbf{r} \rVert}{\lVert \mathbf{e}\rVert} |\vec X_2||\vec X_3|.
\end{eqnarray*}

Here we used the fact that $b\bar b^T+\bar a^Ta\leq \frac{\lVert \mathbf{r} \rVert^2}{\lVert \mathbf{e}\rVert^2}\cdot I_n$.

This quantity is positive since $\frac{\lVert \mathbf{r} \rVert}{\lVert \mathbf{e}\rVert}\leq \sqrt{2}$. However one can estimate more precisely a lower bound.
A straightforward analysis yields that a global minimum of the quantity $\frac{2X^2+Y^2-\alpha XY}{2(X^2+Y^2)}$ for $X,Y \in \R$ with $(X,Y)\neq (0,0)$ given a fixed $\alpha \in[0,2\sqrt{2}]$ is obtained for $X=\sqrt{1+\alpha^2}-1$, $Y=\alpha$, and is equal to $\frac{3-\sqrt{1+\alpha^2}}{4}$. Note that $\lVert\vec X\rVert^2=2(|\vec X_2|^2+|\vec X_3|^2)$ since $[\vec X]\in \partial C_v$. We apply this inequality, to obtain that \[(*)\geq \frac{3-\sqrt{1+4\frac{\lVert \mathbf{r}\rVert^2}{\lVert \mathbf{e} \rVert^2}}}{4}\lVert\vec X\rVert^2=\frac{3-\sqrt{9+2K_{g_\pi}}}{4}\lVert\vec X\rVert^2\geq -\frac{1}{12}K_{g_\pi}\lVert\vec X\rVert^2.\]

Now we consider the case when $\vec X\in  \partial \mathcal{H}_v:$
\begin{eqnarray*}
(*)=1+(2|\vec X_2|^2+|\vec X_3|^2)-\sqrt{2}\Re(\overline{\vec X_2}^T\cdot (-b+\bar a^T)\cdot \vec X_3).
\end{eqnarray*}
This quantity is greater than $1$ as in the previous statement. If $K_{g_\pi}<0$, then we have a strict inequality except if $\vec X=\vec X_1$, which means exactly that $\hat\Psi_0(\vec X)=\vec X$.
\end{proof}

\subsection{Maximal \texorpdfstring{$\SO_o(2,n+1)$}{SOo(2,n+1)}-Higgs bundles}
\label{subsec:MaximalSO2N}
Let $\mathcal I$ be a holomorphic line bundle satisfying $\mathcal I^2=\mathcal O$, $(\mathcal V, Q_{\mathcal V})$ be an orthogonal bundle. Let $\gamma\in H^0(X, \Hom(\mathcal I\mathcal K^{-1}, \mathcal V)\otimes \mathcal K)$ and $\gamma^{\dagger}\in H^0(X, \Hom(\mathcal V, \mathcal I\mathcal K)\otimes \mathcal K)$ satisfying 
\[Q_0(\gamma^{\dagger} s, t)=Q_{\mathcal V}(s, \gamma t).\] Here $Q_0=\begin{pmatrix}0&1&0\\1&0&0\\0&0&1\end{pmatrix}$ is the orthogonal structure on $\mathcal I\mathcal K\oplus \mathcal I\mathcal K^{-1}\oplus \mathcal I$. Let $(\mathcal E, \phi)$ be the following Higgs bundle:
\begin{equation}
\begin{tikzcd}
&\mathcal{V}\arrow[dl, "\gamma^{\dagger}"']&\\
\mathcal{IK}\arrow[r, "1"]&\mathcal{I}\arrow[r, "1"] & \mathcal{IK}^{-1}\arrow[lu, "\gamma"'] &
\end{tikzcd}.
\end{equation}
\begin{assumption}\label{ass:MaximalSO(2,n+1)}
1. The harmonic metric 
\[H=H_{\mathcal{V}}\oplus g_\pi\oplus 1\oplus g_\pi^{-1},\] where $g_\pi$ is a Hermitian metric on $\mathcal{IK}^{-1}$ and $H_{\mathcal V}$ is compatible with $Q_{\mathcal V}$-structure.\\
2. The metric $g_\pi$ is complete. 
\end{assumption}

The following lemma was proved in \cite{LT} except the last line and later reproved in \cite{DaiLiBounded} using a different argument. For completeness, we reproduce the proof from \cite{DaiLiBounded}, adding only the last line to provide improved control in the nilpotent case.
\begin{lem}\label{lem:NormBoundMaximalSO(2,n+1)}Under Assumption \ref{ass:MaximalSO(2,n+1)}, \[K_{g_\pi}=2\|\gamma\|^2-2\leq 0.\] If the inequality holds at one point, it holds everywhere. \\
Moreover, the following are equivalent:
\begin{itemize}
    \item $X$ is equipped with a  complete conformal hyperbolic metric $g_{hyp}$\\ and $\sup_{X}|\tr(\varphi^{4})|_{g_{hyp}}$ is bounded;
    \item There exists a constant $\epsilon>0$ such that $K_{g_\pi}\leq -\epsilon.$
\end{itemize}
Furthermore, the constants $\sup_{X}|\tr(\varphi^{4})|_{g_{hyp}}$ and $\epsilon$ depend only on each other. 

If in addition, $\tr(\varphi^4)=0, (i.e., \gamma^{\dagger}\gamma=0),$ then $K_{g_\pi}\leq -\frac{2}{3}.$
\end{lem}
\begin{proof}
The Hitchin equation gives
\[i\Lambda_{g_\pi} (F_H|_{\mathcal{IK}^{-1}})=-i\Lambda_{g_\pi}(\gamma^*\wedge\gamma+1\wedge 1^*)=\|\gamma\|^2-1.\]
Thus \[K_{g_\pi}=2i\Lambda_{g_\pi}(F(g_\pi))=2i\Lambda_{g_\pi} (F_H|_{\mathcal{IK}^{-1}})=2\|\gamma\|^2-2.\]
Let $s_\gamma$ be the section of $\End(\mathcal E)\otimes \mathcal K$ given by $\gamma$.
From Lemma \ref{lem:EquationSection}, we have 
\[\Delta_{g_\pi}\log\|\gamma\|^2\geq \frac{\|[s_{\gamma},\varphi^*]\|^2-\|[s_{\gamma},\varphi]\|^2}{\|\gamma\|^2}+\frac{1}{2}K_{g_\pi}=2\|\gamma\|^2-\|1\|^2-\tfrac{\|\gamma^{\dagger}\gamma\|^2}{\|\gamma\|^2}+\frac{1}{2}K_{g_\pi}.\]

Therefore, using the fact that $\|AB\|^2\leq \|A\|^2\|B\|^2$ for two matrices $A, B$, we have 
\[\Delta_{g_\pi}\log\|\gamma\|^2\geq 3\|\gamma\|^2- 3.\]
Since $g_\pi$ is complete and $K_{g_\pi}\geq -2$, we can apply the Cheng-Yau maximum principle and obtain $\|\gamma\|^2\leq 1.$ 
We continue to calculate $\Delta_{g_\pi}K_{g_\pi}$:
\begin{eqnarray*}\Delta_{g_\pi}K_{g_\pi}&=&2\Delta_{g_\pi}(\|\gamma\|^2)\geq2\|\gamma\|^2\cdot \Delta_{g_\pi}\log(\|\gamma\|^2)\\
&\geq& 6\|\gamma\|^2\cdot (\|\gamma\|^2-1)\\
&=&\frac{3}{2}(K_{g_\pi}+2)\cdot K_{g_\pi}
\end{eqnarray*}
By the continuity, the inequality extends to the whole surface $X$. So we obtain either $K_{g_\pi}\equiv 0$ or $K_{g_\pi}<0$. 

Since $(\mathcal E, \varphi)$ is an $SO(2,n+1)$-Higgs bundle, $\tr(\varphi^{k})=0$ unless $4$ divides $k$ and $\tr(\varphi^{4k})=c_k(\tr(\varphi^4))^k$ for some constant $c_k$ only depending on $k$. Applying Lemma \ref{lem:NegativeCurvature1}, we obtain the statement. 

If we have additional condition $\gamma^{\dagger}\gamma=0$, then \[\Delta_{g_\pi}\log\|\gamma\|^2\geq 3\|\gamma\|^2- 2.\]
Applying the Cheng-Yau maximum principle, we obtain $\|\gamma\|^2\leq \frac{2}{3}.$
\end{proof}

\subsection{Nestedness of cones for \texorpdfstring{$\SO_o(2,n+1)$}{SOo(2,n+1)}-maximal bundles}
\label{subsec:Nestedness of cones for SON}
Given a unit vector for the metric $g_\pi$, $v\in T^1X$, in some orthonormal basis for $H$ we obtain :
\[\Psi(v)=\varphi(v)+\varphi^{*H}(v)=\begin{pmatrix}0&0&1&a^T\\0&0&1&\bar a^T\\1&1&0&0\\\bar a&a&0&0\end{pmatrix}.\]
In this expression $a=\alpha(v)$ and $b=\beta(v)$. Let us write the Fuchsian part of $\Psi$: \[\Psi_0(v)=\begin{pmatrix}0&0&1&0\\0&0&1&0\\1&1&0&0\\0&0&0&0\end{pmatrix}.\] We consider the following change of basis diagonalizing the Fuchsian part $\Psi_0$ of $\Psi$ the Higgs field:\[P=\begin{pmatrix}\frac{1}{2}&\frac{1}{2}&\tfrac{1}{\sqrt{2}}&0\\\frac{1}{2}&\frac{1}{2}&-\tfrac{1}{\sqrt{2}}&0\\
\tfrac{1}{\sqrt{2}}&-\tfrac{1}{\sqrt{2}}&0&0\\
0&0&0&I_n\end{pmatrix},\] where $a$ is a column vector. Set
\[\hat\Psi(v):=P\Psi(v)P^{-1}=
\begin{pmatrix}
\sqrt{2}&0&0&\Re(a^T)\\
0&-\sqrt{2}&0&\Re(a^T)\\
0&0&0&\sqrt{2}i\Im(a^T)\\
\Re(a)&\Re(a)&-\sqrt{2}i\Im(a)&0\end{pmatrix},\]\[\Psi_0(v):=P\Phi_0(v)P^{-1}.\]We work in this new basis, in which every vector $\vec{X}$ decomposes as $X=\vec{X_1}+\vec{X_2}+\vec{X_3}+\vec{Y}$. Using this basis depending on the unit vector $v\in TX$, we define the following hermitian form on $\mathcal{IK}\oplus \mathcal{IK}^{-1}\oplus \mathcal{I}\oplus \mathcal{V}$:
\begin{equation}
\label{eq:Definition qv hyperbolic}
q_v(\vec X, \vec X')=\Re(\vec X_1\cdot \overline{\vec X_1'}-\vec X_2\cdot \overline{\vec X_2'}-\vec X_3\cdot \overline{\vec X_3'}-\vec Y\cdot \overline{\vec Y}).
\end{equation}

This hermitian form $q_v$ does not depend on any choice of basis, as it is characterized by the hermitian metric $H$ and the eigenspaces of the endomorphism $\hat\Psi_0(v)$.
It will allow us to define cones $C_v\subset \mathbb{P}(\C^{n+3})$ and filled hyperboloids $\mathcal{H}_v\subset \C^{n+3}$ as follows:
\begin{equation}
\label{eq:Defn Cv,Hv}
C_v=\{[\vec X]|q_v(\vec X, \vec X)>0\}, \,\mathcal{H}_v=\{\vec X|q_v(\vec X, \vec X)\geq 1\}.
\end{equation}
These cones, when viewed in the same ambient space by using the flat connection associated to the harmonic bundle, are nested along the geodesic flow. This will be a consequence of the parallelism and of the following property:
\begin{prop}\label{prop:Nestedness SO2n} Under Assumption \ref{ass:MaximalSO(2,n+1)}, for all $\vec X\in \partial C_v$ meaning $q_v(\vec X,\vec X)=0$:    
    $$q_v(\hat\Psi(v) \vec X,\vec X)=\frac{1-\sqrt{1+\frac{1}{2}K_{g_\pi}}}{\sqrt{2}}\lVert\vec X\rVert^2\geq -\frac{1}{4\sqrt{2}}K_{g_\pi}\lVert\vec X\rVert^2.$$
    
Under Assumption \ref{ass:MaximalSO(2,n+1)}, for all $\vec X\in \partial \mathcal{H}_v$ meaning $q_v(\vec X,\vec X)=1$:    
    \begin{equation}\label{eq:q_v Hyperboloidmaximal}
        q_v(\hat\Psi(v)\vec X,\vec X)\geq \sqrt{2}. 
    \end{equation}
    If moreover if $K_{g_\pi}<0$ equality in \eqref{eq:q_v Hyperboloidmaximal} can only occur if $\hat\Psi_0(\vec X)=\vec X$.
\end{prop}
In the terminology that we introduce in Section \ref{sec:Dynamics}, the endomorphism $\Psi(v)\in\mathfrak{so}_{n+3}(\C)$ is at least $\frac{1}{4\sqrt{2}}|K_{g_\pi}|$-nesting the cone $C_v$ and $\sqrt{2}$-nesting the filled hyperboloid $\mathcal{H}_v$.

\begin{proof}
By a direct calculation, 
\begin{eqnarray*}
(*)=q_v(\hat\Psi(v)\vec X,\vec X)=\sqrt{2}(|x_1|^2+|x_2|^2)-2\Re(\Re(a)^T\cdot \vec Y\cdot \bar x_2+\sqrt{2}i\Im(a)^T\cdot \vec Y\cdot \bar x_3).
\end{eqnarray*}
By Lemma \ref{lem:NormBoundMaximalSO(2,n+1)}, $|a|^2\leq 1$.

For $\vec X$ satisfying $q(\vec X, \vec X)=0$,
\begin{eqnarray*}
(*)&=&\sqrt{2}(2|x_2|^2+|x_3|^2+|\vec Y|^2)-2\Re(\Re(a)^T\cdot \vec Y\cdot \bar x_2+\sqrt{2}i\Im(a)^T\cdot \vec Y\cdot \bar x_3)\\
&\geq&\sqrt{2}(2|x_2|^2+|x_3|^2+|\vec Y|^2)-2(|\Re(a)|\cdot |\vec Y|\cdot |x_2|+\sqrt{2}|\Im(a)|\cdot |\vec Y|\cdot |x_3|).\end{eqnarray*}

We consider the function $Q$:
\[Q=\frac{2\sqrt{2}X^2+\sqrt{2}Y^2+\sqrt{2}Z^2-2(cXZ+\sqrt{2} dYZ)}{2(X^2+Y^2+Z^2)}\] for $X,Y,Z\in \R$ with $(X,Y,Z)\neq (0,0,0)$ given a fixed pair of $c,d$ satisfying $c^2+d^2\in [0,1].$ By a direct analysis of the minimum eigenvalue of the associated symmetric $3\times 3$ matrix, the global minimum of $Q$ is $\frac{1-\sqrt{c^2+d^2}}{\sqrt{2}}.$ We apply this inequality to obtain \[(*)\geq \frac{1-|a|}{\sqrt{2}}\lVert\vec X\rVert^2=\frac{1-\sqrt{1+\frac{1}{2}K_{g_\pi}}}{\sqrt{2}}\lVert\vec X\rVert^2\geq -\frac{1}{4\sqrt{2}}K_{g_\pi}\lVert\vec X\rVert^2,\]since $K_{g_\pi}=2|a|^2-2.$

Now we consider the case when $\vec X\in  \partial \mathcal{H}_v:$
\begin{eqnarray*}
(*)=\sqrt{2}+\sqrt{2}(2|x_2|^2+|x_3|^2+|\vec Y|^2)-2\Re(\Re(a)^T\cdot \vec Y\cdot \bar x_2+\sqrt{2}i\Im(a)^T\cdot \vec Y\cdot \bar x_3).
\end{eqnarray*}
This quantity is greater than $\sqrt{2}$ as in the previous statement. If $K_{g_\pi}<0$, then we have a strict inequality except if $\vec X=\vec X_1$, which means exactly that $\hat\Psi_0(\vec X)=\vec X$.
\end{proof}

\subsection{Nilpotent \texorpdfstring{$\SL(n+4,\mathbb{R})$}{SL(n+4,R)}-Higgs bundles}
\label{subsec:NilpotentSLn+4} 
Let $(\mathcal V, Q_{\mathcal V})$ be an orthogonal rank $n$ vector bundle, $\mathcal L$ be a holomorphic line bundle. Consider the following rank $n+4$ orthogonal bundle
\[
(\mathcal E,Q_{\mathcal E})=(\mathcal V,Q_{\mathcal V})\oplus 
(\mathcal L\mathcal K\oplus \mathcal L\oplus \mathcal L^{-1}\oplus \mathcal L^{-1}\mathcal K^{-1}, Q_0=\begin{pmatrix}0&0&0&1\\0&0&1&0\\0&1&0&0\\1&0&0&0\end{pmatrix}).
\] Let \[\mu\in H^0(X,\mathcal L^{-2}\mathcal K),~ 
\nu\in H^0(X,\Hom(\mathcal V,\mathcal L\mathcal K)\otimes \mathcal K),
~
\nu^\dagger \in H^0(X,\Hom(\mathcal L^{-1}\mathcal K^{-1},\mathcal V)\otimes \mathcal K),
\]
satisfying 
\[Q_V(\nu^{\dagger}s,t)=Q_0(s,\nu t),\quad 
\nu\nu^\dagger=0.
\]
Let $(\mathcal E, \varphi)$ be the following Higgs bundle:
\begin{equation}
\label{eq:cyclicSO(n+4)}
\begin{tikzcd}
\mathcal{V}\arrow[r, "\nu"]&\mathcal{KL}\arrow[r, "1"] & \mathcal{L}\arrow[r, "\mu"] &\mathcal{L}^{-1}\arrow[r, "1"]&\mathcal{L}^{-1}\mathcal K^{-1}\arrow[bend right,llll, "\nu^{\dagger}"]&
\end{tikzcd}.
\end{equation}
Thus $(\mathcal E, \varphi)$ is a nilpotent $\SL(n+4,\mathbb R)$-Higgs bundle.

\begin{rem}
 When $(\mathcal V, Q_{\mathcal V})$  is an $\SO(2,\mathbb C)$-bundle $(\mathcal N\oplus\mathcal N^{-1},\begin{pmatrix}0&1\\1&0\end{pmatrix})$, it gives an $\SL(6,\R)$-VHS:
 \begin{equation}
\label{eq:cyclicSL6}
\begin{tikzcd}
\mathcal{N}\arrow[r, "\nu"]&\mathcal{KL}\arrow[r, "1"] & \mathcal{L}\arrow[r, "\mu"] &\mathcal{L}^{-1}\arrow[r, "1"]&\mathcal{L}^{-1}\mathcal K^{-1}\arrow[r, "\nu"]&\mathcal N^{-1}&
\end{tikzcd}
\end{equation}
\end{rem}
\begin{assumption}\label{ass:NilpotentSO(n+4)}
1. Assume the harmonic metric splits $H=H|_{\mathcal V}\oplus H|_{\mathcal{LK}}\oplus H|_{\mathcal L}\oplus H|_{\mathcal L^{-1}}\oplus H|_{\mathcal{L}^{-1}\mathcal K^{-1}}$, and $H$ is compatible with the $Q_{\mathcal E}$-structure, i.e. \[H=H_{\mathcal V}\oplus H_1\oplus H_1g_\pi\oplus H_1^{-1}g_\pi^{-1}\oplus H_1^{-1},\] where $H_{\mathcal V}$ is compatible with the $Q_{\mathcal V}$-structure, $H_1$ is a Hermitian metric on $\mathcal{KL}$, and $g_\pi$ is a Hermitian metric on $\mathcal K^{-1}$.\\
2. The metric $g_\pi$ is complete.
\end{assumption}
\begin{lem}\label{lem:NormBoundSO(n+4)}
Under Assumption \ref{ass:NilpotentSO(n+4)}, $K_{g_\pi}<-\frac{1}{6}$, equivalently   \[\|\mu\|^2+\|\nu\|^2< \frac{23}{12}.\]
\end{lem}
\begin{proof}
The Hitchin equation gives
\[i\Lambda_{g_\pi} (F_H|_{\mathcal{KL}})=-i\Lambda_{g_\pi}(\nu\wedge\nu^*+1^*\wedge 1)=1-\|\nu\|^2.\]
\[i\Lambda_{g_\pi} (F_H|_{\mathcal{L}})=-i\Lambda_{g_\pi}(\mu^*\wedge\mu+1\wedge 1^*)=\|\mu\|^2-1.\]
Thus \[K_{g_\pi}=2i\Lambda_{g_\pi}(F(g_\pi))=2i\Lambda_{g_\pi} (-F_H|_{\mathcal{KL}}+F_H|_{\mathcal{L}})=2\|\mu\|^2+2\|\nu\|^2-4.\]

Let $s_\mu, s_\nu$ be the section of $\End(\mathcal E)\otimes \mathcal K$ given by $\mu$, $\nu$ respectively. 
From Lemma \ref{lem:EquationSection},
\[\Delta_g\log\|\mu\|^2\geq \frac{\|[s_{\mu},\phi^*]\|^2-\|[s_{\mu},\phi]\|^2}{\|\mu\|^2}+\frac{1}{2}K_g=2\|\mu\|^2-2+\frac{1}{2}K_g\geq 3\|\mu\|^2-4.\]
\[\Delta_g\log\|\nu\|^2\geq \frac{\|[s_{\nu},\phi^*]\|^2-\|[s_{\nu},\phi]\|^2}{\|\nu\|^2}+\frac{1}{2}K_g=2\|\nu\|^2-1+\frac{1}{2}K_g\geq 3\|\nu\|^2-3.\]
Note the equation of $\Delta_g\log\|\nu\|^2$ make uses of the assumption $\nu\nu^{\dagger}=0$.
Since $g$ is complete and $K_g\geq -4$,
we can apply the Cheng-Yau maximum principle to both equations and obtain $\|\mu\|^2\leq \frac{4}{3}$ and $\|\nu\|^2\leq 1.$

Let $s$ be the section of $\End(\mathcal E)\otimes \mathcal K$ given by $\mu$ and $\nu$. 
From Lemma \ref{lem:EquationSection},
\begin{eqnarray*}
\Delta_{g_\pi}\log\|s\|^2&\geq& \frac{\|[s,\varphi^*]\|^2-\|[s,\varphi]\|^2}{\|s\|^2}+\frac{1}{2}K_{g_\pi}\\
&=&\frac{2(\|\mu\|^4+\|\nu\|^4)-(2\|\mu\|^2+\|\nu\|^2)}{\|\mu\|^2+\|\nu\|^2}+\frac{1}{2}K_{g_\pi}\\
&\geq& 2(\|\mu\|^2+\|\nu\|^2)-3-\frac{\|\mu\|^2}{\|\mu\|^2+\|\nu\|^2}\\
&\geq& 2(\|\mu\|^2+\|\nu\|^2)-3-\frac{\frac{4}{3}}{\|\mu\|^2+\|\nu\|^2}
\end{eqnarray*}
Note that the second equation uses the assumption $\nu\nu^{\dagger}=0$ and the last equation uses $\|\mu\|^2\leq \frac{4}{3}.$
Since $g_\pi$ is complete and $K_{g_\pi}\geq -4$,
we can apply the Cheng-Yau maximum principle to both equations and obtain $\|\mu\|^2+\|\nu\|^2\leq \frac{9+\sqrt{177}}{12}<\frac{23}{12}.$
\end{proof}

At an orthonormal basis $\mathbf e$ for $H$ and a unit vector $v$ of $TX$ with respect to $g_\pi$, we obtain 
\[
\Psi(v)
=
\begin{pmatrix}
0&\bar a^T&0&0&a^T\\
a&0&1&0&0\\
0&1&0&\bar b&0\\
0&0&b&0&1\\
\bar a&0&0&1&0
\end{pmatrix},
\]
where $a$ is a row vector satisfying $a^Ta=0$. Let $P=\text{diag}(Id, P_0,P_0)$, where $P_0=\begin{pmatrix}\frac{1}{\sqrt{2}}&\frac{1}{\sqrt{2}}\\-\frac{1}{\sqrt{2}}&\frac{1}{\sqrt{2}}\end{pmatrix}$. 
With respect to the orthonormal frame $P\mathbf e$, write the fiber element in terms of coordinate
\[
\vec X=(x_1,\dots,x_{n+4})^T.\]
Define the cone $C_v$ as follows:
\[
C_v=
\left\{\vec X|q(\vec X, \vec X)>0\right\}\]
where $q(\vec X, \vec X')=
\Re(x_{n+1}\overline{x}_{n+1}'+x_{n+3}\overline{x}_{n+3}'-
\vec Y\cdot \overline{\vec Y'}
-x_{n+2}\bar x_{n+2}'
-x_{n+4}\bar x_{n+4}').$

\begin{prop}\label{prop:Nestedness SLn+4}
Under Assumption \ref{ass:NilpotentSO(n+4)}, for all $\vec X\in \partial C_v$ satisfying $\|\vec X\|=1$, 
    \[q(P\Psi(v)P^{-1} \vec X,\vec X)\geq 0.01. \]
\end{prop}
\begin{proof}
By a direct calculation, \[P\Psi(v)P^{-1}=
\begin{pmatrix}
 0&\frac{1}{\sqrt{2}}\bar a^T & -\frac{1}{\sqrt{2}}\bar a^T&\frac{1}{\sqrt{2}} a^T&\frac{1}{\sqrt{2}} a^T\\
 \frac{1}{\sqrt{2}} a&1&0&\frac{\bar b}{2}&-\frac{\bar b}{2}\\
  -\frac{1}{\sqrt{2}} a&0&-1&\frac{\bar b}{2}&-\frac{\bar b}{2}\\
   \frac{1}{\sqrt{2}} \bar a&\frac{b}{2}&\frac{b}{2}&1&0\\
      \frac{1}{\sqrt{2}} \bar a&-\frac{b}{2}&-\frac{b}{2}&0&-1\\
\end{pmatrix}.
\]
Let
$\vec Y=(x_1,\dots,x_n)^T.$ Then
\begin{eqnarray*}
(*)&=&q(P\Psi(v)P^{-1} \vec X,\vec X)\\
&=&
|x_{n+1}|^2
+|x_{n+2}|^2
+|x_{n+3}|^2
+|x_{n+4}|^2\\
&&+\Re(\bar bx_{n+3}\bar x_{n+1}+\sqrt{2} a\cdot \vec Y\cdot \bar x_{n+2}-\sqrt{2}\bar a\cdot \vec Y\cdot \bar x_{n+4}+bx_{n+2}\bar x_{n+4}).
\end{eqnarray*}
By Lemma \ref{lem:NormBoundSO(n+4)}, 
$|a|^2+|b|^2\leq \frac{23}{12}$. Then we have 
\[
(*)
\geq (1-\tfrac{|b|}{2})(|x_{n+1}|^2
+|x_{n+2}|^2
+|x_{n+3}|^2
+|x_{n+4}|^2)\\
+\Re(\sqrt{2} a\cdot \vec Y\cdot \bar x_{n+2}-\sqrt{2}\bar a\cdot \vec Y\cdot \bar x_{n+4}).
\]
Since $\vec X\in \partial C_v$ and $|\vec X|=1$, we have 
\begin{eqnarray*}
(*)
&\geq& 0.01+(1.98-|b|)(|x_{n+2}|^2
+|x_{n+4}|^2)+(0.98-\tfrac{|b|}{2})|\vec Y|^2\\
&&+\Re(\sqrt{2} a\cdot \vec Y\cdot \bar x_{n+2}-\sqrt{2}\bar a\cdot \vec Y\cdot \bar x_{n+4})\\
&=& 0.01+
(0.98-\tfrac{|b|}{2})
\Big|\vec Y+\frac{\bar a^Tx_{n+2}}{\sqrt{2}(0.98-\tfrac{|b|}{2})}-\frac{a^Tx_{n+4}}{\sqrt{2}(0.98-\tfrac{|b|}{2})}\Big|^2\\
&&+(1.98-|b|-\frac{|a|^2}{2(0.98-\tfrac{|b|}{2})})(|x_{n+2}|^2+|x_{n+4}|^2)
\end{eqnarray*}
Here, we use the fact that $aa^T=0$. By Lemma \ref{lem:NormBoundSO(n+4)}, $|a|^2+|b|^2\leq \frac{23}{12}$ and thus \[1.98-|b|-\frac{|a|^2}{2(0.98-\tfrac{|b|}{2})}\geq 1.98-|b|-\frac{\frac{23}{12}-|b|^2}{2(0.98-\tfrac{|b|}{2})}>0,\]
by using $|b|\leq 2/\sqrt{3}$.
Therefore, $(*)\geq 0.01.$
\end{proof}

\subsection{Yano holomorphic curves}
 We now discuss the differential geometric properties of the equivariant map defining the parallel distributions of planes, as well as an associated map in a smaller homogeneous space.

While the whole proof of our main theorems can be written down using only the parallel distribution and exhibiting the family of multicones, we believe it is insightful to motivate the construction via
an auxiliary harmonic map to a pseudo-Riemannian symmetric space.

\subsubsection{Yano structures on homogeneous spaces}
Here we describe the notion of Yano-holomorphic map,
and relate it to harmonicity. Yano structures were first considered by Yano~\cite{Yan63}. They were of particular use for the study of minimal surfaces in compact homogeneous spaces, see for instance Rawnsley~\cite{Raw85}. This Section aims at describing better the equivariant objects that we construct.

\begin{defn}
Let $M$ be a manifold. A Yano structure on $M$
is a bundle endomorphism~$J\in\rm{End}(TM)$ satisfying
\begin{equation*}
    J^3+J=0\,.
\end{equation*}
A map $f:(M_1,J_1)\rightarrow(M_2,J_2)$ between
manifolds endowed with Yano structures is Yano-holomorphic (or $(J_1,J_2)$-holomorphic)
if it satisfies
\begin{equation*}
    df\circ J_1=J_2\circ df\,.
\end{equation*}
\end{defn}
Note that trivially any almost complex structure
on $M$ defines a Yano structure.
A Yano structure leads to a splitting of $TM$
into $\rm{Ker}(J)$ and $\rm{Ker}(J^2+\rm{Id})$,
and $J$ restricts to an almost-complex structure on the
second subbundle.

\medskip

Let $L\subset G$ be a reductive Lie subgroup with Lie algebra $\mathfrak{l}$.  

\begin{defn}
    Let $(G/L,J)$ be a reductive $G$-homogeneous space equipped with a Yano structure, and let $J_+,J_-$ be the $i$ and $-i$ eigenspace of $J$ on the orthogonal in $\mathfrak{g}_\C$ of $\mathfrak{l}_\C$. It is \emph{horizontal} if:
    $$[J_+,J_-]\subset \mathfrak{l}\otimes\C.$$
\end{defn}

A horizontal Yano structure on a reductive $G$-homogeneous space satisfies the condition~$(d^\nabla J)^{(1,1)}=0$ for any Levi-civita connection $\nabla$ of a $G$-invariant pseudo-metric on $G/L$.
As such, it satisfies the following properties, see Black~\cite{Bla91}, for a detailed discussion.
\begin{prop}
Let $(M,g,\omega)$ be a Kähler-manifold, $\rho:\pi_1M\rightarrow G$
and $f:\widetilde M\rightarrow G/L$ an equivariant
$J$-holomorphic map. Then
\begin{enumerate}
\item
For any invariant metric $\lambda$ on $G/L$,
$f$ is harmonic. (We say that $f$ is \emph{equi-harmonic}).
\item
For any reductive subgroup $L'$ containing $L$ and any invariant
metric $\lambda$ on $G/L'$, the composition $\pi\circ f$ is a harmonic equivariant map to $G/L'$.
\end{enumerate}
\end{prop}

\subsubsection{Minimal Slodowy slice}

Let $\lbrace e,h,f\rbrace$ be a minimal $\mathfrak{sl}_2$-triple in $\mathfrak{g}$, let $T,C$ be the associated subgroups of $G$ and let $K$ be a maximal compact subgroup containing a maximal compact subgroup of $TC$. We construct a Yano structure on $G/(TC\cap K)$ so that the parallel distribution of planes coming from parallelisable harmonic bundles correspond exactly to Yano-holomorphic maps.

\medskip

Since $\mathfrak g$ is a complex reductive Lie algebra, $\mathfrak k_\C$
is isomorphic as complex Lie algebra to $\mathfrak g$, via the map $u+Iv\in\mathfrak k_\C\mapsto u+iv\in \mathfrak g$.
The space $\mathfrak p_\C$ is identified, via the map $u+Iv\mapsto u-iv$ to $\mathfrak g$.
Under this antilinear identification, the bracket of two elements in $\mathfrak{p}_\C$ which belongs to $\mathfrak{k}_\C$, is the standard bracket on $\mathfrak{g}$.
We now consider the Yano structure $J$ defined such that:
\begin{flalign*}
    J_+=\rm{Ker}(J-i\rm{Id})&=\mathfrak g_{-2}\oplus \mathfrak g_1\subset\mathfrak p_\C\\
    J_-= \rm{Ker}(J+i\rm{Id})&=\mathfrak g_2\oplus \mathfrak
    g_{-1}\subset\mathfrak p_\C\\
    \rm{Ker}(J)&=(\mathfrak g_0\subset \mathfrak p_\C) \oplus \mathfrak{k}_\C\
\end{flalign*}
The conjugation on $\mathfrak p_\C$ sends $J_+=\rm{Ker}(J-i\rm{Id})$
to $J_-=\rm{Ker}(J+i\rm{Id})$, so $J$ is a real endomorphism of $T(G/(TC\cap K)$.
Note that $[\mathfrak{g}_2,\mathfrak{g}_1]=[\mathfrak{g}_{-2},\mathfrak{g}_{-1}]=0$ and that $[\mathfrak{g}_1,\mathfrak{g}_{-1}],[\mathfrak{g}_2,\mathfrak{g}_{-2}]\subset \mathfrak{g}_0$. We deduce that $[J_+,J_-]\subset \mathfrak k_\C\cap\mathfrak g_0=\text{Lie}(TC\cap K)_\C$,
hence this Yano structure is horizontal.
\begin{prop}
    Let $(\calE,\varphi=\mathbf{e}+\mathbf{r})$ be a Higgs bundle in one of the minimal Slodowy slices discussed in Section 2. Then the parallel distribution of hyperbolic planes yields a horizontal lift $Y:\widetilde X\rightarrow G/(TC\cap K)$ that is Yano-holomorphic.
    Its differential satisfies
    \begin{equation}
    (dY)_\C=\varphi+\varphi^\ast=\mathbf{e}+\mathbf{e^\ast}+\mathbf{r}+\mathbf{r^\ast}\,.
    \end{equation}
\end{prop}
\begin{proof}
    Let $D=\nabla+\varphi+\varphi^\ast$ be the flat connection.
    The harmonic map to the symmetric space is given by
    the metric $(H_x)_{x\in X}$, whose Lie algebras
    are locally given by $\mathfrak{su}(H_x)\cap\mathfrak g$. The lift $Y$ is locally given by the
    Lie algebras $\mathfrak{g}_{-2}\oplus\mathfrak g_0\oplus\mathfrak g_2\subset\mathfrak k_\C\approx\mathfrak g$.
    Projecting $D$ on the orthogonal of these subalgebras
    gives the explicit formula for $(dY)_\C$, as desired.
    Observe that the fact that $\nabla$ is valued in this space is precisely the parallelity condition on the family of planes.

    By construction, under the identification $\mathfrak p_\C\approx\mathfrak g$, we have $J_+=\mathfrak g_{-2}\oplus\mathfrak g_1$. As $\mathbf{e}\in\mathfrak g_{-2}$ and $\mathbf{r}\in\mathfrak g_1$, we have
    $(dY)_\C^{1,0}\in J_+$, so $Y$ is Yano-holomorphic, as claimed.
\end{proof}
\begin{rem}
    There is a slight abuse of notation here, as
    we consider $\varphi$ as a one-form valued in $((\mathfrak t\oplus\mathfrak c)\cap\mathfrak k)^\perp$, where $\mathfrak t=\langle h\rangle$ is the Lie algebra of $T$. This is rigorous as we previously fixed
    an identification $\mathfrak p\hookrightarrow\mathfrak g$.
\end{rem}
\begin{cor}
    The composition $u=\pi\circ Y$ with $\pi:G/(TC\cap K)\rightarrow G/G^\sigma$ yields a maximal spacelike harmonic immersion, whose induced metric is the projected metric, and which is pointwise tangent to a totally geodesic hyperbolic plane in $G/G^\sigma$. 
\end{cor}
\begin{proof}
    The composition of $dY$ with $d\pi$ yields the expression
    \begin{equation*}
        (du)_\C=\mathbf{e}+\mathbf{e^\ast}
    \end{equation*}
from which we deduce that it is an immersion,
spacelike for the Killing pseudo-metric
as $\mathbf{e}+\mathbf{e^\ast}\in\mathfrak p_\C$.

The harmonicity is equivalent to the fact that $\mathbf{e}$
is holomorphic, or can be deduced from the fact that $Y$ is Yano-holomorphic for a horizontal Yano structure.

Its induced metric is $\mathbf e\wedge\mathbf e^\ast$, which is precisely the projected metric. As it is a conformal harmonic spacelike immersion, it is a maximal immersion.

Finally, it is pointwise tangent to a hyperbolic plane,
which is the exponential of the $\mathfrak{sl}_2$-subalgebra given by $\mathfrak t\subset\mathfrak g^\sigma,\,\mathfrak g_2\oplus \mathfrak g_{-2}\subset\mathfrak p_\C$.
\end{proof}

In order to clarify the general picture, and in particular Figure \ref{fig:Diagram} from the introduction, note the following.

\begin{fact}
    The cone map $\mathcal V:T^1\widetilde X\rightarrow\{\text{Cones in }\mathbb P(V)\}$, $v\mapsto C_v$
    factors through a $G$-equivariant map
    $G/C\cap K\rightarrow\{\text{Cones in }\mathbb P(V)\}$.
\end{fact}
\begin{proof}
    For the construction of the cones, we use a representation~$G\hookrightarrow \rm{SL}(V)$
    such that $V=V_{-1}\oplus V_1\oplus V_{0}$. For every $v\in T^1\widetilde{X}$, $\mathbf{e}(v)$ defines a conformal linear homomorphism $e:V_{1}\rightarrow V_{-1}$ and $\mathbf{r}$ defines a map $r\in\rm{Hom}(V_1,V_0)\oplus\rm{Hom}(V_0,V_{-1})$.
    Since $e$ is conformal, this yields a totally geodesic embedding $G/TC\cap K\rightarrow \SL(V)/\rm{S}(U(1)\times U(V_1)\times U(V_0))$
    where $U(1)\times U(V_1)\times U(V_0)\rightarrow\rm{GL}(V)$ via $(\lambda,A,B)\mapsto A\oplus\lambda A\oplus B$.
    The subgroup $C\cap K$ is mapped to matrices of
    the form $A\oplus A\oplus B$,
    which yields a totally geodesic embedding
    $G/(C\cap K)\rightarrow\rm{SL}(V)/\rm{S}(U(V_1)\times U(V_0))$.

    The cone map is now constructed in the following way: A point $x\in\rm{SL}(V)/S(U(V_1)\times U(V_0)$
    is represented by a decomposition $V=V_1\oplus V_{-1}\oplus V_0$ and a quadruple $(e^*,H_1,H_{-1},H_0)$
    where $H_i$ is a Hermitian metric on $V_i$ and $e^*:V_{-1}\rightarrow V_{1}$ is a linear isometry.
    The associated cone $C_v\subset\mathbb P(V)$ where $x=\mathcal{V}(v)$ can be rewritten the following way:
    \begin{multline*}
        C_v=\left\{\vec{X}\in V \middle\vert \frac{1}{2}\|e^*(\vec{X}_{-1})+\vec{X}_{{1}}\|_{H_{1}}^2\geq   
        \frac{1}{2}\|e^*(\vec{X}_{-1})-\vec{X}_{{1}}\|_{H_{1}}^2+\|\vec{X}_{0}\|_{H_V}^2\right\}.
    \end{multline*}

Here we decompose $\vec{X}=\vec{X}_{1}+\vec{X}_{-1}+\vec{X}_0$ with $\vec{X}_i\in V_i$.
\end{proof}

\section{Regularity and domination}
\label{sec:Dynamics}

In this section we show how to obtain dynamical properties and weight domination for our representations, as well as the $\Theta$-undistorted property and quasi-isometric embedding of the equivariant surfaces. We first develop the notion of nestedness for filled hyperboloids, and then develop the notion of nested cones from \cite{DE26} that we extend in order to get more precise estimates also in the non-compact case.

\medskip

Throughout this section we denote by $\mu:\X\times\X\to\mathfrak{a}^+$ the \emph{Cartan projection}, also called vector-valued distance. Fixing a point $o\in\X$ with stabilizer $K$ and a Cartan decomposition $G=K\exp(\mathfrak{a}^+)K$, it is defined by $\mu(g\cdot o,g'\cdot o)=a$ whenever $g^{-1}g'\in K\exp(a)K$ with $a\in\mathfrak{a}^+$. When $G=\SL(n,\C)$ and $x\in\X$ corresponds to a hermitian metric on $\C^n$, $\mu(x,g\cdot x)$ records the logarithms of the singular values of $g$ with respect to this metric. We denote by $\lambda:G\to\mathfrak{a}^+$ the Jordan projection, which satisfies $\lambda(g)=\lim_{k\to+\infty}\frac{1}{k}\mu(x,g^k\cdot x)$ for any $x\in\X$.

\subsection{Nestedness of hyperboloids}

In order to obtain domination of certain weights, we will consider the nestedness of filled hyperboloids in $\C^n$. Let $1\leq d\leq n$ be an integer. A (filled) $d$-hyperboloid is a set of the form $\mathcal{H}=\lbrace v\in \C^n\mid Q(v,v)\geq 1\rbrace$ for some hermitian bilinear form $Q$  of signature $(d,n-d)$. 

\medskip

Given two hyperboloids $\mathcal{H}_1,\mathcal{H}_2\subset \C^n$, we define the nestedness $\Nest(\mathcal{H}_1,\mathcal{H}_2)$ as the supremum of the set of $C\geq 1$ such that $\mathcal{H}_2\subset C \mathcal{H}_1$. 

Let $\lambda:\SL(n,\C)\to \mathfrak{a}^+$ be the Jordan projection, and $\lambda_1,\cdots, \lambda_n$ the characteristic values of $g$ ordered with decreasing modulus.

\begin{prop}
\label{prop:HyperboloidNestingControl}
    Let $g\in G$ and $\mathcal{H}\subset \C^n$ be a $d$-hyperboloid such that $\log\Nest(\mathcal{H},g\cdot \mathcal{H})>0$. Then
    $\lambda_d(g)\geq \log\Nest(\mathcal{H},g\cdot \mathcal{H})$.
\end{prop}

\begin{proof}

Let us write $\eta=\log\Nest(\mathcal{H},g\cdot \mathcal{H})$. Consider the anti-linear pairing $q$ such that $\mathcal{H}=\lbrace x\in \C^n\mid q(x,x)\geq 1\rbrace$. Consider a hermitian metric $h$ on $\C^n$ compatible with $q$, meaning that in some $h$-unitary basis:
$$q(x,x)=|x_1|^2+\cdots +|x_d|^2-|x_{d+1}|^2-\cdots-|x_n|^2.$$
Let $E\subset \C^n$ be the space spanned by the first $d$ vectors in this basis. For every $x\in E$ with $h(x,x)=1$, then $g(x)\in e^\eta \mathcal{H}$, meaning that $q(g(x),g(x))\geq e^{2\eta}$. Hence $h(g(x),g(x))\geq e^{2\eta}$ for all such $x$, which means that the logarithm of the $d$ largest singular values $\sigma_1(g)\geq \cdots\geq \sigma_d(g)$ of $g$ with respect to the hermitian metric $h$ satisfy $\sigma_d(g)\geq \eta$. 

We apply the previous argument to $g^n$. Hence $\frac{1}{n}\sigma_d(g^n)\geq \eta$, but $\frac{1}{n}\sigma_d(g^n)$ converges to $\lambda_d$, which concludes the argument.
\end{proof}

\begin{defn}
    \label{defn:Infinitesimal Nesting Cone}
    Given an hyperboloid $\mathcal{H}=\lbrace x\mid q(x,x)\geq 1\rbrace$ for a hermitian form of signature $(d,n-d)$. We say that an endomorphism $A$ of $\C^n$ is $\epsilon$-nesting $\mathcal{H}$ if for all $x\in \C^n$ such that $q(x,x)=1$:
    $$q(Ax,x)\geq \epsilon.$$
\end{defn}

\begin{prop}
\label{prop:HyperboloidNestingInfinitesimal}
    Let $(g_t)_{t\in \R}$ be a curve in $\SL(n,\C)$, $\mathcal{H}_0$ be a $d$-hyperboloid defined by a Hermitian form $Q_0$ and define $Q_t=g_t\cdot Q_0$, $\mathcal{H}_t=g_t\cdot \mathcal{H}_0$ for all $t\in \R$. Suppose that for all $t_0\in \R$, $g_{t_0}$ is $\epsilon_{t_0}$-nesting $\mathcal{H}_{t_0}$.

Then for all $t<t'$:
$$\log\Nest(\mathcal{H}_t,\mathcal{H}_{t'})\geq\int_{t}^{t'}\epsilon_s\mathrm{d}s.$$
\end{prop}

\begin{proof}

By homogeneity, our assumption implies that for all $t_0\in\R$ and all $v\in\C^n$ with $Q_{t_0}(v,v)>0$, $\frac{\mathrm{d}}{\mathrm{d}t}_{\mid t_0}Q_{t}(v,v)\leq -2\epsilon_{t_0}Q_{t_0}(v,v).$ Let $t<t'$ and $v\in \mathcal{H}_{t'}$. The function $f(s)=Q_s(v,v)$ satisfies $f(t')\geq 1$, and the functional inequality $f'\leq -2\epsilon_sf$ wherever $f>0$; in particular $f$ is non-increasing on any interval where it is positive, so $f>0$ on $[t,t']$, and \[f(t)\geq f(t')\exp\left(\int_{t}^{t'}2\epsilon_s\mathrm{d}s\right)\geq \exp\left(\int_{t}^{t'}2\epsilon_s\mathrm{d}s\right).\] Since $e^{\eta}\mathcal{H}_t=\lbrace w\mid Q_t(w,w)\geq e^{2\eta}\rbrace$, this means that $v\in e^{\int_{t}^{t'}\epsilon_s\mathrm{d}s}\mathcal{H}_t$. Hence $\mathcal{H}_{t'}\subset e^{\int_{t}^{t'}\epsilon_s\mathrm{d}s}\mathcal{H}_t$, which is the required control on $\Nest(\mathcal{H}_t,\mathcal{H}_{t'})$.
\end{proof}

\begin{rem}
    Even though this notion is way simpler than nestedness of cones in projective space, it controls the $d$-th singula value instead of the roots, and hence it does not allow us to prove the Anosov property.
\end{rem}

\subsection{Cross ratio nestedness of cones}
Let us fix $n\geq 2$ and $1\leq d<n$. In what follows we will call a \emph{$d$-cone} a subset $C\subset \mathbb{CP}^{n-1}$ which can be written for some basis of $\C^n$:
$$C=\lbrace [x_1,\cdots,x_n]\in \mathbb{CP}^{n-1} \mid |x_1|^2+\cdots +|x_d|^2>|x_{d+1}|^2+\cdots +|x_{n}|^2\rbrace.$$

\begin{rem}
In the present paper, $d$-cones will always be such quadratic cones. However, one can also obtain eigenvalue gaps by the nesting of general open sets containing the projectivisation of a $d$-dimensional subset and whose closure is disjoint from a $(n-d)$-dimensional set. 
\end{rem}

In other words $C$ is the locus of positive lines for a pseudo-hermitian metric of signature $(d,n-d)$ on $\C^n$. We say that a hermitian scalar product $h$ is \emph{compatible} if $C$ can be written as above for a $h$-unitary metric.

\medskip

A \emph{real projective line} $L$ in $\mathbb{CP}^n$ is the projectivization of a real plane $P\subset \C^n$ that is not a complex line. Given any $4$ points $[x_1,x_2,y_2,y_1]$ on such a real projective line we can define their \emph{cross ratio} which is invariant by projective transformations. If $x_1,x_2,y_2,y_1$ correspond in an affine chart to $X_1,X_2,Y_2,Y_1\in \R$ this cross ratio is equal to:
$$[x_1,x_2,y_2,y_1]=\frac{Y_2-X_1}{Y_2-X_2}\times\frac{Y_1-X_2}{Y_1-X_1}.$$

\medskip

First let us measure the nestedness of a cone along a real projective line.

\begin{defn}
Let $C_1,C_2$ be $d$-cones with $\overline{C_2}\subset C_1$, and let $L$ be a projective line \emph{crossing} $(C_1,C_2)$, \ie such that $L$ intersects $C_1$ and the complement of $\overline{C_2}$. This line intersects $\partial C_1$ and $\partial C_2$ respectively in $(x_1,x_2,y_2,y_1)\in L$, which we can assume are cyclically in this order. We define:
$${\log\Nest}_L(C_1,C_2)=\log[x_1,x_2,y_2,y_1].$$
\end{defn}

\begin{defn}
Let $C_1,C_2$ be $d$-cones with $\overline{C_2}\subset C_1$, we define their \emph{nestedness} as the infimum ${\log\Nest}(C_1,C_2)$ for all real projective lines $L$ crossing $(C_1,C_2)$ of ${\log\Nest}_L(C_1,C_2)$.
\end{defn}

\begin{rem}
    As in \cite{DE26b} we will simply write $\Nest=\exp({\log\Nest})$.
\end{rem}

This definition is convenient as it satisfies the following:

\begin{prop}[{\cite{Dav26,DE26}}]
\label{prop:Additivity Nestedness}
    Let $C_1,C_2$ be $d$-cones with $\overline{C_2}\subset C_1$, then ${\log\Nest}(C_1,C_2)>0$. Let $C_1,C_2,C_3$ be $d$-cones with $\overline{C_2}\subset C_1$ and $\overline{C_3}\subset C_2$, then $${\log\Nest}(C_1,C_2)+{\log\Nest}(C_2,C_3)\leq {\log\Nest}(C_1,C_3).$$ 
\end{prop}

If an element of $\GL_n(\C)$ nests a cone then one can control some of its eigenvalue gap.

\begin{prop}[{\cite{DE26}}]
\label{prop:Gap eigenvalues from nestedness}
    Let $C$ be a $d$-cone and let $g\in \GL_n(\C)$ be such that $\overline{g\cdot C}\subset C$. If $(\lambda_i)$ are the characteristic values of $g$ with $|\lambda_1|\geq |\lambda_2|\geq \cdots \geq |\lambda_n|$:
    $$\log|\lambda_d|-\log|\lambda_{d+1}|\geq {\log\Nest}(C,g\cdot C).$$
\end{prop}

One can also control the singular value gap of an element nesting a cone. In the proof of our main theorem we will use this this precise control of singular values, but the simple control of eigenvalue would suffice to just show the Anosov property.

First note the following geometric characterization of roots of Cartan projections. We only need this general statement to control precisely the singular value gaps for the exceptional groups.

\begin{prop}
    \label{prop:Genreal fact Singular values}
    Let $G$ be a simple complex group and $\mathcal{F}$ a flag manifold of $G$, associated to a simple root $\alpha$. Let us fix a point $o\in\X$ with stabilizer $K$. Fix the $K$-invariant Riemannian metric on $\mathcal{F}$ induced by the Killing form and let $g\in G$. For all $x\in \mathcal{F}$:
    $$\alpha(\mu(o,g\cdot o))\geq-\log(\lVert d_x g \rVert).$$
    In this expression $\lVert d_xg\rVert=\lVert d_xg\rVert_\text{op}$ is the operator norm for the metric on $\mathcal{F}$, defined as the supremum of $\lVert d_xg(v)\rVert/\lVert v\rVert$ for all non-zero $v\in T_x\mathcal{F}$.
\end{prop}

This inequality is an equality for some $x\in \mathcal{F}$ \cite[{Theorem 4.7}]{GGKW2017}. This reference also shows that the inequality holds at $x\in \mathcal{F}$ whenever $g$ belongs to the Levi subgroup associated to $x$ and its opposite flag $x^\text{opp}$ (where $x^\text{opp}$ is taken with respect to the maximum compact $K$), and in particular as a corollary the inequality also holds whenever $g\cdot x$ and $g\cdot x^\text{opp}$ are opposite flags with respect to $K$. We now check its holds for any $g\in G$ and $x\in \mathcal{F}$.

\begin{proof}
We equip $\mathfrak{g}$ with the hermitian pairing $\langle\cdot, \cdot \rangle =B(\cdot,- \theta \cdot)$ induced by the Killing form and the Cartan involution associated to $K$. Let us write $g=k_1 e^hk_2$ with $h=\mu(o,g\cdot o)\in \mathfrak{a}^+$ and $k_1,k_2\in K$. The action of $k_1,k_2$ does not change the operator norm, so we consider only the action of $e^h$ on $\mathcal{F}$. 

Let $x\in \mathcal{F}$. We need to a find a tangent vector to the flag manifold that is not too contracted. Let $X\perp\mathfrak{p}_x$ be non-zero vector in the corresponding parabolic subalgebra, corresponding to a vector $[X]\in T_x\mathcal{F}\simeq \mathfrak{g}/\mathfrak{p}_x$, and let us set $Y=\Ad(e^{-h})X$\footnote{Here $Y$ can be viewed as the image of $X$ by the action of $e^h$ on the dual $\mathfrak{g}^*$, which is identified with $\mathfrak{g}$ via the Killing form.}. For $Z\in\mathfrak{p}_x$, $\langle Y,\Ad(e^h)Z\rangle=\langle X,Z\rangle=0$, so $Y\perp\Ad(e^h)\cdot\mathfrak{p}_x$. Therefore, if we consider the norm of the image of $[X]$, viewed as a tangent vector to $T_{e^h\cdot x}\mathcal{F}$:
$$\lVert d_xe^h([X])\rVert=\lVert[\Ad(e^h)X]\rVert\geq\frac{|\langle\Ad(e^h)X,Y\rangle|}{\lVert Y\rVert}=\frac{\lVert X\rVert^2}{\lVert Y\rVert}.$$

Indeed the norm of an element $[\Ad(e^h)X]\in \mathfrak{g}/(\Ad(e^h)\cdot\mathfrak{p}_x)$ is the norm of its projection to $(\Ad(e^h)\cdot\mathfrak{p}_x)^\perp$, which is at least the norm of its orthogonal projection to the line spanned by $Y$ in $(\Ad(e^h)\cdot\mathfrak{p}_x)^\perp$. Since we took $X\in \mathfrak{p}_x^\perp$, one has $\lVert X\rVert=\lVert [X]\rVert$ and hence:
$$\lVert d_xe^h\rVert\geq\frac{\lVert X\rVert}{\lVert Y\rVert}.$$

Let $\mathfrak{q}=\mathfrak{h}\oplus\g_{-\alpha}\oplus\bigoplus_{\beta\in \Sigma^+}\g_\beta\subset \mathfrak{g}$ where $\mathfrak{h}$ is the Cartan subalgebra, and $\Sigma^+$ the set of positive roots. The operator $\Ad(e^{-h})$ acts trivially on $\mathfrak{h}$, by $e^{-\beta(h)}\leq1$ on the $\g_\beta$  with $beta\in \Sigma^+$ because $h\in\mathfrak{a}^+$, and by $e^{\alpha(h)}$ on $\g_{-\alpha}$. Hence $\lVert Y\rVert=\lVert\Ad(e^{-h})X\rVert\leq e^{\alpha(h)}\lVert X\rVert$ whenever $X\in\mathfrak{q}$, and then $\lVert d_xe^h\rVert\geq e^{-\alpha(h)}$. It is therefore enough to find a non-zero $X\in\mathfrak{p}_x^\perp\cap\mathfrak{q}$ to conclude.

\medskip

We now check that $\mathfrak{p}_x^\perp\cap\mathfrak{q}\neq\lbrace0\rbrace$. For $x\in\mathcal{F}$, let $x^{\mathrm{opp}}$ be the point of the opposite flag manifold $\mathcal{F}^{\mathrm{opp}}$ whose stabilizer has Lie algebra $\theta(\mathfrak{p}_x)$, and let $\mathfrak{u}_{x^{\mathrm{opp}}}$ be the nilradical of $\theta(\mathfrak{p}_x)$. Then $\mathfrak{p}_x^\perp=\mathfrak{u}_{x^{\mathrm{opp}}}$.

It is therefore enough to show that $\mathfrak{u}_{x^\text{opp}}\cap\mathfrak{q}\neq\lbrace0\rbrace$ for every $x^\text{opp}\in\mathcal{F}^{\mathrm{opp}}$. If $x^\text{opp}$ is transverse to $x_0$, then $x^\text{opp}=u\cdot x_0^{\mathrm{opp}}$ for some $u$ in the unipotent radical $U^+$ of the stabilizer of $x_0$, so $\mathfrak{u}_{x^\text{opp}}=\Ad(u)\cdot \mathfrak{u}_{x^\text{opp}_0}$. Recall that the Lie algebra of $U^+$ is the sum of root spaces $\mathfrak{g}_\beta$ for all positive roots $\beta\in \Sigma^+$ having at least one component in $\alpha$, and in particular $U^+$ preserves $\mathfrak{q}$. On the other hand $\mathfrak{u}_{x_O}^\text{opp}$ is the sum of root spaces $\mathfrak{g}_{-\beta}$ where $\beta\in \Sigma^+$ has at least one component in $\alpha$.

Hence $\Ad(u)\cdot \mathfrak{g}_{-\alpha}$ belongs to the intersection of $\mathfrak{q}$ and $\mathfrak{u}_{x^\text{opp}}=\mathfrak{p}_x^\perp$. In general, the points of $\mathcal{F}$ transverse to $x_0$ are dense, so the inequality extends to all points $x\in \mathcal{F}$.
\end{proof}

Let us also add the following characterization of the operator norm of the action of an element of $\GL(V)$ on $\Gr_d(V)$.

\begin{lem}
    \label{lem:Operator norm}
    Let $E\oplus F=\C^n$ be a hermitian decompostion for a fixed hermitian metric and let $g\in \GL_n(\C)$. Let $\pi$ be the orthogonal projection to the orthogonal for the hermitian metric to $g\cdot E$. The operator norm $\lVert d_Eg\rVert$ for the $\mathrm{U}(n)$-invariant Riemanian metric on $\Gr_d(\C^n)$ induced by the Killing form is equal to the maximum of the following ratio, for all non-zero $e\in E$ and $f\in F$:
$$\frac{\lVert \pi(g(f))\rVert}{\lVert g(e)\rVert}\frac{\lVert e\rVert}{\lVert f \rVert}.$$
\end{lem}

\begin{proof}
    A tangent vector $[\mathrm{v}]\in T_E\Gr_d(\C^n)$ is a linear map $[\mathrm{v}]:E\to \C^n/E$. It is mapped by $g$ to a linear map $[g\cdot\mathrm{v}]:g\cdot E\to \C^n/g\cdot E$. For the standard metric on $\Gr_d(\C^d)$ the norm of $[\mathrm{v}]$ is the $\ell^2$-norm of its representative $\mathrm{v}:E\to F$, since $E\perp F$. The norm of $[g\cdot \mathrm{v}]$ is the $\ell^2$-norm $\lVert g\cdot \mathrm{v}\rVert_2$ of its representative $g\cdot\mathrm{v}: g(e)\mapsto \pi(g(\mathrm{v}(e)))$.

    \medskip

One has  $g\cdot\mathrm{v}=B\circ \mathrm{v} \circ A^{-1}$ where $A:E\to g\cdot E$ is the action of $g$ and $B:F\to (g\cdot E)^\perp$ is equal to $\pi\circ g$. The ratio between $\ell^2$-norm of $g\cdot\mathrm{v}$ and $\mathrm{v}$ is bounded from above by the product of the operator norms of $A^{-1}$ and $B$. Indeed for two square matrices $M,N$ one has $\lVert MN\rVert_2,\lVert NM\rVert_2\leq \lVert M\rVert_2\lVert N\rVert_{\text{op}}$. The operator norms of $A^{-1}$ and $B$ are respectively the maximum of $\frac{\lVert e\rVert}{\lVert g(e)\rVert}$ and $\frac{\lVert \pi(g(f))\rVert}{\lVert f\rVert}$.
\end{proof}

We relate the nestedness with gaps in singular values, and more generally in terms of contraction of flag manifolds.

\begin{prop}
\label{prop:Gap singular values from nestedness}
    Let $C$ be a $d$-cone, and let $g\in \GL_n(\C)$ be such that $\overline{g\cdot C}\subset C$. Then if $(\sigma_i)$ are the singular values of $g$ with $\sigma_1\geq \sigma_2\geq \cdots \geq \sigma_n$, with respect to a metric compatible with $C$:
    $$\log\sigma_d-\log\sigma_{d+1}\geq {\log\Nest}(C,g\cdot C).$$

Moreover we show that for the metric on $\Gr_d(V)$ induced by a compatible metric $h$ on $V$, $\lVert d_E g \rVert\leq {\Nest}(C,g\cdot C)^{-1}$, where $g$ is viewed as a diffeomorphism of $\Gr_d(V)$.
\end{prop}

\begin{proof}
The minimax characterization of singular values tells us that the ratio $\sigma_d/\sigma_{d+1}$ is characterized as the maximum for every $d$-dimensional subspace $E\subset \C^n$ and every supplement space $F$ of the following ratio for $e\in E$ and $f\in F$ non-zero:
$$\frac{\lVert g(f)\rVert}{\lVert g(e)\rVert}\frac{\lVert e\rVert}{\lVert f \rVert}.$$

Since the orthogonal projection reduce the norm, it suffices to bound from below the ratio from Lemma \ref{lem:Operator norm} by the nestedness, and hence it suffices to prove the moreover part.

Let us take $E$, $F$ respectively the span of the first $d$ basis vector and the last $n-d$ basis vectors, in a basis such that: 
$$C=\lbrace [x_1,\cdots,x_n]\in \mathbb{CP}^{n-1} \mid |x_1|^2+\cdots +|x_d|^2>|x_{d+1}|^2+\cdots +|x_{n}|^2\rbrace.$$
Let $q$ be the hermitian form of signature $(d,n-d)$ such that for $x=(x_1,\cdots, x_n)\in \C^n$: $$q(x):=q(x,x)=|x_1|^2+\cdots +|x_d|^2-(|x_{d+1}|^2+\cdots +|x_{n}|^2).$$
Let $E'=g\cdot E$ and $F'$ be the orthogonal with respect to $E'$ for the hermitian form $q$.

 Let $e\in E$ and $f\in F$ be non-zero vectors, without any loss of generality we will assume that they both have norm $1$. We have $g(e)=e'\in E'$ and $g(f)=f'+\eta$ where $f'\in F'$ and $\eta\in E'$. Note that this is is not the projection $\pi$ : the image of the projection is different but both projections are parallel to $E'$. Hence $f'$ is not in general equal to $\pi(g(f))$, but $\pi(f')=\pi(g(f))$. We want to bound from above the ratio $\frac{\lVert \pi(g(f))\rVert}{\lVert g(e)\rVert}=\frac{\lVert \pi(f')\rVert}{\lVert e'\rVert}$. Now let us write $\xi=g^{-1}(\eta)$. Up to multiplying $e$ by a unit complex number, we have $\text{Re}\langle e,\xi\rangle=0$, which will be convenient for us. Let us consider the real projective line $L$ spanned by $e'$ and $f'$.

The intersection $\partial C \cap L$ in this case is the pair: $(x_1,y_1)=(ae'+bf', ae'-bf')$, where:
$$a=\frac{1}{\sqrt{q(e')}}\geq \frac{1}{\lVert e' \rVert },$$ 
$$b=\frac{1}{\sqrt{-q(f')}}\leq \frac{1}{\lVert \pi(f')\rVert}.$$
For the second inequality, we need to introduce yet another projection: write $f'=\eta'+f''$ with $\eta'\in E'$ and $f''\in F$, coming from the direct sum $E'\oplus F=\C^n$ (this sum is direct since $\mathbb{P}(E')=g\cdot \mathbb{P}(E)\subset C$, $\mathbb{P}(F)\cap C=\emptyset$). One has $q(f',\eta')=0$ since $f'$ is already projected onto the $q$-orthogonal of $E'$. The norm of the orthogonal projection $\lVert \pi(f')\rVert$ is bounded above by $\lVert f'-\eta'\rVert$ since $\eta'\in E'$. But $f'-\eta'\in F$ so:
\[\lVert f'-\eta'\rVert=\sqrt{-q(f'-\eta')}=\sqrt{-q(f')+q(\eta')}\leq \sqrt{-q(f')}\,.\]
Here we used that $q(f',\eta')=0$ and $\eta'\in E'$ so $q(\eta')>0$ (this last point comes from the fact that $\mathbb{P}(E')=g\cdot \mathbb{P}(E)\subset g\cdot C\subset C$, so $q>0$ on $E'$.

The intersection $g(\partial{C})\cap L$ is then the pair $(x_2,y_2)=(e'+c f', e'-cf')$, where:
$$c=\frac{1}{\sqrt{1-\lVert \xi\rVert^2}}.$$

Indeed $c$ and $-c$ are the zeroes of the following expressions in $t$: 
$$q(g^{-1}(e'+tf'))=q(e+t(f-g^{-1}(\eta)))=1-t^2(1-\lVert \xi\rVert^2).$$ 
Here we use the fact that $\text{Re}\langle e,\xi\rangle=0$.

\medskip

 The points $(x_1,x_2,e',y_2,y_1,f')$ are in this order on $L$: $(y_1,y_2)$ and $(x_2,x_1)$ are in opposite components of $L\setminus \lbrace e',f'\rbrace$ and the fact that $\overline{g(C)}\subset C$ implies that $c< b/a$, describing the orderings of these pairs on the arcs of $L\setminus \lbrace e',f'\rbrace$.
 
 Hence the $4$ points of interest they are in the following cyclic order $(x_1,x_2,y_2,y_1)$. 
 We moreover have that:
${\Nest}(C,g\cdot C)\leq {\Nest}_L(C,g\cdot C)=  [x_1,x_2,y_2,y_1]$. Recall $[x_1,x_2,y_2,y_1]=\frac{Y_2-X_1}{Y_2-X_2}\frac{Y_1-X_2}{Y_1-X_1}$. By working in an affine chart where $e'+Xf'$ corresponds to $X\in \R$ we obtain:
$$[x_1,x_2,y_2,y_1]=\frac{-c-b/a}{-2c}\times\frac{-b/a-c}{-b/a-b/a}=\frac{(c+b/a)^2}{4c(b/a)}\leq \frac{b}{a}\leq \frac{\lVert e'\rVert}{\lVert \pi(f')\rVert}.$$ 
In the first inequality we use that $1\leq c\leq b/a$. Hence $\Nest(C, g\cdot C)\leq \lVert d_Eg\rVert^{-1}$ by Lemma \ref{lem:Operator norm}. This concludes the proof. 
\end{proof}

\subsection{Precise estimate}
In this subsection we study in more details how infinitesimal nestedness of cones provides a bound on the nesting constant of cones.

\begin{defn}
\label{defn:Infinitesimal Nesting Cone Precise}
    Let $C$ be a $d$-cone which is the set of positive lines for a pseudo-hermitian pairing $q$ and let $h$ be a compatible hermitian metric and let $\epsilon\geq0$. We say that an element $g\in \mathfrak{gl}_n(\C)$ is \emph{$\epsilon$-nesting $C$ with respect to $h$} if for all $v\in \partial C$:
    $$q( g\cdot v, v)\geq \epsilon h(v,v).$$
\end{defn}

 We first verify that this local nestingness implies that globally the cones are nested with an explicit bound.

\begin{prop}
\label{prop:Differential inequation nesting}
 Let $(g_t)$ be a smooth curve in $\GL_n(\C)$ with $g_0=\Id$. Let $C_0$ be a $d$-cone and $h_0$ a compatible metric. Define $C_t=g_t\cdot C_0$ and $h_t=g_t\cdot h_0$. Suppose that for all $t\geq 0$ the derivative $g'_t=\frac{\mathrm{d}}{\mathrm{d}s}_{\mid s=t}g_sg_t^{-1}\in \mathfrak{gl}_n(\C)$ is $\epsilon_t$-nesting the cone $C_t$ with respect to $h_t$, for some $\epsilon_t\geq 0$ varying continuously. Then $F(t)=\text{Nest}(C_0,C_{t})$ satisfies the differential inequality:
 $$F'(t)\geq\epsilon_t\sqrt{\left(F(t)-1/2\right)^2-1/4}.$$
 \end{prop}

 \begin{proof}
Let us fix $t\geq 0$ and a real projective line $L$ that crosses $(C_0,C_{t})$. Let $x_0,y_0\in \partial C_0\cap L$ and $x_t,y_t\in \partial C_{t}\cap L$ so that $(x_0,x_t,x_{t'},y_{t'},y_t,y_0)$ are cyclically in this order on $L$ for all $0<t<t'$. We fix an affine chart on $L$ so that $x_t$ corresponds to $X_t=-1$, $y_t$ corresponds to $Y_t=1$, and in general $x_{t'},y_{t'}$ correspond to $X_{t'},Y_{t'}$. We have:
$$\text{Nest}_L(C_0,C_{t})=\frac{(Y_0-X_t)(Y_t-X_0)}{(Y_t-X_t)(Y_0-X_0)}=\frac{(Y_0+1)(1-X_0)}{2(Y_0-X_0)}=\frac{1-X_0Y_0}{2(Y_0-X_0)}+\frac{1}{2}.$$
We write for simplicity $A=\frac{1-X_0Y_0}{2(Y_0-X_0)}$.

\medskip

On $L$ the form $q_t$ defining $C_t$ writes as $X\mapsto 1-X^2$, (up to a multiplicative constant that can be assume to be $1$ up to choosing an appropriate changing the affine chart). The compatible metric $h_t$ writes as $aX^2+bX+c$ for some $a,b,c\in \R$. The compatibility implies that $aX^2+bX+c\geq |X^2-1|$. Up to changing the affine chart while fixing $-1$ and $1$ we can assume that $b=0$. Then $c\geq 1$ and $a\geq 1$ so $h(X)\geq X^2+1$ in this basis.

\medskip

 Fix $t\leq s$. Note that since $y_s\in \partial C_s$ and $C_s=g_sg_t^{-1}C_t$, the derivative ${\frac{\mathrm{d}}{\mathrm{d}s}y_s}_{\mid s=t}$ decomposes as $g'_t\cdot y_t+v$ where $v\in T_{y_t}C_t$. Let $\hat{y_s}\in \C^n$ be the vector representing $y_s$ in the affine chart of $L$, then $q_t(\frac{\mathrm{d}}{\mathrm{d}s}\hat{y}_s,\hat{y}_t)=q_t(g'_t\cdot \hat{y}_t,\hat{y}_t)$ as the contribution of $v$ vanishes. Since $g'_t$ is $\epsilon_t$-nesting $C_t$ with respect to $h_t$:
$$\frac{\mathrm{d}}{\mathrm{d}s}_{\mid s=t}q_t(\hat{y}_s,\hat{y}_s)=2q_t(g'_t\cdot \hat{y}_t,\hat{y}_t)\geq 2\epsilon_th_t(\hat{y}_t,\hat{y}_t)\geq 2\epsilon_t.$$
Hence for all $s<t$, $1-Y_s^2=q_t(\hat{y}_s,\hat{y}_s)\leq-2\epsilon_t (s-t)+o(s-t)$.
Since the derivative of $X\mapsto 1-X^2$ at $X=1$ is $-2$, the intersection point $Y_{t+\delta}$ of $\partial C_{t+\delta}$ with $L$ close to $Y_t$ satisfies, for $\delta>0$ small enough, $Y_{t+\delta}\leq 1-\delta \epsilon_t+o(\delta)$. Similarly we have:
$$X_{t+\delta}\geq -1+\delta \epsilon_t+o(\delta).$$

As long as $X\leq X_0<Y_0\leq Y$, the quantity $\frac{(Y_0-X)(Y-X_0)}{(Y-X)(Y_0-X_0)}$ is increasing in $X$ and decreasing in $Y$. Hence for $\delta$ small enough:
\begin{align*}
    \text{Nest}_L(C_0,C_{t+\delta})\geq \frac{(Y_0+1-\delta \epsilon_t)(1-X_0-\delta \epsilon_t)}{(2-2\delta \epsilon_t)(Y_0-X_0)}+o(\delta)\\
=\text{Nest}_L(C_0,C_{t}) + \epsilon_t \delta B+o(\delta).
\end{align*}

Here $B=\frac{(1-X_0)(Y_0+1)-(1-X_0)-(Y_0+1)}{2(Y_0-X_0)}=-\frac{1+X_0Y_0}{2(Y_0-X_0)}$. This quantity is positive since $X_0<-1$ and $Y_0>1$.

We compute $A^2-B^2=\frac{-4X_0Y_0}{4(Y_0-X_0)^2}$.

We conclude that $A^2-B^2\leq \frac{1}{4}$. Indeed $(X_0-Y_0)^2= X_0^2+Y_0^2-2X_0Y_0\geq -4X_0Y_0$. This concludes the proof as $F'(t)\geq \epsilon_t B\geq \epsilon_t\sqrt{A^2-1/4}$ and $A=F(t)-\frac{1}{2}$.

\end{proof}

\begin{prop}
\label{prop:Nestedness Infinitesimal to global simple}
     Under the assumptions of Proposition \ref{prop:Differential inequation nesting},
     for all $t>0$:
\begin{equation}
    \Nest(C_0,C_t)\geq \big(\cosh(\frac{1}{2}\int_{s=0}^t\epsilon_s)\big)^2\,.
\end{equation}
 In particular, for all $t>0$,
 \begin{equation}
     {\log\Nest}(C_0,C_t)\geq \int_{s=0}^t\epsilon_s-2\log 2\,.
 \end{equation}
 
\end{prop}
\begin{proof}
    Let $F$ be the nesting
    as before.
    First, for $t>0$ small enough, the strict nesting
    yields that $F>1$.
    As such, we integrate the inequality
    \begin{equation*}
    \frac{F'(t)}{\sqrt{F(t)(F(t)-1})}\geq\epsilon_t\,,
    \end{equation*}
    to obtain
    \begin{equation*}
    2\text{arccosh}(\sqrt{F(t)})\geq \int_0^{t}\epsilon_sds\,.
    \end{equation*}
    This ensures that $F(t)>1$ for any $t>0$, and so the inequality holds for all $t>0$.
    Equivalently, this shows that
    \begin{equation*}
        F(t)\geq\big(\cosh(\frac{1}{2}\int_{s=0}^t\epsilon_s)\big)^2\,,
    \end{equation*}
    as desired.
\end{proof}

\subsection{Anosov representations and undistorted surfaces.}

We now define the Anosov property, giving the characterization from \cite{GW12,KLP17,BPS2019}. Let $\Theta$ be a set of simple roots associated with the Lie group $G$. Let $\Gamma$ be a hyperbolic group, that we equip with a word metric $|\cdot|_\Gamma$. Recall that $\mu:\X\times\X\to\mathfrak{a}^+$ denotes the Cartan projection. 

\begin{defn}
    A representation $\rho:\Gamma\to G$ is \emph{$\Theta$-Anosov} if for one (equivalently, any) point $x\in \mathbb X$ there exist $A,B>0$ such that for all $\gamma\in\Gamma$ and $\alpha\in\Theta$:
    $$\alpha(\mu(x,\rho(\gamma)\cdot x))\geq A|\gamma|_\Gamma-B.$$
\end{defn}

When we consider surfaces inside the symmetric space $\mathbb{X}$ that are not equivariant with respect to a representation, we will use the following notion of $\Theta$-undistorted maps.

\begin{defn}  
We say that a Lipschitz map $u:E\to \mathbb{X}$ from a metric space $(E,d)$ is $\Theta$-undistorted if there exist $A,B>0$ such that for all $\alpha\in \Theta$:
$$\alpha(\mu(u(x),u(y)))\geq A\cdot d(x,y)-B.$$
\end{defn}

Any Lipschitz map $u:E\to \mathbb{X}$ that is equivariant with respect to a representation $\rho:\Gamma\to G$, for a cocompact and proper action of $\Gamma$ on $E$ is $\Theta$-undistorted if and only if the representation $\rho$ is $\Theta$-Anosov.  

\subsection{Proof of the main theorems.}

Let $X$ be a not necessarily compact Riemann surface. Let $G$ be a complex simple group of finite center and let $\Theta$ be the set of simple roots that do not vanish on $h$, described in Figure \ref{table:Slodowy slices}, which contains only one root except if $G$ is of type $A_n$. Recall that $\mu:\X\times \X\to \mathfrak{a}^+$ is the Cartan projection and $\lambda:G\to \mathfrak{a}^+$ the Jordan projection into the model Weyl chamber of the restricted (real) Cartan subalgebra.

 \begin{thm}
\label{thm:Unsdistorted On the Unit Disk}
 Let $G$ be a complex simple Lie group with finite center and of type different from $E_8$. Let $(\mathcal{E},\varphi,H)$ be a parallelisable harmonic bundle over $\mathbb{D}$ in the restricted Slodowy slice of the minimal $\mathfrak{sl}_2$-triple with complete projected metric $g_\pi$ and bounded differential $\mathbf{r}$. 
 
 Consider the corresponding minimal immersion $\harm:\mathbb{D}\to \X$ into the symmetric space induced by the harmonic metric $H$ and the flat connection $\nabla=\nabla^H+\varphi+\varphi^{*H}$. Let $\Theta\subset \Sigma$ be the set of simple roots that do not vanish on $h$, and let $-\epsilon$ be the maximum of the curvature of the projected metric $g_\pi$, with $\epsilon>0$. 

 Then $h$ is $\Theta$-undistorted and embedded. More precisely, 
 for all $x,y\in \mathbb{D}$ and all $\alpha\in \Theta$:
$$\alpha(\mu(\harm(x),\harm(y)))\geq \frac{\epsilon}{12} d_\pi(x,y)-B.$$
Here $d_\pi$ denotes the distance associated to the projected metric on $\mathbb{D}$, and $B>0$ is a constant depending only on $G$. 

Let $g\in G$ be the realization of the parallel transport in $\X$ along the image by $\harm$ of the geodesic segment for the projected metric joining $x$ and $y$:
$$\alpha(\lambda(g))\geq \frac{\epsilon}{12} d_\pi(x,y).$$
 \end{thm}

 \begin{rem}
     In this theorem, the metric we place on $\X$ has an influence on the constant in front of the projected metric. Here we keep the conventions for the normalisation of the projected metric from the section on $3$-cyclic Higgs bundles in Section \ref{subsec:3-cyclic}. In particular for this metric the totally geodesic copies of $\H^2$ in $\X$ associated to the minimal $\mathfrak{sl}_2$ triple have curvature $-4$.
 \end{rem}

\begin{proof}
    We can always lift the harmonic bundle to a $G'$-Higgs bundle where $G'$ is any finite extension of $G$. Up to this operation, that does not change the values of $\mu$ and of the Jordan projection $\lambda$, we take a cover of $G$ so that if $V$ is the smallest linear representation of the Lie algebra $\mathfrak{g}$ of $G$, it extends to a representation $G\to \SL(V)$. Let $V=V_1\oplus V_0\oplus V_{-1}$ be the weight decomposition for the longest coroot $h$ in $\mathfrak{g}$, and let $d=\dim(V_1)$.

    The fact that the differential $\mathbf{r}$ is bounded and the projected metric is complete implies that the curvature $K_{g_\pi}$ of the projected metric satisfies $K_{g_\pi}\leq -\epsilon$ for some $\epsilon>0$, by Theorem \ref{thm:NormBound3-cyclic}.

    \medskip

    We consider the minimal immersion $\harm:\mathbb{D}\to \X$ associated to the flat connection of the harmonic bundle and the harmonic metric $H$. Let $\mathscr{V}:T\mathbb{D}\to T\X$ be the map from Section \ref{subsec:parralel distribution of planes}, defining the parallel distribution of planes. Let us fix two points $x,y\in \mathbb{D}$ and consider the geodesic $\gamma:[0,d_\pi(x,y)]\to \mathbb{D}$ for the projected metric $g_\pi$ between them. Let $f\in G$ be the element realizing the parallel transport along $\harm\circ \gamma$ in $\X$. We define for each $v\in T\mathbb{D}$ a $d$-cone $C_v\subset \mathbb{P}(V)$, as in \eqref{eq:Defn Cv,Hv} in Section \ref{subsec:Nestedness}.
    Corollary \ref{cor:ParallelTransport Distribution planes} implies that $f(C_{v_0})=C_{v_1}$, where $v_0=\gamma'(0)$, $v_1=\gamma'(d_\pi(x,y))$. Moreover $\mu(\harm(x),\harm(y))=\mu(\harm(x),f\cdot\harm(x))$. 

    We now apply Proposition \ref{prop:Nestedness Infinitesimal to global simple}, with the computations of the infinitesimal nesting from Proposition \ref{prop:3cycliclowerbound}. Together with the fact that $K_{g_\pi}\leq -\epsilon$, we obtain that ${\log\Nest}(C_{v_0},C_{v_1})\geq \frac{\epsilon}{12} d_\pi(x,y)-2\log(2)$.

    Hence it only remains to control $\alpha(\mu(\harm(x),f\cdot \harm(x)))$ and $\alpha(\lambda(f))$ for $\alpha \in \Theta$ to obtain the desired result, since $\harm(y)=f\cdot \harm(x)$. For the classical groups, one can directly link the gap in singular values with the evaluation of a root to the Cartan projection. Hence we immediately get the following inequality from Proposition \ref{prop:Gap singular values from nestedness}, applied with the hermitian metric $\harm(x)$ which is compatible with $C_{v_0}$:\[\alpha(\mu(\harm(x),f\cdot \harm(x))),\alpha(\lambda(f))\geq \frac{\epsilon}{12} d_\pi(x,y)-2\log(2).\]

If $G$ is of type $A_n$ we consider only the root $\alpha_1$, and then by a similar argument looking at the dual action on $\C^{n+1}$ we obtain the result for the other root $\alpha_n$.

    \medskip

    For general simple complex groups, in particular for the exceptional cases, we are working with the smallest linear representation $V$, and so we now detail the link between the Cartan projections in $\SL(V)$ and $G$. Let $\mathcal{F}$ be the flag manifold associated to one of the simple $\Theta$ that do not vanish on the longest coroot $h$ if $G$ is not of type $A_n$. Let $V_1\subset V$ be the $1$-weight space for $h$ and let $P_0$ be its stabilizer in $G$. The lie algebra of $P_0$ is exactly the parabolic Lie algebra associated with $\Theta$. Hence the $G$-orbit of $V_1$ in $\Gr_d(V)$, with $d=\dim(V_1)$ is, up to finite cover, identified with the flag manifold $\mathcal{F}$. 
    
    Using the second part of Proposition \ref{prop:Gap singular values from nestedness} we found a point $E\in \Gr_d(V)$ in the orbit of $V_1$ such that at this point $\lVert\mathrm{d}_Ef \rVert\leq \Nest(C_{v_0},C_{v_1})^{-1}$ for the operator norm for the $\mathrm{U}(n)$-invariant metric. The operator norm on the full Grassmannian is greater than that on the smallest orbit, meaning that we have found $E\in\mathcal{F}$ such that $\lVert\mathrm{d}_Ef \rVert\leq e^{-{\log\Nest}(C_{v_0},C_{v_1})}$, with the operator norm for a $K$-invariant metric, where $K$ is the stabilizer of $\harm(x)$. This metric is comparable with the metric induced by the Killing form, so the operator norms does only change by a constant factor. Hence for some $C>0$ depending only on $G$, we use the criterion from Proposition \ref{prop:Genreal fact Singular values} to obtain:\footnote{For the classical groups, one can obtain the result with $C=0$ using directly Proposition \ref{prop:Gap singular values from nestedness}.}
    $$\alpha(\mu(\harm(x),f\cdot \harm(x)))\geq {\log\Nest}(C_{v_0},C_{v_1})-C\geq \frac{\epsilon}{12}d_\pi(x,y)-2\log(2)-C.$$
    For the Jordan projection, we use the fact that:
    $$\lambda(f)=\lim_{n\to +\infty}\frac{1}{n}\mu(\harm(x),f^n\cdot \harm(x)).$$
    Note that $f^n\cdot C_{v_0}\subset C_{v_0}$ and that by the additivity of the nestedness, from Proposition \ref{prop:Additivity Nestedness}, $\log\Nest(C_{v_0},f^n\cdot C_{v_0})\geq n\log\Nest(C_{v_0},C_{v_1})$. Therefore:
    $$\alpha(\lambda(f))=\lim_{n\to +\infty}\frac{1}{n}\alpha(\mu(\harm(x),f^n\cdot \harm(x)))\geq \lim_{n\to +\infty}\frac{\epsilon}{12}d_\pi(x,y)-\frac{1}{n}(2\log(2)+C),$$
    $$\alpha(\lambda(f))\geq \frac{\epsilon}{12}d_\pi(x,y).$$

This concludes the proof.
\end{proof}

We also obtain a better control of some weight applied to the Cartan and Jordan projections, without the hypothesis of having a bounded differential and without assuming that the universal cover is the disk. Let $\omega_V$ be the highest weight of the representation $V$.

\begin{thm}
\label{thm:Undistorted weight on the disk}
    Let $G$ be a complex simple Lie group with finite center and of type different from $E_8$. Let $(\mathcal{E},\varphi,H)$ be a parallelisable harmonic bundle over $X=\mathbb{D}$ or $X=\C$ in the restricted Slodowy slice of the minimal $\mathfrak{sl}_2$-triple with complete projected metric $g_\pi$ defining a distance $d_\pi$ on $X$. Let $x,y\in X$ :
    $$\omega_V(\mu(\harm(x),\harm(y)))\geq d_\pi(x,y).$$
    If $g\in G$ is the realization of the parallel transport from $\harm(x)$ to $\harm(y)$ along the image by $\harm$ of a geodesic for the projected metric, then :
    $$\omega_V(\lambda(g))\geq d_\pi(x,y).$$
\end{thm}

\begin{proof}
The proof is analogous to that of Theorem \ref{thm:Unsdistorted On the Unit Disk}. We define $\gamma:[0,d_\pi(x,y)]\to X$ the corresponding geodesic in $X$ for the projected metric $g_\pi$. Let $v_0=\gamma'(0)$ and $v_1=\gamma'(d_\pi(x,y))$. Let $f\in G$ be the element realizing the parallel transport along $\harm\circ \gamma$ in $\X$. We define for each $v\in T\mathbb{D}$ a filled hyperboloid $\mathcal{H}_v\subset V$, as in \eqref{eq:Defn Cv,Hv} in Section \ref{subsec:Nestedness}. Corollary \ref{cor:ParallelTransport Distribution planes} implies that $f(\mathcal{H}_{v_0})=\mathcal{H}_{v_1}$, where $v_0=\gamma'(0)$, $v_1=\gamma'(d_\pi(x,y))$. Moreover $\mu(\harm(x),\harm(y))=\mu(\harm(x),f\cdot \harm(x))$. 

We now apply Proposition \ref{prop:HyperboloidNestingInfinitesimal}, with the computations of the infinitesimal nesting from Proposition \ref{prop:3cycliclowerbound} and obtain that ${\log\Nest}(\mathcal{H}_{v_0},\mathcal{H}_{v_1})\geq d_\pi(x,y)$.
    
We apply Proposition \ref{prop:HyperboloidNestingControl} to get that $\lambda_d(f)\geq d_\pi(x,y)$ where $\lambda_d$ is the logarithm of the modulus of the $d$-th largest eigenvalue and $d=\dim(V_1)$. This directly implies that $\omega_V(\lambda(f))\geq d_\pi(x,y)$, where $\omega_V$ is the highest weigth of the linear representation from Proposition \ref{prop:RepresentationV}. Since this is a fundamental weight, one has the general inequality $\omega_V(\mu(\harm(x),f\cdot \harm(x)))\geq \omega_V(\lambda(f))$, which concludes the proof.
\end{proof}

The previous result can be interpreted as a domination with respect to the Fuchsian case.
Let $\harm_0:\mathbb{D}\to \X$ be the totally geodesic conformal embedding associated with any harmonic bundle $(\mathcal{E}_T,\mathbf{e},H)$ whose projected metric is the hyperbolic metric on $\mathbb{D}$. Let $f\in G$ be the realization of the parallel transport from $\harm_0(x)$ to $\harm_0(y)$ along the geodesic for the projected metric, which is the hyperbolic metric, with distance $d_\text{hyp}$:
\[\omega_V(\lambda(f))=d_\text{hyp}(x,y).\]

\medskip

In particular for closed Riemann surfaces, we obtain the following.

\begin{cor}
\label{cor:Anosov closed surface groups}
    Let $X$ be a compact Riemann surface of genus $g\geq 2$. Let $\rho:\pi_1(X)\to G$ be the representation associated to a Higgs bundle $(\mathcal{E}_T\star \mathcal{E}_C,\mathbf{e+r})$ in the minimal Slodowy slice with a $TC$-compatible harmonic metric. Let $\rho_0:\pi_1(X)\to G$ be the uniformizing representation, associated to $(\mathcal{E}_{T},\mathbf{e})$.
    Then the following properties are true:
    \begin{itemize}
        \item[](Anosov property) The representation $\rho$ is $\Theta$-Anosov.
        \item[](Domination)For all $\gamma\in \pi_1(X)$:
        $$\omega_V(\lambda(\rho(\gamma)))\geq\omega_V(\lambda(\rho_0(\gamma))).$$
        If moreover $\mathbf{r}\neq 0$ then for some $\epsilon>0$ for all $\gamma\in \pi_1S$:
        $$\omega_V(\lambda(\rho(\gamma)))\geq (1+\epsilon)\omega_V(\lambda(\rho_0(\gamma))).$$
    \end{itemize}
\end{cor}

Note that the harmonic metric for such a polystable Higgs pair is automatically compatible with the $TC$-reduction, as proven later in Corollary \ref{cor:SolutionHitchin equation compacible}.

\begin{proof}

The Anosov property is a direct consequence of Theorem \ref{thm:Unsdistorted On the Unit Disk}.
    \medskip

    The domination is a consequence of Theorem \ref{thm:Undistorted weight on the disk}. Note that the projected metric on a compact surface is greater than the hyperbolic metric by Theorem \ref{thm:NormBound3-cyclic}, which gives the domination. For the strict domination in the case with $\mathbf{r}\neq 0$ one can use the fact that $g_\pi\geq(1+\epsilon) g_\text{hyp}$ on the compact surface for some $\epsilon>0$. 
\end{proof}

We now show our results for the two additional slices we consider. Consider complete maximal spacelike surfaces $M$ in $\mathbb{H}^{2,n}$, which is the space of negative lines in $\R^{2,n+1}$, with curvature bounded from above by $-\epsilon<0$. By \cite[Proposition~3.10]{LTW24}, $M$ is an entire graph and hence is simply connected. The curvature bound implies that \(M\) is quasiperiodic; by \cite[Theorem~B]{LT}, its induced conformal structure is that of the unit disk $\mathbb D$. The following orthogonal splitting then yields a conformal $\mathrm{SO}_0(2,n+1)$-Higgs bundle on $\mathbb D$, equipped with a harmonic metric:
\[
\underline{\mathbb R}^{2,n+1}=TM\oplus\ell\oplus NM.
\]
Its Higgs field records, in particular, the $(1,0)$-part of the second fundamental form \cite[Sections 3.3--3.4]{CTT19}, see also \cite[Section 7.1]{DaiLiBounded} for this construction on the disk. The Gauss map of $M$ is the corresponding minimal immersion in the symmetric space $\X$ associated to this harmonic bundle, see  \cite[Section 3.3]{CTT19}. 

\begin{thm}
\label{thm:CompleteMaximalSurface}
    Let $M\subset \mathbb{H}^{2,n}$ be a complete maximal spacelike surface whose induced metric has sectional curvature lower than $-\epsilon$ for some fixed $\epsilon>0$. Let $\mathcal{G}:M\to \X$ the corresponding Gauss map into the Riemannian symmetric space of $\SO_o(2,n+1)$. This map is $\lbrace \alpha_1\rbrace$-undistorted. More precisely:
$$\alpha_1(\mu(\mathcal{G}(x),\mathcal{G}(y)))\geq 2\log\cosh (\frac{\epsilon}{8\sqrt{2}} d_\pi(x,y))\geq \frac{\epsilon}{4\sqrt{2}}d_\pi(x,y)-2\log(2).$$
\end{thm}
Here $d_\pi$ is the induced metric on the spacelike surface $M$, which is also the projected metric of the associated parallel distribution of planes. The convention we take to normalize this metric is so that when $M$ is a totally geodesic spacelike $\mathbb{H}^2$ in $\mathbb{H}^{2,n}$ its curvature is $-2$.
\begin{proof}
As in Section \ref{subsec:MaximalSO2N}, the harmonic Higgs bundle obtained from a maximal spacelike surface takes the following form:
\begin{equation}
\begin{tikzcd}
&\mathcal{V}\arrow[dl, "\gamma^{\dagger}"']&\\
\mathcal{IK}\arrow[r, "1"]&\mathcal{I}\arrow[r, "1"] & \mathcal{IK}^{-1}\arrow[lu, "\gamma"'] &
\end{tikzcd}.
\end{equation}

The associated metric is diagonal, and compatible with the orthogonal structure on the bundle.

\medskip

The proof of the $\lbrace \alpha_1\rbrace$-undistorted property is analogous to the proof of Theorem \ref{thm:Unsdistorted On the Unit Disk}. Indeed this Higgs bundle with its harmonic metric defines a parallelisable harmonic bundle. Hence we can also apply Corollary \ref{cor:ParallelTransport Distribution planes}. We apply the maximum principle from Lemma \ref{lem:NormBoundMaximalSO(2,n+1)}. We obtain nested $1$-cones from Proposition \ref{prop:Nestedness SO2n}, which in turn allows us to show the $\lbrace \alpha_1\rbrace$-undistorted property of the corresponding map in the symmetric space, by controlling the nestingness of the cones $C_v$.
\end{proof}

We finally prove Theorem \ref{mainthm:SL5}, which concerns the last slice of $\SL(n+4,\R)$-Higgs bundles that we study. Recall from Section~\ref{subsec:NilpotentSLn+4} that these Higgs bundles have the following form:
\begin{equation}\label{eq:SliceSL5}
\begin{tikzcd}
\mathcal{V}\arrow[r, "\nu"]&\mathcal{KL}\arrow[r, "1"] & \mathcal{L}\arrow[r, "\mu"] &\mathcal{L}^{-1}\arrow[r, "1"]&\mathcal{L}^{-1}\mathcal K^{-1}\arrow[bend right,llll, "\nu^{\dagger}"]&
\end{tikzcd}.
\end{equation}

\begin{thm}
\label{thm:SL5}
The representation $\rho:\pi_1S \to \SL(n+4,\R)$ associated to polystable Higgs bundles of the form \eqref{eq:SliceSL5} over a Riemann surface structure $X$ on $S$ are Anosov with respect to the stabilizer of a plane in $\R^{n+4}$.
\end{thm}

\begin{proof}
    The proof is completely analogous to that of Corollary \ref{cor:Anosov closed surface groups}, using the analog of Theorem \ref{thm:Unsdistorted On the Unit Disk}. Indeed this Higgs bundle with its harmonic metric defines a parallelisable harmonic bundle, since the harmonic metric is necessarily diagonal and compatible with the real structure. Hence we can also apply to it Corollary \ref{cor:ParallelTransport Distribution planes}. We apply the maximum principle from Lemma \ref{lem:NormBoundSO(n+4)}. We obtain nested $2$-cones from Proposition \ref{prop:Nestedness SLn+4}, which in turn allows us to show the $\lbrace \alpha_2\rbrace$-undistorted property of the corresponding map in the symmetric space. This in turn proves that the corresponding representations are Anosov with respect to the stabilizer of a plane.
\end{proof}

\subsection{Embeddedness of associated equivariant maps}

In this part we explain how the nesting of multicones
results in the fact that the harmonic maps $\harm:\widetilde X\rightarrow G/K$ and $u:\widetilde X\rightarrow G/G^\sigma$ are embedded in their respective homogeneous space, whenever one takes a parallelisable harmonic Higgs bundle in the minimal Slodowy slice with complete projected metric, when $G$ is not of type $E_8$.

\medskip

Recall in Section \ref{subsec:Nestedness} we showed how to associate to $v\in T^1\widetilde{X}$ a filled hyperboloid $\mathcal{H}_v\subset V$ where $V$ is the smallest linear representation of $\mathfrak{g}$.

\medskip

Then note the following. The subgroup $TC\cap K$ preserves a Hermitian metric $H$ on $V$ and the decomposition $V=V_1\oplus V_0\oplus V_{-1}$. Consider the antilinear pairing $H'=H_{|V_1}\oplus H_{|V_{-1}}\oplus(-H_{|V_0})$. This allows us to associate to every $[g]\in G/(TC\cap K)$ the corresponding antilinear pairing $H'_{[g]}=g\cdot H'$.

\begin{lem}
The hermitian pairing $H'_{[g]}=g\cdot H'$ on $V$ depends only on $[g]\in G/G^\sigma$.
\end{lem}

\begin{proof}
The subgroup $G^\sigma$ is generated by $TC\cap K$ and $(\mathfrak{g}_1\oplus \mathfrak{g}_{-1})\cap \mathfrak{p}$.
\end{proof}

\begin{thm}\label{thm:GaussEmbedded}
    Let $u:\widetilde X\rightarrow G/G^\sigma$ be the map associated to a parallelisable harmonic Higgs bundle over $X$ with complete projected metric. It is a smooth embedding.
\end{thm}
In particular the lift $Y:\widetilde X\to G/(TC\cap K)$ is also an embedding.

\begin{proof}

Take any geodesic $\gamma:[0,d_\pi(x,y)]\to \widetilde{X}$ for the projected metric between points $x,y$. As in the proof Theorem \ref{thm:Undistorted weight on the disk} one can associate to each $t\in [0,d_\pi(x,y)]$ a filled hyperboloid $\mathcal{H}_{-\gamma'(t)}$\footnote{The minus sign is due to the fact that we need here to look at filled hyperboloids that are getting bigger and bigger, instead of smaller and smaller.}. Note that the hyperboloid consists of vectors $z\in V$ such that $H'_{u(\gamma(t))}(z,z)\geq 1$, with equality for all $z\in V_1$ such that $H(z,z)=1$.  Moreover the nestedness estimates of these hyperboloids imply that:
$$\mathcal{H}_{-\gamma'(0)}\subset e^{t}\cdot \mathcal{H}_{-\gamma'(t)}.$$

Let us take some $x\in \partial \mathcal{H}_{-\gamma'(0)}$ such that $H'_{u(\gamma(0))}(z,z)= 1$. For $0\leq t\leq d_\pi(x,y)$, $H'_{u(\gamma(t))}(e^{-t}z,e^{-t}z)\geq 1$, meaning that $H'_{u(\gamma(t))}(z,z)\geq e^{2t}=e^{d(x, \gamma(t))}H'_{u(x)}(z,z)$. This holds for every geodesic $\gamma$, and implies that $u$ is an immersion and is injective. Since the projected metric is complete, $u$ is also proper when restricted to any geodesic. Moreover this properness on geodesics is uniform, as the vectors $z$ for which we measure $H'_{u(\gamma(t))}$ all live in a compact set. Hence $u$ is proper and is hence an embedding.
\end{proof}

\section{Slices in character varieties}
\label{sec:ModuliSpace}
Let $X$ be a compact Riemann surface of genus $g\ge 2$ with canonical bundle
$\mathcal{K}$. Let $G$ be a connected semisimple complex Lie group with Lie algebra
$\mathfrak g$.

\subsection{Twisted Higgs pairs and stability}
In this subsection, we recall the basic notion of Higgs pairs as a generalization of Higgs bundles. One may refer \cite{GG, GGM} for details.
Let $G_0$ denote a connected complex reductive Lie group with Lie algebra $\mathfrak g_0$ and $\rho:G_0\rightarrow \GL(V)$ a holomorphic representation. Let $L$ be a holomorphic line bundle over $X$.
\begin{defn}
A $L$-twisted $(G_0,V)$-Higgs pair is a pair $(\mathcal{E},\varphi)$, where $\mathcal{E}$ is a holomorphic principal $G_0$-bundle on $X$, and $\varphi\in H^0(X, \mathcal{E}[V]\otimes L)$, where $\mathcal{E}[V]:=\mathcal{E}\times_G V$ is the associated vector bundle to $\mathcal{E}$ via $\rho$.
\end{defn}  
When $L$ is $\mathcal K$, we omit to say ``$L$-twisted".

\begin{defn}
Let $G$ be a complex reductive Lie group, $G_1\leq G$ a Lie group, $\mathcal{E}$ a principal $G$-bundle. A holomorphic reduction of structure of $\mathcal{E}$ to $G_1$ is a holomorphic section $\sigma\in H^0(X,\mathcal{E}[G/G_1])$, where $\mathcal{E}[G/G_1]=\mathcal{E}\times_G G/G_1$. 

The natural map $\mathcal{E}\rightarrow  \mathcal{E}[G/G_1]$ has a $G_1$-bundle structure. A reduction $\sigma$ pulls back the $G$-bundle $\mathcal{E}$ to a $G_1$-bundle $\mathcal{E}_\sigma:=\sigma^*\mathcal{E}$ on $X$. Moreover, there is a canonical isomorphism $\mathcal{E}_\sigma\times_{G_1}G\cong \mathcal{E}$.
\end{defn}

Fix a maximal compact subgroup $K_0\subset G_0$. Let $\mathfrak k_0$ be its Lie algebra. For $s\in i\mathfrak k_0$, define the spaces
\[V_s^0=\{v\in V: \forall t\in \mathbb R, \rho(e^{ts})(v)=v\}, \]\[V_s=\{v\in V: \forall t\in \mathbb R, \text{$\rho(e^{ts})(v)$ is bounded as $t\rightarrow \infty$}\}.\]
The corresponding Levi and parabolic subgroups are:
\[L_s=\{g\in G_0:Ad(g)(s)=s\}, \quad P_s:=\{g\in G_0: \text{$e^{ts}ge^{-ts}$ is bounded as $t\rightarrow \infty$}\}.\]
The corresponding Lie algebras are $\mathfrak g_s^0:=\text{Lie}(L_s)$ and $\mathfrak g_s:=\text{Lie}(P_s)$. Define a character $\chi_s:\mathfrak g_s\rightarrow \mathbb C$ given by $\chi_s(x)=B(s,x)$ where $B$ is the Killing form on $\mathfrak g$. The subgroups $L_s, P_s$ acts on $V_s^0$, $V_s$ respectively via $\rho$.

Let $\mathcal{E}$ be a principal $G_0$-bundle and let $\sigma$ be a reduction of $\mathcal{E}$ to $P_s$. If a multiple $q\cdot \chi_s$ for some $q\in \mathbb Z_{\geq 0}$ lifts to a character $\tilde{\chi}_s:P_s\rightarrow \mathbb C^*$, we define the degree of the reduction as:
\[\deg \mathcal{E}(\sigma, s)=\deg_{\chi_s}(\mathcal E_\sigma):=\frac{1}{q} \deg(\mathcal{E}_\sigma\times_{\tilde{\chi}_s}\mathbb C^*).\] 

It is also possible to define this quantity when no multiple
of the character lifts to the group, using differential geometric
techniques. Given a reduction $\sigma$ of $E$ to $P_s$,
there is a further reduction $\sigma'$ to $K_s := K \cap L_s$,
the maximal compact of $L_s$. Let $A$ be a connection on
$E_{\sigma'}$ and consider its curvature
$F_A \in \Omega^2(X,E_{\sigma'}(\mathfrak{k}_s))$.
We have that $\chi_s(F_A) \in \Omega^2(X,i\mathbb{R})$,
and the degree is defined as
\[
\deg E(\sigma,s)
:= \frac{i}{2\pi}\int_X \chi_s(F_A).
\]

Given the representation $\rho:G_0\rightarrow \GL(V)$, let $d\rho:\mathfrak g_0\rightarrow \mathfrak{gl}(V)$ be its differential.
\begin{defn}
A $L$-twisted $(G_0,V)$-Higgs pair $(\mathcal E, \varphi)$ is:
\begin{itemize}
\item semistable, if for any element $s\in i\mathfrak k_0$ and reduction $\sigma\in H^0(X, \mathcal{E}[G_0/P_s])$ such that $\varphi\in H^0(\mathcal{E}_\sigma[V_s]\otimes L)$, we have $\deg \mathcal{E}(\sigma, s)\geq 0$.
\item stable, if it is semistable and, for any element $s\in i\mathfrak k_0\setminus \ker(d\rho)$ and reduction $\sigma\in H^0(\mathcal{E}[G_0/P_s])$ such that $\varphi\in H^0(\mathcal{E}_\sigma[V_s]\otimes L)$, we have $\deg \mathcal{E}(\sigma, s)>0$.
\item polystable, if it is semistable and, for any element $s\in i\mathfrak k_0$ and reduction $\sigma\in H^0( \mathcal{E}[G_0/P_s])$ such that $\varphi\in H^0(\mathcal{E}_\sigma[V_s]\otimes L)$ and  we have $\deg \mathcal{E}(\sigma, s)\geq 0$, there exists a reduction $\sigma'\in H^0(\mathcal{E}_\sigma[P_s/L_s])$ of $\mathcal{E}_\sigma$ to $L_s$ such that $\varphi\in H^0(\mathcal{E}_{\sigma'}[V_s^0]\otimes L)$. 
\end{itemize}
\end{defn}
\begin{defn}
The moduli space of polystable $L$-twisted $(G_0,V)$-Higgs pairs over $X$ parametrizes isomorphism classes of polystable $L$-twisted $(G_0,V)$-Higgs pairs. We denote it by $\mathcal M^{L}(G_0,V)$. If $L$ is $\mathcal K$, we will omit the superscript $L$.
\end{defn} 
By the construction of the moduli space via Geometric Invariant theory by Schmitt \cite{Schmitt}, we have $\mathcal M^L(G_0,V)$ is a quasi-projective complex variety. Fixing a topological type $\delta\in \pi_1(G_0)$ of $G_0$--bundles selects an open and
closed union of its components, written $\mathcal M_\delta^L(G_0,V)$. Let $K_0$ be a maximal compact subgroup of $G_0$. 

We finish this section with a definition.

\begin{defn}
Let $(\mathcal E,\varphi)$ be a $L$-twisted $(G_0,V)$-Higgs pair. Its automorphism
group is
\[
\operatorname{Aut}(\mathcal E,\varphi)
=
\left\{
g\in H^0\bigl(X,\operatorname{Ad}(\mathcal E)\bigr):
\rho(g)\varphi=\varphi
\right\}.
\]
The infinitesimal automorphism algebra is 
\[
\operatorname{aut}(\mathcal E,\varphi)
=
\left\{
u\in H^0\bigl(X,\mathcal E[\mathfrak g_0]\bigr):
d\rho(u)\varphi=0
\right\}.
\]
Here $\operatorname{Ad}(\mathcal E), \mathcal E(\g_0)$ are associated bundles of $\mathcal E$ via adjoint action of $G_0$ on $G_0$ and $\g_0$ respectively. 
\end{defn} 

\subsection{Harmonic metrics and Higgs pairs}
\begin{prop}(\cite[{Theorem 2.24}]{GGM})
\label{prop:Hitchin and polystability equivalent}
For a Higgs pair $(\mathcal{E},\varphi)$, it is polystable if and only if there exists a reduction of $\mathcal{E}$ to $K_0$ such that the Hitchin equation holds. Moreover, it is stable if and only if the reduction to $K_0$ is unique.
\end{prop}

Let $X$ be a compact connected Riemann surface and $G$ a connected
complex semisimple group. Let $s$ be an automorphism
of $G$ of finite order $m$. Set $\zeta=\exp(2\pi i/m)$ and write
\[
  \mathfrak g=\bigoplus_{j\in\mathbb Z/m\mathbb Z}\mathfrak g_{\bar j},
  \qquad ds|_{\mathfrak g_{\bar j}}=\zeta^j\,\mathrm{id},
  \qquad G_0=(G^s)^0.
\]
\begin{prop}(\cite[Proposition 5.7 (i)(ii)]{GPR})\label{prop:polystableEquiv}
Let $(\mathcal E,\varphi)$ be a $(G^s,\mathfrak g_{\bar 1})$-Higgs pair and $F$ denote its  structure extension map by $F(\mathcal E,\varphi)=(\mathcal E\times_{G_s}G, \varphi)$ as a $G$-Higgs bundle. Then $(\mathcal E,\varphi)$ is polystable if and only if $F(\mathcal E,\varphi)$ is polystable. 
\end{prop}

\begin{cor}
\label{cor:SolutionHitchin equation compacibleGeneral}
Let $(\mathcal{E},\varphi)$ be a polystable $(G_0,\mathfrak g_{\bar 1})$-Higgs pair. Then automatically any harmonic metric $H$ on this bundle is compatible with $G_0$.

The $L$-twisted Higgs pair~$(\mathcal E,\varphi)$ is polystable if and only if $F(\mathcal E,\varphi)$ is polystable. 

Moreover, if it is stable, then the harmonic metric is unique.
\end{cor}
\begin{proof}
The first statement follows from Proposition \ref{prop:Hitchin and polystability equivalent}.

For the second statement: 
Suppose $F(\mathcal E, \varphi)$ is polystable. Set 
\[K^s=K\cap G^s,\quad K_0=K\cap G_0.\] By  Proposition \ref{prop:polystableEquiv}, $(\mathcal E, \varphi)$ is polystable as a $(G^s, \mathfrak g_{\bar 1})$-Higgs pair and admits a harmonic metric $H$ compatible with $K_s$. Since the metric $H$ is a section of the associated bundle $\mathcal E\times_{G^s}(G^s/K^s)$ and 
\[G^s/K^s\cong G_0/K_0,\] it is also a section of the associated bundle $\mathcal E\times_{G^s}(G_0/K_0)\cong \mathcal E\times_{G_0}(G_0/K_0)$. Thus  $H$ is also compatible with $G_0$ and thus $(\mathcal E, \varphi)$ is polystable. Lastly, the uniqueness for stable Higgs bundles is also from Proposition~\ref{prop:Hitchin and polystability equivalent}. 
\end{proof}

A direct consequence in our case of interest will be Corollary \ref{cor:SolutionHitchin equation compacible}. 

\subsection{\texorpdfstring{$(TC,\langle e\rangle+R)$}{TC,<e>+R)}-Higgs pairs}
Recall the notions in Section \ref{sec:Slodowy}. We consider the minimal $\mathfrak{sl}_2$-triple $\lbrace e,h,f\rbrace$ in the Lie algebra $\mathfrak{g}$. Recall that the weights of $h$ define a (quaternionic) grading:
\[
\mathfrak g=\mathfrak g_{-2}\oplus\mathfrak g_{-1}\oplus\mathfrak g_0
          \oplus\mathfrak g_{1}\oplus\mathfrak g_{2} ,
\]

Let $\zeta\in\mathrm{Aut}(G)$ be the order--$3$ automorphism
inducing the associated $\mathbb Z/3$--grading, so that
$\langle e\rangle \oplus R=\mathfrak g_{-2}\oplus\mathfrak g_{1}$ is the $e^{2i\pi/3}$--eigenspace. Let $TC\subset G$ be the product subgroup of $T=\exp(\langle h\rangle)$ and the centralizer $C$ of the triple $\lbrace e,h,f\rbrace$. Recall that the Lie algebra of $TC$ is $\mathfrak g_0$. We consider 
$(TC,\langle e\rangle\oplus R)$--Higgs pairs $(\mathcal{E},\mathbf{e}+\mathbf{r})$, where
\[
\mathbf{e}\in H^0(\mathcal{E}[\mathfrak g_{-2}]\otimes \mathcal{K})=H^0(\mathcal{E}[\langle e\rangle]\otimes \mathcal{K}),\quad
\mathbf{r}\in H^0(\mathcal{E}(\mathfrak g_{1})\otimes \mathcal{K})=H^0(\mathcal{E}[R]\otimes \mathcal{K}).
\]

\begin{lem}\label{lem:deg}
In this setting we have the following:
\begin{enumerate}
\item[(i)] $\operatorname{ad}(e):\ \mathfrak g_{1}\rightarrow\mathfrak g_{-1}$ is an isomorphism.
\item[(ii)] $\dim\mathfrak g_{\pm2}=1$,
hence $\dim\mathfrak g=\dim\mathfrak g_0+2\dim\mathfrak g_1+2$.
\item[(iii)] For any $TC$--bundle $\mathcal{E}$:
\[
\deg \mathcal{E}[\mathfrak g_1]=-\tfrac{1}{2}\dim\mathfrak g_1\cdot\deg  \mathcal{E}[\mathfrak g_{-2}],
\qquad
\deg  \mathcal{E}[\mathfrak g_0]=0 .
\]
\end{enumerate}
\end{lem}

\begin{proof}
(i)  Consider the $\mathfrak{sl}_2$--triple $\mathfrak s=\langle e,h,f\rangle$. Decompose $\mathfrak g$ into irreducible
$\mathfrak s$-modules. Since all $\operatorname{ad}h$-weights lie in $[-2,2]$, each summand is some $V(n)$ with $n\le 2$. As $V(n)$ has weights $n,n-2,\dots,-n$, the
weight $\pm1$ occurs only in copies of $V(1)$ and the weight $\pm2$ only in copies
of $V(2)$; in particular $\mathfrak g_{\pm1}$ lies in the $V(1)$-isotypic part. On a
copy of $V(1)$ the operator $\operatorname{ad}(e)$ sends the weight-$1$ line
isomorphically onto the weight-$(-1)$ line, so summing over the isotypic component
\begin{equation}\label{eq:adf-iso}
\operatorname{ad}(e):\ \mathfrak g_{1}\xrightarrow{\ \sim\ }\mathfrak g_{-1}
\quad\text{is an isomorphism.}
\end{equation}

(ii) $\mathfrak g_{\pm2}$ is one--dimensional, spanned by $e, f$ respectively. Summing gives the stated identity.

(iii) Consider the $TC$--equivariant alternating form
$\omega\colon\Lambda^2\mathfrak g_1\to\mathfrak g_2$,
$\omega(x,y)=[x,y]$. From Part (i), $\operatorname{ad}(e):\mathfrak g_{1}\xrightarrow{\ \sim\ }\mathfrak g_{-1}$ is an isomorphism. Together with $B([x,y], e)=-B(x,[e,y])$, we see that $\omega$ is nondegenerate.

The nondegenerate $\mathfrak g_2$--valued symplectic form $\omega$ on $\mathfrak g_1$ 
trivializes $\det\mathfrak g_1\otimes\mathfrak g_2^{-m}$ ($m=\tfrac12\dim\mathfrak g_1$) as a $TC$--module. For a $TC$-bundle $\mathcal E$, the map
$\chi\mapsto\deg\mathcal E[\chi]$ is additive in the character $\chi$, and
$\deg\mathcal E[V]=\deg\mathcal E[\det(\operatorname{Ad}|_V)]$ for any representation
$V$. Thus
\[
\deg\mathcal E[\mathfrak g_1]=m\cdot\deg\mathcal E[\mathfrak g_2].
\]
The Killing form is an $TC$-invariant perfect pairing
$\mathfrak g_2\otimes\mathfrak g_{-2}\to\mathbb C$, so
$\det(\operatorname{Ad}|_{\mathfrak g_{-2}})=\det(\operatorname{Ad}|_{\mathfrak g_2})^{-1}$
and $\deg\mathcal E[\mathfrak g_{-2}]=-\deg\mathcal E[\mathfrak g_2]$. Hence
\[
\deg\mathcal E[\mathfrak g_1]=-m\cdot\deg\mathcal E[\mathfrak g_{-2}]
=-\tfrac12\dim\mathfrak g_1\cdot\deg\mathcal E[\mathfrak g_{-2}].
\]
The restriction
to $\mathfrak g_0$ of the Killing form of $\mathfrak g$ is nondegenerate. Hence $B_{\mathfrak g}|_{\mathfrak g_0}$ gives an $TC$-equivariant
identification $\mathfrak g_0\simeq \mathfrak g_0^\vee$. Thus,
$\mathcal E[\mathfrak g_0]$ is self-dual, and 
$\deg \mathcal E[\mathfrak g_0]=0$.
\end{proof}

\begin{prop}\label{prop:nv}
Let $(\mathcal{E},\mathbf{e}+\mathbf{r})$ be a polystable $(TC,\langle e\rangle\oplus R)$--Higgs pair. Then 
$\deg \mathcal{E}[\mathfrak g_{-2}]=2-2g$ if and only if $\mathbf{e}$ is nowhere vanishing. 

In such case, $\mathcal{E}[\mathfrak g_{-2}]\cong \mathcal{K}^{-1}$ and $\mathcal{E}[\mathfrak g_{2}]\cong \mathcal{K}$.
\end{prop}

\begin{proof}
Set $L:=\mathcal{E}[\mathfrak g_{-2}]\otimes \mathcal{K}$, so $\mathbf{e}\in H^0(X,L)$ and
$\deg L=0$. So either $\mathbf{e}\equiv 0$,
or $\mathbf{e}$ is nowhere vanishing and $L\cong\mathcal O_X$.

Suppose $\mathbf{e}\equiv 0$.
We use semistability against $s=-h\in i\mathfrak k_0$. Let $V=\langle e\rangle \oplus R$. Since
$h\in\mathfrak z(\mathfrak g_0)$, $P_{-h}=TC$ and the corresponding reduction is
trivial; we have:
\[V_{-h}=\bigoplus_{j\ge0}(\mathfrak g_j\cap V)=\mathfrak g_1,\] so
$\varphi=\mathbf r\in H^0(\mathcal{E}|V_{-h}])\otimes \mathcal{K})$ and the semistability requires
$\deg \mathcal{E}(\sigma,-h)\ge0$. Since we know the adjoint action of $h$ we see that this degree is
\[\deg \mathcal{E}(\sigma,-h)=-\sum_j j\deg \mathcal{E}[\mathfrak g_j]=-(2\deg \mathcal{E}[\mathfrak g_1]+4\deg \mathcal{E}[\mathfrak g_2]).\]
With $\deg \mathcal{E}[\mathfrak g_{-2}]=2-2g$, Lemma~\ref{lem:deg} gives 
$\deg \mathcal{E}(\sigma,-h)=(\dim\mathfrak g_1+4)(2-2g)<0$, a contradiction. So $\mathbf{e}\not\equiv0$,
hence nowhere vanishing.
\end{proof}

The moduli space associated to the restricted Slodowy slice is a subset:
$$\mathcal{M}_\mathbf{e}(TC, \langle e\rangle \oplus R)\subset \mathcal{M}(TC, \langle e\rangle \oplus R).$$
This subset is characterized by the property that the topological type of $T$ is such that $\mathcal{E}[\langle e \rangle]\simeq \mathcal{K}^{-1}$ complex line bundle, or equivalently that $\mathbf{e}$ is nowhere vanishing. 

\medskip

This moduli space decomposes into spaces depending of the topological type $\delta$ of $\mathcal{E}$:
$$\mathcal{M}_\mathbf{e}(TC, \langle e\rangle \oplus R)=\bigcup_{\delta} \mathcal{M}_\delta(TC, \langle e\rangle \oplus R).$$

\subsubsection{Toledo invariant}\label{sec:Toledo}
Assume $(\mathcal E, \varphi)\in \mathcal M(TC, \langle e\rangle\oplus R)$. The Toledo character $\chi_T:\mathfrak g_0\rightarrow \mathbb C$ is defined by
\[\chi_T(x)=B(\zeta, x)B(\gamma, \gamma),\] where $\zeta=T_{\alpha}$ is the coroot of the long root $\alpha$ and $\gamma$ is the longest root in $\mathfrak g_1=R$. When $G$ is simple complex, $B(\gamma, \gamma)=\kappa$ where $\kappa=2$ if $\mathfrak g$ is not of type C and $\kappa=1$ if $\mathfrak g$ is of type C. From \cite[Section 3.3]{GG} (also \cite{BC+23}), the Toledo invariant is defined as 
\[\tau(\mathcal E, \varphi)=\deg_{\chi_T}(\mathcal E),\] and from \cite[Section 5.3]{GG}, it satisfies  
\begin{equation}\label{eqn:Bound}-2\kappa(2g-2)\leq \tau(\mathcal E, \varphi)\leq \kappa(2g-2).\end{equation}

Note that $e\in \mathfrak g_{-\alpha}=\mathfrak g_{-2}$. Let $\eta$ be the character for the action of $TC$ acting on $\mathfrak g_{-2}=\langle e\rangle$. Thus $d\eta=-\alpha.$ Since $B(T_{\alpha}, x)=\alpha(x)$, $\chi_T=\kappa\alpha=-\kappa d\eta$. And \[\deg_{\chi_T}(\mathcal E)=\deg_{-\kappa d\eta}(\mathcal E)=-\kappa \deg_{d\eta}(\mathcal E)=-\kappa\deg(\mathcal E(\mathfrak g_{-2})).\]

Now we assume $(\mathcal E, \mathbf e+\mathbf r)\in \mathcal M_{\mathbf e}(TC, \langle e\rangle\oplus R)$. By Proposition \ref{prop:nv}, $\mathcal E[\mathfrak g_{-2}]=\mathcal K^{-1}$. Hence, $\tau(\mathcal E,\mathbf e+\mathbf r)=\kappa(2g-2)$, i.e. it achieves the upper bound in (\ref{eqn:Bound}). 

Garcia-Prada and González in \cite[Section 5.2]{GG} called the Higgs bundles which achieve the lower bound in (\ref{eqn:Bound})  maximal, and derived a Cayley correspondence between the moduli space of maximal $(TC, \langle e\rangle\oplus R)$-Higgs pairs with the moduli space of $\mathcal K^3$-twisted $(C_1, \langle e\rangle)$-Higgs pairs where $C_1$ is the stabilizer of a generic element of $\mathfrak g_1=R$.

Finally in the case of $(TC,\langle e\rangle \oplus R)$-Higgs pair, Corollary \ref{cor:SolutionHitchin equation compacibleGeneral} implies the following.
\begin{cor}
\label{cor:SolutionHitchin equation compacible}
Let $X$ be a closed hyperbolic surface and $(\mathcal{E}_{TC},\mathbf{e}+\mathbf{r})$ being  a polystable $(TC,\langle e\rangle \oplus R)$-Higgs pair arising from  the minimal Slodowy slice, where the triple $\lbrace e,h,f\rbrace$ is parallelisable. Then automatically any harmonic metric $H$ on this bundle is compatible with $TC$.

If its extension as a $G$-Higgs bundle is polystable, then there exists a harmonic metric $H$ on this bundle compatible with $TC$. Moreover, if it is stable, then the harmonic metric is unique.
\end{cor}

\subsection{Slodowy slice map and generalized Cayley correspondence}
 Fix a square root of the canonical bundle $\mathcal{K}^{1/2}$, and let $\mathcal{E}_T^o$ be a $T$-bundle such that $\mathcal{E}_T^o[\langle e\rangle]\simeq \mathcal{K}^{-1}$. Every two such $T$-bundle differ only by a $T\cap C$-bundle. Let $\mathcal{M}^{\mathcal{K}^{3/2}}(C,R)$ be the moduli space of $\mathcal{K}^{3/2}$-twisted $(C,R)$-pairs, meaning of pairs $(\mathcal{E}_C,\mathbf{r})$ with $\mathbf{r}\in H^0(\mathcal{E}_C[R]\otimes \mathcal{K}^{3/2})$. 

Because $T$ and $C$ commute, the multiplication $T\times C\to TC$ is a homomorphism and
one may form the associated $TC$--bundle $\mathcal{E}_T^o\star \mathcal{E}_C:=(\mathcal E_T^0\times
\mathcal E_C)\times_{T\times C}TC$ for any $C$--bundle $\mathcal E_C$. Define
\begin{equation}\label{eq:slodowymap}
\hat\Psi:\ (\mathcal E_C,\mathbf r)\ \longmapsto\
\bigl(\mathcal{E}_T^o\star \mathcal{E}_C,\mathbf e+\mathbf r\bigr),\qquad
\mathbf r\in H^0(\mathcal E_C[R]\otimes \mathcal K^{3/2}),
\end{equation}
from $\mathcal K^{3/2}$--twisted $(C,R)$--pairs to 
$(TC,\langle e\rangle\oplus R)$--Higgs pairs. We call $\hat\Psi$ the Slodowy slice map. 

\begin{lem}\label{lem:objbij}
The map $\hat\Psi$ is a bijection onto the $(TC,\langle e\rangle\oplus R)$--pairs whose
$\mathbf e$ is nowhere vanishing, with inverse
$(\mathcal E,\mathbf e+\mathbf r)\mapsto(\mathcal E_C,\mathbf r)$ obtained by reducing $(\mathcal E_T^0)^{-1}\otimes \mathcal E$ to $C$.
Moreover $\operatorname{Aut}(\mathcal E,\mathbf e+\mathbf r)\cong\operatorname{Aut}(\mathcal E_C,\mathbf r)$.
\end{lem}

\begin{proof}
Given $(\mathcal E_C,\mathbf r)$, the $\mathfrak g_{-2}$--component of $\hat\Psi$ is the
constant section $e$ (constant because $C$ centralizes $e$), which is nowhere
vanishing, and the $\mathfrak g_1$--component lies in $H^0(\mathcal{E}_C[\mathfrak g_1]\otimes \mathcal K^{3/2})=
H^0((\mathcal E_T^0\star\mathcal E_C)[\mathfrak g_1]\otimes \mathcal K)$. Conversely, let $(\mathcal E,\mathbf e+\mathbf r)$ have~$\mathbf e$ nowhere
vanishing. As $\mathfrak g_{-2}$ is one-dimensional, $C=\mathrm{Stab}_{TC}(e)$ and the open
$TC$--orbit of $\mathfrak g_{-2}$ is $\Omega=\g_{-2}\setminus\{0\}\cong TC/C$. From
$\mathcal E_T^0[\g_{-2}]\otimes \mathcal K\cong\mathcal O$ we get $\mathcal E[\mathfrak g_{-2}]\otimes \mathcal K\cong((\mathcal E_T^0)^{-1}\otimes
\mathcal E)[\mathfrak g_{-2}]$, whence $\mathbf e\in H^0\bigl(((\mathcal E_T^0)^{-1}\otimes \mathcal E)[\Omega]\bigr)$ is a reduction
of $(\mathcal E_T^0)^{-1}\otimes \mathcal E$ from $TC$ to $C$. Writing $\mathcal E_C$ for the reduced bundle,
$\mathcal E\cong \mathcal{E}_T^o\star \mathcal{E}_C$ and $\mathbf r\in H^0(\mathcal{E}[\g_1]\otimes \mathcal K)=H^0(\mathcal E_C[\mathfrak g_1]\otimes
\mathcal K^{3/2})$. The two assignments are mutually inverse.

For automorphisms, any $\psi\in\operatorname{Aut}(\mathcal E,\mathbf e+\mathbf r)$ preserves the $TC$--invariant splitting $\g_{-2}\oplus\g_1$ and hence $\langle e\rangle$ and $R$ separately; preserving $\mathbf e$ forces $\psi$ to reduce to $C=\mathrm{Stab}_{TC}(e)$, after which the condition on $\mathbf r$ is exactly that defining $\operatorname{Aut}(\mathcal E_C, \mathbf r)$. Thus $\operatorname{Aut}(\mathcal E,\mathbf e+\mathbf r)=\operatorname{Aut}(\mathcal E_C, \mathbf r)$.
\end{proof}

We will upgrade the above object bijection to the isomorphism between moduli spaces, Cayley correspondence. Before showing the isomorphism, we state two lemmas which will be used to compare the stability of the two twisted pairs. 

\begin{lem}\label{lem:testdecomp}
Let $\mathcal E=\mathcal{E}_T^o\star \mathcal{E}_C$ and $s\in i\mathfrak k_{TC}$, written as
$s=ch+s'$ with $c\in\mathbb R$ and $s'\in i\mathfrak k_{\cc}$. Then:
\begin{enumerate}
\item[(i)] $P^{TC}_s=T\cdot P^{C}_{s'}$ and $L^{TC}_s=T\cdot L^{C}_{s'}$; consequently
reductions $\sigma$ of $\mathcal E$ to $P^{TC}_s$ correspond bijectively to reductions $\sigma'$ of
$\mathcal E_C$ to $P^{C}_{s'}$ (carrying the fixed factor $\mathcal E_T^0$), and likewise for Levi
reductions.
\item[(ii)] $e\in(\g_{-2})_s$ if and only if $c\ge0$; hence every $\varphi$--admissible
$TC$--test of $(\mathcal E,\mathbf e+\mathbf r)$ has $c\ge0$. Here we call $(\sigma, s)$ is  ``$\varphi$--admissible”
if $\varphi=\mathbf e+\mathbf r\in H^0(\mathcal{E}_\sigma[(\g_{-2}\oplus\g_1)_s]\otimes\mathcal K)$. 
\item[(iii)] If $w\in\g_1$ has $\ad(s')$--weight $\mu$, then its $\ad(s)$--weight is
$c+\mu$, so $(\g_1)_s=\{\,w:\mu\le-c\,\}$. Thus when $c=0$ one has
$(\g_1)_s=(\g_1)_{s'}$ and $(\g_1)^0_s=(\g_1)^0_{s'}$, while when $c>0$ one has
$(\g_1)_s\subseteq(\g_1)_{s'}$.
\end{enumerate}
\end{lem}

\begin{proof}
(i) The parabolic $P^{TC}_s$ depends only on $\ad(s)|_{\g_0}$; since $h$ is central in
$\g_0$, $\ad(s)|_{\g_0}=\ad(s')|_{\g_0}$, so $P^{TC}_s=P^{TC}_{s'}$, and as $s'\in\cc$
and $T$ is central this equals $T\cdot P^{C}_{s'}$ (similarly for the Levi). The
identification of reductions follows from $\mathcal E(TC/P^{TC}_s)\cong \mathcal E_C(C/P^{C}_{s'})$,
the $T$--factor being fixed.\\
(ii) $\rho(e^{ts})e=e^{t\ad(s)}e=e^{t(c[h,\,\cdot\,]+[s',\,\cdot\,])}e=e^{-2ct}e$, using
$[h,e]=-2e$ and $[s',e]=0$; this is bounded as $t\to+\infty$ iff $c\ge0$.\\
(iii) On $\g_1$, $\ad(h)$ acts with weight $+1$, so $\ad(s)w=(c+\mu)w$; boundedness as
$t\to+\infty$ is $c+\mu\le0$.
\end{proof}

\begin{lem}\label{lem:degsplit}
With $s=ch+s'$ and corresponding reductions $\sigma\leftrightarrow\sigma'$ as in
Lemma~\ref{lem:testdecomp}(i),
\begin{equation}\label{eq:degsplit}
\deg \mathcal{E}(\sigma,s)\;=\;c\,\lambda\;+\;\deg \mathcal{E}_C(\sigma',s'),\qquad
\lambda:=B(h,h)\,(g-1)>0 .
\end{equation}
\end{lem}
\begin{proof}
The character $\chi_s=B(s,\,\cdot\,)$ splits as $\chi_s=c\,B(h,\,\cdot\,)+B(s',\,\cdot\,)$,
and $B(h,\cc)=0$ because $B(h,s')=B([f,e],s')=B(f,[e,s'])=0$ for $s'\in\cc$. Hence the
two summands are supported respectively on the $\mathcal E_T^0$-- and $\mathcal E_C$--factors of $\mathcal E$, and
their degrees add:
\[
\deg \mathcal E(\sigma,s)=\deg\!\bigl(\mathcal E_T^0\times_{cB(h,\cdot)}\mathbb C^\ast\bigr)+\deg \mathcal E_C(\sigma',s').
\]
The character $cB(h,\,\cdot\,)$ has weight $\chi_s(h)=c\,B(h,h)$ along $h$; with the
normalization (weight one $\mapsto \mathcal K^{1/2}$) it produces
$\mathcal K^{\,cB(h,h)/2}$, of degree $cB(h,h)(g-1)=c\lambda$. Finally $B(h,h)>0$ since $h\in
i\mathfrak k$ and $B$ is positive there.
\end{proof}

\begin{thm}\label{thm:cayley}
The Slodowy slice map $\hat\Psi$ in \eqref{eq:slodowymap} induces an isomorphism of the following two moduli spaces:
\[
\Psi:\ \mathcal M^{\mathcal K^{3/2}}(C,R)\ \longrightarrow\ \mathcal M_{\mathbf e}(TC,\langle e\rangle\oplus R).
\]
\end{thm}

\begin{proof}
The object bijection and the identification of automorphisms are
Lemma~\ref{lem:objbij}. It remains to check the equivalence between stability. 

\textit{For Semistability:}
Suppose $(\mathcal E_C,\mathbf r)$ is semistable and let $(s,\sigma)$ be a $\varphi$--admissible
$TC$--test of $(\mathcal E,\mathbf e+\mathbf r)$, $s=ch+s'$. By Lemma~\ref{lem:testdecomp}(ii), $c\ge0$. If
$c=0$ then by Lemma~\ref{lem:testdecomp}(iii) the test is exactly the $C$--test $(s',\sigma')$, which is
admissible ($\mathbf r\in H^0(\mathcal E_{\sigma'}[(\g_1)_{s'}]\otimes \mathcal K^{3/2})$), and
\eqref{eq:degsplit} gives $\deg \mathcal E(\sigma,s)=\deg \mathcal E_C(\sigma',s')\ge0$. If $c>0$ then
$(\g_1)_s\subseteq(\g_1)_{s'}$ by Lemma~\ref{lem:testdecomp}(iii), so $r$ is still $s'$--admissible and target
semistability gives $\deg \mathcal{E}_C(\sigma',s')\ge0$; hence by \eqref{eq:degsplit},
\[
\deg \mathcal{E}(\sigma,s)=c\lambda+\deg \mathcal{E}_C(\sigma',s')\ \ge\ c\lambda\ >\ 0 .
\]
Thus $(\mathcal E,\mathbf e+\mathbf r)$ is semistable. Conversely, if $(\mathcal E,\mathbf e+\mathbf r)$ is semistable and
$(s',\sigma')$ is an admissible $C$--test, take $s=s'$ ($c=0$): then $\mathbf e$ has
$s$--weight $0$ and $\mathbf r$ is admissible, so $(s,\sigma)$ is an admissible $TC$--test and
$0\le\deg \mathcal E(\sigma,s)=\deg \mathcal E_C(\sigma',s')$. Hence $(\mathcal E_C,\mathbf r)$ is semistable.

\textit{For stability:}
The strict tests range over $i\mathfrak k\setminus \ker(d\rho)$, determined by the representations
\emph{defining the pairs}: $\rho_{TC}$ on $\g_{-2}\oplus\g_1$ and $\rho_C$ on $\g_1$.
For $s'\in\cc$ one has $d\rho_{TC}(s')=0\iff d\rho_C(s')=0$, since $C$ acts trivially
on $\g_{-2}$ and by the same representation on $\g_1$; and $d\rho_{TC}(h)\neq0$ because
$h$ acts on $e$ with weight $-2$. So the strict $C$--directions ($c=0$) coincide on
the two sides. For a strict
$TC$--test $s=ch+s'$ with $c=0$ one gets $\deg \mathcal{E}(\sigma,s)=\deg \mathcal{E}_C(\sigma',s')>0$ by
target stability, and for $c>0$ one gets $\ge c\lambda>0$; conversely strict
$C$--tests are strict $TC$--tests with $c=0$. Hence $(\mathcal E,\mathbf e+\mathbf r)$ is stable iff
$(\mathcal E_C,\mathbf r)$ is.

\textit{For polystability:}
Suppose $(\mathcal E_C,\mathbf r)$ is polystable; it is then semistable, so $(\mathcal E,\mathbf e+\mathbf r)$ is
semistable by the above. Let $(s,\sigma)$ be a $\varphi$--admissible $TC$--test
with $\deg \mathcal{E}(\sigma,s)=0$ and $s=ch+s'$. If $c>0$ then \eqref{eq:degsplit} gives
$\deg \mathcal{E}(\sigma,s)\ge c\lambda>0$, a contradiction; so $c=0$, the test is the $C$--test
$(s',\sigma')$, and $\deg \mathcal{E}_C(\sigma',s')=0$. Polystability of $(\mathcal E_C, \mathbf r)$ provides a
reduction $\sigma'_0$ of $\mathcal{E}_C$ to $L^{C}_{s'}$ with $\mathbf r\in
H^0(\mathcal{E}_{C,\sigma'_0}[(\g_1)^0_{s'}]\otimes \mathcal K^{3/2})$. By Lemma~\ref{lem:testdecomp}(i)
this extends, through $L^{TC}_s=T\cdot L^{C}_{s'}$ and the fixed $\mathcal E_T^0$, to a reduction
$\sigma_0$ of $\mathcal E$ to $L^{TC}_s$; and $e$ has $s$--weight $0$ (as $c=0$ and $[s',e]=0$),
so $\mathbf e+\mathbf r\in H^0(\mathcal{E}_{\sigma_0}[(\g_{-2}\oplus\g_1)^0_s]\otimes K)$ ---showing the polystability
of $(\mathcal E,\mathbf e+\mathbf r)$. Conversely, given $(\mathcal E,\mathbf e+\mathbf r)$ polystable and an admissible
$C$--test $(s',\sigma')$ with $\deg \mathcal{E}_C(\sigma',s')=0$, take $s=s'$ ($c=0$); then
$\deg \mathcal{E}(\sigma,s)=0$, source polystability yields a Levi reduction $\sigma_0$ of $\mathcal{E}$
to $L^{TC}_s=T\cdot L^{C}_{s'}$ adapted to $\mathbf e+\mathbf r$, and discarding the fixed $\mathcal E_T^0$
leaves a reduction $\sigma'_0$ of $\mathcal E_C$ to $L^{C}_{s'}$ with $\mathbf r\in
H^0(\mathcal{E}_{C,\sigma'_0}[(\g_1)^0_{s'}]\otimes \mathcal K^{3/2})$. Hence $(\mathcal E,\mathbf e+\mathbf r)$ is
polystable iff $(\mathcal E_C, \mathbf r)$ is.

By the above discussion, both $\hat\Psi$ and its inverse (Lemma~\ref{lem:objbij}), $\Psi$ is an isomorphism of
the moduli spaces.
\end{proof}

From this isomorphism, we can decompose the moduli space decomposes into spaces depending of the topological type $\delta$ of the corresponding $C$-bundle $\mathcal{E}_C$:
$$\mathcal{M}_\mathbf{e}(TC, \langle e\rangle \oplus R)=\bigcup_{\delta} \mathcal{M}_\delta(TC, \langle e\rangle \oplus R).$$

\subsection{Local deformation structure of \texorpdfstring{$\mathcal M_e$}{Me}}
 By \cite{BiswasRamanan}, the space of infinitesimal deformations are
isomorphic to the hypercohomology group $\mathbb{H}^1(C^{\bullet})$ of the deformation complex:
\[C^{\bullet}\colon C^0=\mathcal{E}[\mathfrak g_0]\xrightarrow{\mathrm{ad}(\varphi)}C^1=\mathcal{E}[\langle e\rangle \oplus R]\otimes \mathcal{K}\]
The deformation complex induces the long exact sequence:
\[
\begin{aligned}
0 \longrightarrow \mathbb H^0(C^\bullet)
&\longrightarrow H^0(\mathcal{E}[\mathfrak g_0])
\xrightarrow{\operatorname{ad}(\varphi)}
H^0(\mathcal{E}[\langle e\rangle \oplus R])\otimes \mathcal{K})
\longrightarrow \mathbb H^1(C^\bullet) \\
&\longrightarrow H^1(\mathcal{E}[\mathfrak g_0])
\xrightarrow{\operatorname{ad}(\varphi)} H^1(\mathcal{E}[\langle e\rangle \oplus R]\otimes \mathcal{K})
\longrightarrow \mathbb H^2(C^\bullet)
\longrightarrow 0 .
\end{aligned}
\]
In particular, $\mathbb H^0(C^\bullet)$ is isomorphic to the space $\operatorname{aut}(\mathcal E,\varphi)$. Then 
\[-\chi(C^{\bullet})=\dim\mathbb H^1-\dim\mathbb H^0-\dim\mathbb H^2,\] is the so-called ``expected" dimension of the moduli space. 
\begin{prop}\label{prop:dim}
Fix a topological type $\delta$ of $TC$--bundles with
$\deg \mathcal{E}[\mathfrak g_{-2}]=2-2g$. Then $\mathcal{M}_\delta(TC, \langle e\rangle \oplus R)$ is of ``expected" complex dimension
\[
(g-1)\bigl(\dim_{\mathbb C}\mathfrak g-3\bigr).
\]
\end{prop}
\begin{proof} We obtain
\begin{eqnarray*}
-\chi(C^\bullet)&=&-\chi(\mathcal E[\g_0])+\chi(\mathcal E[\langle e\rangle\oplus R]\otimes \mathcal K)\\
&=&(g-1)\bigl(\dim\mathfrak g_0+\dim (\langle e\rangle\oplus R)\bigr)+\deg \mathcal{E}[\langle e\rangle\oplus R]-\deg \mathcal{E}[\mathfrak g_0]
\end{eqnarray*}
The nowhere-vanishing
of $\mathbf{e}$ is Proposition~\ref{prop:nv}. By assumption $\deg \mathcal{E}[\mathfrak g_{-2}]=2-2g$, so Lemma~\ref{lem:deg} gives
$\deg \mathcal{E}[\mathfrak g_{1}]=(g-1)\dim\mathfrak g_1$, $\deg \mathcal{E}[\mathfrak g_0]=0$, and
$\deg \mathcal{E}[\langle e\rangle\oplus R]=(g-1)(\dim\mathfrak g_1-2)$. Substituting
$\dim(\langle e\rangle\oplus R)=\dim\mathfrak g_1+1$, we have
\[-\chi(C^\bullet)=(g-1)(\dim\mathfrak g_0+2\dim\mathfrak g_1-1)=(g-1)(\dim\mathfrak g-3).\]
\end{proof}

\subsubsection{The vanishing of \texorpdfstring{$\mathbb H^0(C^\bullet)$}{H0(C.)}}

\begin{lem}\label{lem:FuchsianExample} Assume $(\mathcal E_C, 0)$ semistable and $\mathfrak z_\cc\neq0$ (the type $A$ case), then $(\mathcal E_C, 0)$ and $(\mathcal E,\mathbf e)$ are both 
\emph{strictly} semistable. 
\end{lem}

\begin{proof}

Pick $0\neq z\in\mathfrak z_\cc$. Since the action of $C$ on $R$ is faithful,
$z\notin\ker(d\rho)$, so $z$ is a nontrivial test direction for the $(C,R)$--pair. The
field being zero, both $z$ and $-z$ are admissible. Since
$\deg\mathcal E_C(\sigma_z,-z)=-\deg\mathcal E_C(\sigma_z,z)$, while semistability gives
$\deg\mathcal E_C(\sigma_z,\pm z)\ge0$, we get $\deg\mathcal E_C(\sigma_z,z)=0$. Thus $z$ is a nontrivial
admissible test attaining equality, so $(\mathcal E_C,0)$ is
strictly semistable. By the Cayley correspondence in Theorem \ref{thm:cayley}, $(\mathcal E, \mathbf e)$ is strictly semistable.
\end{proof}

\begin{lem}\label{lem:vanishingH0}
 If $(\mathcal E,\varphi)\in \mathcal M_{\mathbf e}(TC,\langle e\rangle+R)$ is stable, then $\mathbb H^0(C^\bullet)=\operatorname{aut}(\mathcal E, \varphi)=0$.
\end{lem}
\begin{proof}
By \cite[Proposition 2.14]{GGM}, $(\mathcal E,\mathbf e+\mathbf r)$ being stable implies \[\operatorname{aut}(\mathcal E,  \mathbf e+\mathbf r)\subset \{s\in H^0(\mathcal E[\mathfrak{z}_{\mathfrak g_0}]):[s, \mathbf e+\mathbf r]=0\}.\]

Note $\mathfrak z_{\mathfrak{g}_0}=\mathbb C h\oplus \mathfrak{z}_{\mathfrak c}$ and $H^0(\mathcal E[\mathfrak{z}_{\mathfrak g_0}])\cong \mathfrak{z}_{\mathfrak g_0}$.  Suppose $z\in \mathfrak z_{\mathfrak g_0}$ satisfies $[z,\mathbf e+\mathbf r]=0$, then $[z,\mathbf e]=0\in\mathfrak g_{-2}$ and $[z,\mathbf r]=0\in \mathfrak g_1$. Assume $\pi_{\mathbb C h}z=a\cdot h$ for $a\in \mathbb C$. From $[z,\mathbf e]=a[h,\mathbf e]=0$ and thus $a=0$. Thus $z\in \mathfrak z_{\mathfrak c}$.

Therefore, $(\mathcal E,\mathbf e+\mathbf r)$ being stable implies \[\operatorname{aut}(\mathcal E,  \mathbf e+\mathbf r)\subset \{s\in H^0(\mathcal E[\mathfrak{z}_{\mathfrak c}]):[s, \mathbf r]=0\}.\]

Case 1: $\mathfrak g$ is not of type $A$ or $D_3$, then $\mathfrak z_{\mathfrak c}=0$. Case 2: $\mathfrak g$ is of type $A$. Since $D_3\cong A_3$, it is enough to consider $A_n$. Write $\mathcal E[\mathbb C^{n+1}]=\mathcal L\oplus \mathcal L\mathcal K^{-1}\oplus V$, $\rk V=n-1$, with 
\[
\mathbf r=\begin{pmatrix}0&0&\beta\\ 0&0&0\\ 0&\alpha&0\end{pmatrix}\in H^0(\mathcal E[\mathfrak{g}_{1}\oplus\mathfrak g_{-2}]\otimes \mathcal{K}),\] where 
$\alpha\in H^0(\Hom(\mathcal L\mathcal K^{-1},V)\otimes \mathcal K)$ and $\beta\in H^0(\Hom(V,\mathcal L)\otimes \mathcal K)$. Write \[z=\operatorname{diag}(-(n-1)\lambda,-(n-1)\lambda,2\lambda\Id_{V})\in H^0(\mathcal E[\mathfrak z_{\mathfrak c}]).\] Then
\[
[z,\mathbf r]=\begin{pmatrix}0&0&-(n+1)\lambda\beta\\ 0&0&0\\ 0&(n+1)\lambda\alpha&0\end{pmatrix}.
\] Thus $[z,\mathbf r]=0$ if and only if $z=0$ or $\alpha, \beta$ both vanish, in which case $(\mathcal E, \varphi=\mathbf e)$ is strictly semistable, by Lemma \ref{lem:FuchsianExample}. 
\end{proof}

\subsubsection{Stability and the structure group extension}

\begin{prop}\label{prop:StableEquiv}
Let $(\mathcal E,\varphi)\in \mathcal M_{\mathbf e}(TC,\langle e\rangle+R)$. Then it is stable if and only if its
extension $F(\mathcal E,\varphi)$
is stable as a $G$-Higgs bundle.
\end{prop}
\begin{proof}
Suppose that $(\mathcal E,\varphi)$ is stable. By Corollary \ref{cor:SolutionHitchin equation compacibleGeneral}, $F(\mathcal E, \varphi)=(\mathcal E(G), \Phi)$ is polystable and it admits a compatible
harmonic metric $h$.
For $j\in\mathbb Z/3\mathbb Z$, set
\[
  \mathfrak a_{\bar j}
  :=
  \left\{
    u\in H^0\bigl(X,\mathcal E[\mathfrak g_{\bar j}]\bigr)
    :[u,\Phi]=0
  \right\}.
\]
\[
  \operatorname{aut}(F(\mathcal E,\varphi))=\mathfrak a
  =\mathfrak a_{\bar0}\oplus
   \mathfrak a_{\bar1}\oplus
   \mathfrak a_{\bar2}.
\]
The Bochner formula on the compact surface $X$ implies that
every $u\in\mathfrak a$ satisfies
\[
  D_hu=0,\qquad [u,\Phi^{*_h}]=0;
\]
see \cite[Corollary~3.12]{LZ26}.
Consequently, taking the Hermitian adjoint induces an
antilinear isomorphism
\[
  \mathfrak a_{\bar j}
  \longrightarrow
  \mathfrak a_{-\bar j}.
\]

By Lemma \ref{lem:vanishingH0}, we have $\mathfrak a_{\bar0}=\operatorname{aut(\mathcal E,\varphi)}=0$.

Now let $u\in\mathfrak a_{\bar1}$. By the Jacobi identity, \[[[u,u^{*_h}],\Phi]=-[[u^{*_h},\Phi],u]-[[\Phi,u],u^{*_h}]=0,\] we have $[u,u^{*_h}]
  \in\mathfrak a_{\bar0}=0.$
Thus $\operatorname{ad}(u)$ is normal with respect to the
Hermitian form induced by $h$, and since normal matrices are unitarily diagonalizable, $u$ is pointwise semisimple.

By Proposition \ref{prop:nv}, $\mathcal E[\mathfrak g_{-2}]=\mathcal K^{-1}$ and $H^0(X,\mathcal E[\mathfrak g_{-2}])=0.$
The $\mathfrak g_{-2}$-component of $u$ therefore vanishes,
and $u$ takes values in 
$\mathfrak g_1$.
Every element of $\mathfrak g_1$ is nilpotent, since its adjoint
action strictly raises the grading.
Hence $u$ is both semisimple and nilpotent, and therefore $u=0$.
We conclude that $\mathfrak a_{\bar1}=0$; taking adjoints also
gives $\mathfrak a_{\bar2}=0$.

Thus $\operatorname{aut}(F(\mathcal E,\varphi))=\mathfrak a=0$.
Since $F(\mathcal E,\varphi)$ is polystable and $G$ is simple, $\operatorname{aut}(F(\mathcal E,\varphi))=0$ implies that $F(\mathcal E,\varphi)$
is stable; see \cite[Proposition~2.14]{GGM}.

The reverse direction is obvious.
\end{proof}

\subsubsection{The vanishing of \texorpdfstring{$\mathbb H^2(C^\bullet)$}{H2(C*)}}
\begin{lem}\label{prop:vanishingH2}
For $(\mathcal{E},\varphi)\in\mathcal M_{\mathbf e}(TC,\langle e\rangle\oplus R)$, if $F(\mathcal E, \varphi)$ is stable, then we have $\mathbb H^2(C^\bullet)=0$.
\end{lem}
\begin{proof}
By Serre duality $\mathbb H^2(C^{\bullet})\cong\mathbb H^0(C^{\bullet\vee})^{\ast}$,
where $C^{\bullet\vee}\colon \mathcal{E}[(\langle e\rangle \oplus R)^{\ast}]\xrightarrow{-\mathrm{ad}\,\varphi}
\mathcal{E}[\mathfrak g_0]\otimes \mathcal{K}$ and $(\langle e\rangle \oplus R)^{\ast}\cong\mathfrak g_{2}\oplus\mathfrak g_{-1}$.
Hence
\[
\mathbb H^2(C^{\bullet})^{\ast}\cong 
\bigl\{(w_{2},w_{-1})\in H^0(\mathcal{E}[\mathfrak g_{2}])\oplus H^0(\mathcal{E}[\mathfrak g_{-1}]):
[\mathbf{e},w_{2}]+ [\mathbf{r},w_{-1}]=0\bigr\}.
\]
The brackets $[\mathbf{r},w_{2}]$ and $[\mathbf{e},w_{-1}]$ vanishing as they land
in $\mathfrak g_{\pm3}=0$.

We can view $(w_2, w_{-1})$ as the holomorphic section in $H^0(\mathcal E[\g])$. Thus it is $\ad(\varphi)$-invariant. From Proposition \ref{prop:StableEquiv} and the stability of $(\mathcal E,\varphi)$, $F(\mathcal E,\varphi)$ is stable. Thus there is no such section and  the statement follows.
\end{proof}

\subsubsection{Existence of stable Higgs pair and topological invariant}

In most cases, taking a stable bundle $\mathcal{E}_C$ guarantees the stability, ensuring that we have many elements in the slice.

\begin{lem}\label{lem:EC-implies-stable-ss}
Suppose $\mathfrak c=\operatorname{Lie}(C)$ is semisimple
(equivalently $\mathfrak z_{\mathfrak c}=0$); this holds for
$\mathfrak g=\mathfrak g_2,\mathfrak f_4,\mathfrak e_6,\mathfrak e_7,\mathfrak e_8,
\mathfrak{sp}_{2m},\mathfrak{so}_N$ ($N\neq 6$).
If the $C$-bundle $\mathcal E_C$ is stable, then
$(\mathcal E_C,\mathbf r)$ and $(\mathcal E,\mathbf e+\mathbf r)$ are stable.
\end{lem}

\begin{proof}
For the principal bundle, $\mathcal E_C$ is stable iff
$\deg\mathcal E_C(\sigma,s)>0$ for \emph{every} $s$ satisfying $0\neq s\in i\mathfrak k_{\cc}\setminus \ker(\ad_{\cc})=i\mathfrak k_{\cc}$ and associated reduction $\sigma$ to $P_s$. Here, we use the assumption $\cc$ is semisimple and obtain $\ker(\ad_{\cc})=0$.

For the $(C,R)$--Higgs pair, with $\rho:C\to \GL(R)$ and Higgs field
$\mathbf r\in H^0(\mathcal E[R]\otimes \mathcal K^{3/2})$, stability requires
$\deg\mathcal E_C(\sigma,s)>0$ only for those $(\sigma,s)$ that are \emph{admissible}: namely
$s\in i\mathfrak k_{\cc}\setminus\ker(d\rho)$ and $\mathbf r\in H^0\bigl(\mathcal E_{C,\sigma}[V_s]\otimes
\mathcal K^{3/2}\bigr)$, i.e.\ $\mathbf r$ lies in the non--negative part $V_s$. But this is already implied by the stability of $\mathcal E_C$. The stability of $(\mathcal E, \mathbf e+\mathbf r)$ follows from the Cayley correspondence in Theorem \ref{thm:cayley}.
\end{proof}

In the case of type $A$, we proceed more directly to show that there are stable pairs. Recall that the centralizer group $C$ is \[C=\{\operatorname{diag}(a, a, A), a\in \C^*, A\in \GL(n-1,\C), \det(A)=a^{-2}\}, \]

Define the character $\chi_0:C\rightarrow \C^*$ as taking $\operatorname{diag}(a, a, A)$ to $a$. Given $(\mathcal E, \varphi)\in \mathcal M_{\mathbf e}(TC, \langle e\rangle\oplus R)$, consider the associated $C$-bundle $\mathcal E_C$ from the Cayley correspondence, then we define $\mathcal M_{\mathbf e}^d(TC,\langle e\rangle\oplus R)$ the subspace of $\mathcal M_{\mathbf e}(TC,\langle e\rangle\oplus R)$ such that $\deg(\mathcal E_C(\chi_0))=d$.  Since $(\chi_0)_*:\pi_1(C)\rightarrow \mathbb Z$ is injective, $\deg(\mathcal E_T(\chi_1))$ is the complete topology invariant for principal $C$-bundles. Thus $\mathcal M_{\mathbf e}(TC,\langle e\rangle\oplus R)$ is decomposed into 
$$\mathcal{M}_\mathbf{e}(TC, \langle e\rangle \oplus R)=\bigcup_{d} \mathcal M_{\mathbf e}^d(TC, \langle e\rangle \oplus R).$$

Let $(\mathcal E, \varphi)\in \mathcal M_{\mathbf e}(TC,\langle e\rangle\oplus R)$. We consider the associated Higgs vector bundle of $(\mathcal E,\varphi)$ by the natural representation of $\SL(n+1,\mathbb C)$ on $\mathbb C^{n+1}$. Write $\mathcal E[\mathbb C^{n+1}]=\mathcal L\oplus \mathcal L\mathcal K^{-1}\oplus V$, $\rk V=n-1$, with 
\[\mathbf e=\begin{pmatrix}0&0&0\\ 1&0&0\\ 0&0&0\end{pmatrix}\in H^0(\mathcal E[g_{-2}]\otimes \mathcal{K}),\quad
\mathbf r=\begin{pmatrix}0&0&\beta\\ 0&0&0\\ 0&\alpha&0\end{pmatrix}\in H^0(\mathcal E[\mathfrak{g}_{1}]\otimes \mathcal{K}),\] where 
$\alpha\in H^0(\Hom(\mathcal L\mathcal K^{-1},V)\otimes \mathcal K)$ and $\beta\in H^0(\Hom(V,\mathcal L)\otimes \mathcal K)$. Recall the group \[T=\{\operatorname{diag}(b, b^{-1}, I_{n-1}), b\in \C^*\}. \] Let $\chi_1:T\rightarrow \C^*$ be defined as taking $\operatorname{diag}(b, b^{-1}, I_{n-1})$ to $b$. Note that $\deg\mathcal L=\deg(\mathcal E(\chi_0\chi_1))$ for $\chi_0\chi_1:TC\rightarrow \mathbb C^*$. Since $\deg(\mathcal E_T^0(\chi_1))=g-1$, we have \[\deg(\mathcal L)=d+(g-1).\] 

\begin{lem}\label{lem:SLstable}
 For $\mathfrak g=\mathfrak{sl}_{n+1}$, $\mathcal M_{\mathbf e}^0(TC,\langle e\rangle\oplus R)$ contains a stable point. 
\end{lem}
\begin{proof} 
Let $V$ be a stable vector bundle of rank $n-1$ and degree $0$. By Riemann-Roch theorem,   
\[\dim_{\mathbb C}H^0(\mathcal K^{\frac{3}{2}}V)\geq \deg(\mathcal K^{\frac{3}{2}}V)+(n-1)(1-g)=2(g-1)(n-1).\]
\[\dim_{\mathbb C}H^0(V^*\mathcal K^{\frac{3}{2}})\geq \deg(V^*\mathcal K^{\frac{3}{2}})+(n-1)(1-g)=2(g-1)(n-1).\] Let $0\neq \alpha\in H^0(\mathcal K^{\frac{3}{2}}V)$ and $0\neq \beta\in H^0(\mathcal K^{\frac{3}{2}}V^*)$.  On $\mathcal E[\mathbb C^{n+1}] = \mathcal K^{\frac{1}{2}} \oplus\mathcal K^{-\frac{1}{2}} \oplus V$, we consider the Higgs field
\[
  \Phi \;=\;
  \begin{pmatrix}
     0 & 0 & \beta\\
     1 & 0 & 0 \\
     0 & \alpha& 0
  \end{pmatrix}.\] It is the associated Higgs vector bundle of a $(TC,\langle e\rangle\oplus R)$-Higgs pair $(\mathcal E,\varphi)$ by the natural representation of $\SL(n+1,\mathbb C)$ on $\mathbb C^{n+1}$. 

From \cite[Proposition A.9]{LZ26}, to check the stability of $(\mathcal E, \varphi)$, it suffices to check the stability of $(\mathcal E[\mathbb C^{n+1}], \Phi)$ as a cyclic Higgs bundle, that is, for any proper $\Phi$-invariant subbundle $F=\oplus_{i=1}^3 F_i$, where $F_i\subset E_i$. Here $E_1=\mathcal K^{\frac{1}{2}}, E_2=\mathcal K^{-\frac{1}{2}}$ and $E_3=V$.

By the $\Phi$-invariance, $F_1$ and $F_2$ are either both zero or both nonzero. 

Therefore, we only need to check $F$ is either a subbundle of $V$ or $\mathcal K^{\frac{1}{2}}\oplus \mathcal K^{-\frac{1}{2}}\oplus V_1$ where $V_1$ is a subbundle of $V$. Since $\alpha\neq 0, \beta\neq 0$, $F$ cannot be $V$ or $\mathcal K^{\frac{1}{2}}\oplus \mathcal K^{-\frac{1}{2}}$.

From the stability of $V$ and $\deg(V)=0$, we have $\deg F<0$. Thus $(\mathcal E, \varphi)$ is stable as $(TC,\langle e\rangle\oplus R)$-Higgs pair. 
\end{proof}

We have the following  control of the topological invariant $d$. 
\begin{prop}\label{prop:controlDegree}
 For $\mathfrak g=\mathfrak{sl}_{n+1}$, we have $d\in [-(g-1),g-1]$ if and only if $\mathcal M_{\mathbf e}^d(TC,\langle e\rangle\oplus R)\neq \emptyset$. 

 If $n=2$, every element in $\mathcal M_{\mathbf e}^d(TC,\langle e\rangle\oplus R)$ is stable except the case $\mathbf r=0$ and $d=0$. 

 If $n\geq 3$, $\mathcal M_{\mathbf e}^{\pm(g-1)}(TC,\langle e\rangle\oplus R)$ does not contain any stable point. If $n=2$, $\mathcal M_{\mathbf e}^{-(g-1)}(TC,\langle e\rangle\oplus R)$ and $\mathcal M_{\mathbf e}^{(g-1)}(TC,\langle e\rangle\oplus R)$ are mapped onto the Hitchin section for $\SL(3,\C)$ under the structure extension map $F$.
\end{prop}

For the $\SL(3,\C)$ case, $\mathcal M_{\mathbf e}^{0}(TC,\langle e\rangle\oplus R)$ corresponds to the representations which deform from the Barbot representations, containing the $\SL(3,\R)$-representations studied in \cite{BD25}. 

\begin{proof} We use the expression of the Higgs vector bundle $(\mathcal E[\C^{n+1}], \Phi)$ of $(\mathcal E, \varphi)$ before Lemma \ref{lem:SLstable}. Note that $(\mathcal E[\C^{n+1}], \Phi)$ is (poly-)stable if and only if  $(\mathcal E, \varphi)$ is (poly-)stable, from Proposition \ref{prop:StableEquiv}. Denote $l=\deg\mathcal L$. Recall the group \[T=\{\operatorname{diag}(b, b^{-1}, I_{n-1}), b\in \C^*\}. \] Let $\chi_1:T\rightarrow \C^*$ be defined as taking $\operatorname{diag}(b, b^{-1}, I_{n-1})$ to $b$. Note that $\deg\mathcal L=\deg(\mathcal E(\chi_0\chi_1))$ for $\chi_0\chi_1:TC\rightarrow \mathbb C^*$. Since $\deg(\mathcal E_T(\chi_1))=g-1$, we have $l=d+(g-1)$. 

We have $\deg(V)=2g-2-2l.$

Suppose $\alpha\neq 0$. Let $\operatorname{Im}\alpha$ denote the line subbundle in $V$ generated by the image of $\alpha$. Since $\alpha: \mathcal L\mathcal K^{-1}\rightarrow \operatorname{Im}\alpha\otimes \mathcal K$ is nonzero, then $\deg(\operatorname{Im}\alpha)\geq \deg(L\mathcal K^{-2})=l-2(2g-2)$. 

Suppose $\beta\neq 0$. Let $\Ker\beta$ denote the line subbundle in $V$ generated by the kernel of $\beta$. If $n\geq 3$, $\beta$ induces a nonzero map $ V/\Ker\beta\rightarrow \mathcal L\otimes \mathcal K$ and then $\deg(V/\Ker\beta)\leq \deg(L\mathcal K)=l+(2g-2)$. Using $\deg(V)=2g-2-2l,$ we have $\deg(\Ker\beta)\geq -3l.$ 

Case 1: $\alpha\neq 0$, $\beta\neq 0$.
If $n\geq 3$, $\Ker\beta$ and $\mathcal L\oplus \mathcal{LK}^{-1}\oplus \operatorname{Im}\alpha$ are two proper $\Phi$-invariant subbundles of $\mathcal E$. By the polystability, $\deg(V)\leq 0$ and $\deg(\mathcal L\oplus \mathcal{LK}^{-1}\oplus \operatorname{Im}\alpha)\leq 0$. Hence we have $-3l\leq 0$, $3l-3(2g-2)\leq 0$. Thus $l\in [0,2g-2]$. Moreover, if $n\geq 3$, these two subbundles are proper. If we assume the stability, $l\in (0,2g-2)$. 

Case 2: $\alpha\neq 0$, $\beta=0$. Then $V$ and $\mathcal L\oplus \mathcal{LK}^{-1}\oplus \operatorname{Im}\alpha$ are two $\Phi$-invariant subbundles of $\mathcal E$. By the polystability, $\deg(V)\leq 0$ and $\deg(\mathcal L\oplus \mathcal{LK}^{-1}\oplus \operatorname{Im}\alpha)\leq 0$. Hence we have $l\in [g-1,2g-2]$. Moreover, if $n\geq 3$, these two subbundles are proper. If we assume the stability, $l\in (g-1,2g-2)$. 

Case 3: $\alpha=0$, $\beta\neq 0$. Then $\Ker\beta$ and $\mathcal L\oplus \mathcal{LK}^{-1}$ are two proper $\Phi$-invariant subbundles of $\mathcal E$. By the polystability, $\deg(\ker\beta)\leq 0$ and $\deg(\mathcal L\oplus\mathcal{LK}^{-1})\leq 0$. Hence we have $l\in [0,g-1]$. Moreover, if we assume the stability, $l\in (0,g-1)$. 

Case 4: $\alpha=0$, $\beta=0$. Then $V$ and $\mathcal L\oplus \mathcal{LK}^{-1}$ are two proper $\Phi$-invariant subbundles of $\mathcal E$. By the polystability, $\deg(V)\leq 0$ and $\deg(\mathcal L\oplus\mathcal{LK}^{-1})\leq 0$. Hence we have $l=g-1$. 

Note that $d=l-(g-1)$. In conclusion, we obtain the desired control of the degree $d$ for $\mathcal M_{\mathbf e}^d(TC,\langle e\rangle \oplus R)$ to contain a (poly-)stable point. 

On the other hand, for each $d\in [-(g-1),g-1]\cap \mathbb Z$, we can construct a polystable element in $\mathcal M_{\mathbf e}^d(TC, \langle e\rangle\oplus R)$ from a direct sum of an element in $\mathcal M_{\mathbf e}^{d, \SL(3,\mathbb C)}(TC, \langle e\rangle\oplus R)$ with a stable vector bundle $V_1$ of rank $n-2$ and degree $0$.

Suppose $n=2$. Set $l=\deg(L)$.

Case 1: $\alpha\neq 0$, $\beta\neq 0$. From \cite[Proposition A.9]{LZ26}, to check the stability of $(\mathcal E, \varphi)$, it suffices to check any proper $\Phi$-invariant subbundle $F=\oplus_{i=1}^3 F_i$, where $F_i\subset E_i$. Here $E_1=\mathcal K^{\frac{1}{2}}, E_2=\mathcal K^{-\frac{1}{2}}$ and $E_3=V$. By the $\Phi$-invariance, $F_1$ and $F_2$ are either both zero or both nonzero.  Since $V$ is a line bundle, there is no such subbundle to check. Therefore it is stable. 

Case 2: $\alpha\neq 0$, $\beta=0$. The only proper $\Phi$-subbundles are $V$ and $\mathcal{LK}^{-1}\oplus V$. Their degrees are both negative if and only if $0<l<g-1$.

Case 3: $\alpha=0$, $\beta\neq0$. The only proper $\Phi$-subbundles are $\mathcal L\oplus \mathcal{LK}^{-1}$ and $\mathcal{LK}^{-1}$. Their degrees are both negative if and only if $g-1<l<2g-2$.

Case 4: $\alpha=0$, $\beta=0$. Then $V$ and $\mathcal L\oplus \mathcal{LK}^{-1}$ are both $\Phi$-invariant proper subbundles. Thus $l=g-1$ and the Higgs bundle is strictly polystable. 

In conclusion, $F(\mathcal E,\varphi)$ is stable except it is in Case 4: $\alpha=0$, $\beta=0$ which forces $l=g-1$. 

Now we focus on the extremal values of $d$. 
When $d=(g-1)$ ($l=2g-2$), we have $\alpha\in H^0(\mathcal {LK}^{-1}V\mathcal K)$ with $\deg(\mathcal {L}V)=0$. Since $V=\mathcal{L}^{-2}\mathcal K$,  $\mathcal L=\mathcal K$ and $\alpha$ is an isomorphism. Thus $(\mathcal E, \varphi)$ is in the Hitchin section. 

When $d=-(g-1)$ ($l=0$), we have $\beta\in H^0(V^*\mathcal {LK})$ with $\deg(V^*\mathcal{LK})=0$. Since $V=\mathcal{L}^{-2}\mathcal K$,  $\mathcal{L}=\mathcal{O}$ and $\beta$ is an isomorphism. Thus $(\mathcal E, \varphi)$ is in the Hitchin section. 

The surjectivity of Hitchin section follows from every element in Hitchin section is stable. 
\end{proof}

\subsubsection{Summary of Moduli space}
In the following, we summarize the results of moduli spaces we obtain, using the local computations of hypercohomology.
\begin{prop}\label{prop:Moduli}
$\mathcal M_{\mathbf e}(TC,\langle e\rangle\oplus R)$ is a complex quasi-projective variety of pure dimension $(g-1)(\dim\mathfrak g-3)$. Moreover, each stable point $(\mathcal E,\varphi)$ has a neighborhood homeomorphic to $\mathbb H^1\sslash\operatorname{Aut}(\mathcal E,\varphi)$, where $\operatorname{Aut}(\mathcal E,\varphi)$ is finite.  

For $\mathfrak g$ not of type $A$ or $D_3$ and each $\delta\in \pi_1(C)$, the moduli space $\mathcal M_{\delta}(TC,\langle e\rangle\oplus R)$ is a complex quasi-projective variety  of dimension $(g-1)(\dim\mathfrak g-3)$.
\end{prop}
\begin{proof}
For $\mathfrak g$ neither being of type A nor $D_3$, by Lemma \ref{lem:EC-implies-stable-ss}, $\mathcal M_\delta(TC,\langle e\rangle\oplus R)$ has a stable point for each given $\delta$ by choosing a stable $C$--bundle $\mathcal E_C$ of type $\delta$. Thus the moduli space $\mathcal M_{\delta}(TC,\langle e\rangle\oplus R)$ is nonempty and contains a stable point. For $\mathfrak g$ being of type A, by Lemma \ref{lem:SLstable}, $\mathcal M_{\mathbf e}(TC,\langle e\rangle\oplus R)$ has a stable point. Since $A_3\cong D_3$, by Lemma \ref{lem:SLstable}, $\mathcal M_{\mathbf e}(TC,\langle e\rangle\oplus R)$ has a stable point. Thus for each $\mathfrak g$, $\mathcal M_{\mathbf e}(TC,\langle e\rangle\oplus R)$ contains a stable point.

Moreover, $\mathbb H^2=0$ by Proposition~\ref{prop:vanishingH2} at stable points. Using standard slice methods of Kuranishi (see \cite[Chapter 7.3]{kobayashi2014differential} for details for the moduli space of holomorphic bundles and \cite{fan2021analytic} for that of Higgs bundles, which also works for that of general Higgs pairs), when $\mathbb{H}^2=0$, a neighborhood of $(\mathcal{E},\varphi)$ in the moduli space is isomorphic to
\[\mathbb{H}^1\sslash\operatorname{Aut}(\mathcal{E},\varphi).\] For a stable point $(\mathcal E,\varphi)$, $\mathbb H^0=0$ by
Lemma~\ref{lem:vanishingH0} and thus $\operatorname{Aut}(\mathcal E,\varphi)$ is finite. In particular, the moduli space is of dimension
$-\chi(C^{\bullet})=(g-1)(\dim\mathfrak g-3)$ from Proposition~\ref{prop:dim}. 
\end{proof}

In the following, we show the structure extension map is a proper map with finite fiber.

\begin{prop}\label{prop:properFInite}
Let $X$ be a compact connected Riemann surface and $G$ a connected
complex semisimple group. Let $s$ be an automorphism
of $G$ of finite order $m$. Set $\zeta=\exp(2\pi i/m)$ and write
\[
  \mathfrak g=\bigoplus_{j\in\mathbb Z/m\mathbb Z}\mathfrak g_{\bar j},
  \qquad ds|_{\mathfrak g_{\bar j}}=\zeta^j\,\mathrm{id},
  \qquad G_0=(G^s)^0.
\]
Let $\mathcal M(G_0,\mathfrak g_{\bar 1})$ be the moduli space of polystable pairs
$(\mathcal{E}_0,\varphi)$ with
$\varphi\in H^0(X,\mathcal{E}_0(\mathfrak g_{\bar 1})\otimes K_X)$. Extension of structure group defines a proper
map with finite fibers
\[
  F:\mathcal M(G_0,\mathfrak g_{\bar 1})\longrightarrow\mathcal M(G),\qquad
  [\mathcal{E}_0,\varphi]\longmapsto [\mathcal{E}_0\times_{G_0}G,\Phi],
\]
where $\Phi$ is induced by the inclusion $\mathfrak g_{\bar 1}\subset\mathfrak g$.
\end{prop}
\begin{proof}
We show the extension map has finite fiber. 
Choose an $s$-invariant maximal compact subgroup $K\subset G$,
and put $K_0=K\cap G_0$.
Fix $[\mathcal{E},\Phi]$ in the image of $F$, together with a harmonic
metric $h$ induced by one of its polystable preimages.
Represent any other preimage by a holomorphic reduction
$\mathcal{E}_0\subset \mathcal{E}$, and choose a harmonic $K_0$-reduction
$Q_0\subset \mathcal{E}_0$. Its extension $Q_0K$ defines a harmonic
metric $h'$ on $(\mathcal{E},\Phi)$.
By uniqueness of harmonic metrics up to Higgs-bundle
automorphisms, there exists $b\in\operatorname{Aut}(\mathcal{E},\Phi)$
such that $b^*h=h'$.
Thus
\[
\widetilde{\mathcal{E}}_0:=b(\mathcal{E}_0),\qquad
\widetilde Q_0:=b(Q_0),\qquad
P_h=\widetilde Q_0K.
\]
The adjusted reduction represents the same source
isomorphism class.

Define a holomorphic $s$-twisted automorphism by
\[
u(e_0g)=e_0s(g),
\qquad e_0\in\widetilde{\mathcal{E}}_0,\quad g\in G.
\]
Since $s$ fixes $G_0$ pointwise, this is well-defined.
Moreover,
\[
u^m=1,\qquad u(P_h)=P_h,\qquad u_*\Phi=\zeta\Phi,
\]
where $u_*[e,v]=[u(e),ds(v)]$.
The second identity follows from
$P_h=\widetilde Q_0K$ and $s(K)=K$, and the third from
$V\subseteq\mathfrak g_{\bar1}$.

Let $H=\operatorname{Aut}(\mathcal{E},\Phi,h)$, a compact Lie group,
and choose one such twisted automorphism $u_{\mathrm{ref}}$.
It normalizes $H$, and every other $u$ constructed above
satisfies $u\,u_{\mathrm{ref}}^{-1}\in H$.
Hence all these automorphisms belong to the compact Lie group
\[
\widehat H
=\langle H,u_{\mathrm{ref}}\rangle
=\bigcup_{j=0}^{m-1}Hu_{\mathrm{ref}}^j.
\]
By \cite[Proposition~7.1]{Rezk}, the
$\widehat H$-conjugacy classes in
$\operatorname{Hom}(\mathbb Z/m\mathbb Z,\widehat H)$
are open.
Evaluation at a generator identifies this space with the compact set
\[
\{u\in\widehat H:u^m=1\}.
\]
Its quotient by conjugation is therefore compact and discrete,
hence finite.
Since $H$ has finite index in $\widehat H$, there are also
only finitely many relevant $H$-conjugacy classes.

Conjugating $u$ by $a\in H$ replaces $\widetilde{\mathcal{E}}_0$ by
$a(\widetilde{\mathcal{E}}_0)$, an isomorphic source pair.
For each fixed $u$ arising from a reduction, its fixed-point
bundle is
\[
\mathcal{E}^u:=\{e\in \mathcal{E}:u(e)=e\}=\widetilde{\mathcal{E}}_0G^s.
\]
Thus $u$ determines the extended $G^s$-reduction uniquely.
Its $G_0$-reductions correspond to sections of the finite covering
\[
\mathcal{E}^u/G_0\longrightarrow X.
\]
Since $X$ is connected, each section is determined by its value
at one point, so there are at most $|G^s/G_0|$ such reductions.
The prescribed topological type and the $V$-valued Higgs
condition only restrict this finite set.
Consequently, each relevant conjugacy class contributes
finitely many source points, and $F$ has finite fibers.

\medskip

Next, we prove the properness of the extension map. 
We adapt the compactness argument of
\cite[Lemmas~2.7--2.8]{Simpson1992}.
Let $C\subset\mathcal M(G)$ be compact and take a sequence
$x_n=[\mathcal{E}_{0,n},\varphi_n]\in F^{-1}(C)$.
Choose harmonic representatives $(A_n,\varphi_n)$ on a fixed
smooth $K_0$-bundle of topological type $\delta$.
Their extensions solve the $G$-Hitchin equations.

The induced adjoint vector Higgs bundles are therefore harmonic.
Their Higgs fields
$\Psi_n=\operatorname{ad}(\varphi_n)$
have uniformly bounded characteristic coefficients, since these
coefficients depend continuously on $F(x_n)\in C$.
By \cite[Lemma~2.7]{Simpson1992}, the norms
$\|\Psi_n\|_{L^\infty}$ are uniformly bounded.
Since $\mathfrak g$ is semisimple, $\operatorname{ad}$ is
injective, so $\|\varphi_n\|_{L^\infty}$ is uniformly bounded
as well. The equation
\[
F_{A_n}+[\varphi_n,\varphi_n^*]=0
\]
then gives a uniform curvature bound of $A_n$.

Following the proof of \cite[Lemma~2.8]{Simpson1992},
we apply Uhlenbeck compactness with structure group $K_0$.
Together with elliptic regularity, this yields a smoothly
convergent subsequence modulo $K_0$-gauge transformations.
The limit satisfies the same equations on the fixed smooth
bundle. Moreover, the $\mathfrak g_{\bar 1}$-valued condition is preserved,
since $V$ is a closed $K_0$-invariant subspace.
Thus the limit defines a point
$x_\infty\in\mathcal M_\delta(G_0,\mathfrak g_{\bar 1})$.
By continuity, $F(x_\infty)\in C$.

Consequently, $F^{-1}(C)$ is compact, proving the
properness. 
\end{proof}

\subsection{Joint moduli space}
Let $S$ be a closed oriented surface of genus $g\geq 2$. Let $\mathcal{JM}_\mathbf{e}(TC, \langle e\rangle\oplus R)$ be the parameter space of pairs $(X, (\mathcal{E}_{TC}, \mathbf{e+r}))$ where $X$ is a Riemann surface structure on $S$ and $(\mathcal{E}_{TC}, \mathbf{e+r})$ an element of the corresponding moduli space $\mathcal{M}_\mathbf{e}(TC, \langle e\rangle\oplus R)=\mathcal{M}^X_\mathbf{e}(TC, \langle e\rangle\oplus R)$ for this Riemann surface $X$. These pairs are considered up to the action diffeomorphisms of $S$ isotopic to the identity and the action of the Gauge group on $\mathcal{M}_\mathbf{e}(TC, \langle e\rangle\oplus R)$. In order to give this space the structure of a moduli space, we note that each fiber naturally admits a finite-to-one proper map $F$ to $\mathcal{M}(G)=\mathcal{M}^X(G)$ from Proposition \ref{prop:properFInite}. Hence, the joint moduli space admits a finite-to-one proper map $\mathcal{J}F$ to the joint moduli space of $G$-Higgs bundles $\mathcal{JM}(G)$, constructed for instance in \cite{CTW25}. Through $\mathcal JF$, we can constant the following holonomy map.

\begin{thm}
\label{thm:JointModuliHalf}
    The holonomy map, when restricted to the subspace of stable pairs is a half dimensional immersion into the character variety:
$$\mathrm{hol} : \mathcal{JM}_\mathbf{e}(TC, \langle e \rangle\oplus R)\to \mathrm{Hom}(\pi_1S,G)\sslash G.$$

Moreover, the image of stable pairs are infinitesimal irreducible representation into $G$, that is, the Lie algebra of centralizer of the representation in $G$ is trivial. In particular, if $G=\SL(n+1,\C)$, then the image of stable pairs are  irreducible representation into $G$.
\end{thm}

\begin{proof}
    We use  \cite[{Proposition 6.10}]{CTW25}.
    Let $(X, \mathcal{E}_{TC},\mathbf{e}+\mathbf{r})$ be a stable pair. By definition, $(X,\mathcal{E}_{TC})$ reduces to a~$TC$-principal bundle. Hence we need to find a~$TC$-invariant splitting of the Lie algebra $\mathfrak g=U_1\oplus U_2$ such
    that the Higgs field $\mathbf{e}+\mathbf{r}=(\varphi_1,\varphi_2)$
    with $\varphi_1$ nowhere vanishing. The splitting to choose here is simply $\mathfrak g=\mathfrak{g}_2\oplus (\oplus\mathfrak g_i,i=-2,\ldots,1)$, as $\mathbf{e}$ is nowhere vanishing. 

    By Proposition \ref{prop:StableEquiv}, stable $(TC,\langle e\rangle+R)$-Higgs pairs are mapped to stable $G$-Higgs bundles under the structure extension map $F$. From~\cite[{Proposition 6.10}]{CTW25} we deduce that image of stable pairs under $F$ are
    in the open strata $M_0$, and stable cyclic bundles. By Proposition \ref{prop:properFInite}, we obtain the image in $\mathcal{JM}_\mathbf{e}(TC, \langle e \rangle\oplus R)$ under the structure extension map $F$ inside $\mathcal{JM}(G)$ is of the same dimension as $\mathcal{JM}_\mathbf{e}(TC, \langle e \rangle\oplus R)$, which is of half dimension of the character variety $\chi(\pi_1S,G)$ from Proposition \ref{prop:Moduli}. Since all of our bundles are minimal surfaces, we apply~\cite[{Theorem G}]{CTW25} states that the holonomy map restricted on the intersection of $M_0$ with stable minimal surfaces is an immersion, as desired.

    The moreover part follows from the fact that holonomy of stable $G$-Higgs bundles are infinitesimal irreducible representations, see \cite[Theorem 3.32]{GPR}. In particular, when $G=\SL(n+1,\C)$, from \cite[Remark 4.4]{GPR}, every stable $G$-Higgs bundle satisfies $\operatorname{Aut}(\mathcal E,\varphi)=1$ and corresponds to irreducible representations by \cite[Theorem 7.2]{GPR}.
\end{proof}

\subsection{Proper slice in the character variety}

The domination with respect to weights allows us to prove the following. Let $G$ be a simple complex group with finite center of type different from $E_8$.

\begin{thm}\label{thm:JointModuliproper}
The holonomy map is proper:
$$\mathrm{hol}:\mathcal{JM}_\mathbf{e}(TC,\langle e\rangle \oplus R)\to \mathrm{Hom}(\pi_1S,G)\sslash G.$$

As a result, the image is a closed subset. 
\end{thm}

\begin{proof}
Let $(X_k, (\mathcal{E}_k, \varphi_k))$ be a sequence of elements in $\mathcal{JM}_\mathbf{e}(TC,\langle e\rangle \oplus R)$ such that the corresponding holonomies $[\rho_k]$ remain in a compact subset $K$ of the character variety. First assume that $(X_k)$ is diverging in Teichmüller space. Then up to taking a subsequence there exists a simple closed curve $\gamma$ whose length for the hyperbolic metric $\ell_{\text{hyp},X_k}(\gamma)$ on $X_k$ tends to $+\infty$. The projected metric has the property of being larger than the hyperbolic metric, see Theorem \ref{thm:NormBound3-cyclic}. Hence the length of $\gamma$ for the projected metric $g_{\pi,k}$ on $X_k$ satisfies $\ell_{g_{\pi,k}}(\gamma)\geq\ell_{\text{hyp},X_k}(\gamma)$. Finally using Theorem \ref{thm:Undistorted weight on the disk} we see that $\omega_V(\lambda(\rho((\gamma))))\geq \ell_{\text{hyp},X_k}(\gamma)$. 
This contradicts the fact that the representations $([\rho_k])$ remain in a compact subset of the character variety. 

\medskip
Hence $(X_k)$ is bounded in Teichmüller space,
so remains in a compact set $K_1$.
Let's consider $\rm{NAH}^{-1}(K)$ for $\rm{NAH}$ the nonabelian
Hodge map on the joint moduli space of Higgs bundles
considered  in \cite{CTW25}.
From the fact that $\rm{NAH}$ is analytic \cite[{Theorem 4.23}]{CTW25}, we deduce that $\rm{NAH}^{-1}(K)$ is closed.
Also its restriction on each fiber is proper, as the nonabelian Hodge map for $G$-Higgs bundles is a fiberwise homeomorphism and the structure extension map $F$ is proper from Proposition \ref{prop:properFInite}.
So we deduce that the set $\{(X,\mathcal E,\varphi): X\in K_1,\rm{NAH}(X,\mathcal E,\varphi)\in K\}$ is fiberwise compact.
As the map from the joint moduli space to the character variety is a fibration by~\cite[{Theorem A}]{CTW25}, we
deduce that the whole set is compact.

Up to subsequence, we deduce that $(X_k,(\mathcal E_k,\varphi_k))$ converges to $(X,(\mathcal E,\varphi))$,
hence the holonomy map is proper, as desired.
\end{proof}
\begin{rem}
    In \cite[Theorem G]{CTW25}, they also prove that the nonabelian Hodge map is
    a symplectic immersion, so the Slodowy slices define symplectic submanifolds of half-dimension in the character variety of $G$.
\end{rem}
\appendix

\section{Exceptional Lie algebras}
\label{sec:AppendixLieAlgebras}
In this section we describe the Slodowy slices for the root $\mathfrak{sl}_2$-triples for the exceptional Lie algebras. 

\medskip

In particular we describe for the long root $\mathfrak{sl}_2$-triples the adapted linear representation $V$ from Proposition \ref{prop:RepresentationV}, meaning that the weights of the longest coroot in $V$ are $-1,0$ and $1$. We also determine the weight $\omega_V\in \mathfrak{a}^*$ from Proposition \ref{prop:RepresentationV}.

\medskip

We use the coordinates and simple roots of
\cite[Table 1, page 292]{OnishchikVinberg1990}. This reference also gives us the list of fundamental weights. The representation $V$ is the representation  $R(\pi_1)$ in this reference, associated to the first fundamental weight.
All coordinate vectors below are elements of $\mathfrak a^*$.
For each chosen root $\gamma$, take a triple $(e,h,f)$ with
$e\in\mathfrak g_{-\gamma}$, $f\in\mathfrak g_{\gamma}$ and $h$ the coroot of $\gamma$, so that $[h,e]=-2e$ as in Section~\ref{sec:ParallelizableSlodowySlices}. In particular $R=\mathfrak g_1$ is the sum of the root spaces $\mathfrak g_\eta$ with $\eta(h)=1$ for the long root triples. We write $\mathfrak{h}$ for the Cartan subalgebra, which is the complexification of the restricted Cartan subalgebra $\mathfrak{a}$.

For a representation $\rho:\mathfrak g\to\mathfrak{gl}(V)$, we write its Cartan weight spaces and its $h$-eigenspaces as:
\[
 V[\omega]=\{v\in V:\rho(x)v=\omega(x)v\text{ for all }x\in\mathfrak h\},\]
 \[V_m=\ker(\rho(h)-m\operatorname{id}_V)
     =\bigoplus_{\omega(h)=m}V[\omega].
\]

\subsubsection*{\texorpdfstring{$G_2$}{G2}: long root}
\pdfbookmark[2]{G2:longroot}{G2long}
\label{subsec:AppendixG2}
Let $\alpha,\beta$ be the simple roots, with $\alpha$ long and $\beta$
short. The positive roots are
\[
 \Sigma^+=\{\alpha,\beta,\alpha+\beta,\alpha+2\beta,
                 \alpha+3\beta,2\alpha+3\beta\}.
\]
Choose $\gamma=\alpha$. Then :
\[
 \mathfrak c=\mathfrak{h}_\alpha
   \oplus\mathfrak g_{\alpha+2\beta}
   \oplus\mathfrak g_{-\alpha-2\beta}
   \cong\mathfrak{sl}_2(\mathbb C),
\]
\[
 R=\mathfrak g_{-\beta}\oplus\mathfrak g_{-\alpha-3\beta}
   \oplus\mathfrak g_{\alpha+\beta}
   \oplus\mathfrak g_{2\alpha+3\beta}.
\]

Here $\mathfrak{h}_\alpha$ is the kernel of $\alpha$ in $\mathfrak{h}$. This triple is parallelisable as it is the long root $\mathfrak{sl}_2$-triple. The action of $C$ on $R$ is the irreducible action $S^3\C^2$.

\medskip

The smallest nontrivial irreducible representation has dimension $7$.
Its weights, each with multiplicity one, are:
\[
 0,\quad\pm\beta,\quad\pm(\alpha+\beta),\quad\pm(\alpha+2\beta).
\]
Here $\omega_V=\alpha+2\beta$. Consequently:
\[
 V_1=V[\alpha+\beta]\oplus V[\alpha+2\beta].\]
 Hence $(\dim V_{-1},\dim V_0,\dim V_1)=(2,3,2)$.

 \medskip

 The $\alpha$-bundles from \cite{DE26} belong to the restricted Slodowy slice for this $\mathfrak{sl}_2$-triple.

\subsubsection*{\texorpdfstring{$G_2$}{G2}: short root}
Among root $\mathfrak{sl}_2$-triples, this is the only parallelisable
triple that is neither magical nor associated with a long root.
Choose $\gamma=\beta$. Then
\[
 \mathfrak c=\mathfrak h_\beta
   \oplus\mathfrak g_{2\alpha+3\beta}
   \oplus\mathfrak g_{-2\alpha-3\beta}
   \cong\mathfrak{sl}_2(\mathbb C),
 \qquad
 R=\mathfrak g_{-\alpha}\oplus\mathfrak g_{\alpha+3\beta}.
\]
Here $\mathfrak{h}_\beta$ is the kernel of $\beta$ on the Cartan subalgebra $\mathfrak{h}$. This triple is parallelisable as $[R,R]$ lies in the $0$ eigenspace of $h$ on $\mathfrak{g}$, which is equal to $\langle h\rangle \oplus \mathfrak{c}$. This is the only example of parallelisable triple that is not magical nor minimal, see Appendix \ref{appendix:Classification}.

\medskip

The Anosov property is proved for representations associated to certain Hodge bundles in this Slodowy slice in \cite{DE26}, where they are called $\beta$-bundles. However not all $\beta$-bundle are in the corresponding restricted Slodowy slice, as some contain a component of the Higgs field in the centralizer $\mathfrak{c}$.

\subsubsection*{\texorpdfstring{$F_4$}{F4}: long root}
Let $(\epsilon_1,\epsilon_2,\epsilon_3,\epsilon_4)$ be an orthonormal
basis of $\mathfrak a^*$. The $48$ roots are
\[
 \begin{aligned}
 \Sigma_{\mathrm{short}}
   &=\{\pm\epsilon_i:1\leq i\leq4\}
     \cup\left\{\tfrac12(\pm\epsilon_1\pm\epsilon_2
                                      \pm\epsilon_3\pm\epsilon_4)\right\},\\
 \Sigma_{\mathrm{long}}
   &=\{\pm\epsilon_i\pm\epsilon_j:1\leq i<j\leq4\}.
 \end{aligned}
\]
Fix the simple roots
\[
 \alpha_1=\tfrac12(\epsilon_1-\epsilon_2-\epsilon_3-\epsilon_4),\quad
 \alpha_2=\epsilon_4,\quad
 \alpha_3=\epsilon_3-\epsilon_4,\quad
 \alpha_4=\epsilon_2-\epsilon_3,
\]
and let $\overline{\mathfrak a^+}=\{x:\alpha_i(x)\geq0\text{ for all }i\}$.
Choose the highest root $\gamma=\epsilon_1+\epsilon_2$.
Then $R$ is the sum of the $14$ root spaces indexed by :
\[
 \begin{gathered}
 \epsilon_1,\ \epsilon_2,\\
 \epsilon_i\pm\epsilon_j\qquad(i=1,2;\ j=3,4),\\
 \tfrac12(\epsilon_1+\epsilon_2\pm\epsilon_3\pm\epsilon_4).
 \end{gathered}
\]
The centralizer is the kernel $\mathfrak{h}_\gamma$ of $\gamma$ in $\mathfrak{h}$ plus the root spaces associated to the $18$ roots:
\[
 \begin{gathered}
 \pm\epsilon_3,\ \pm\epsilon_4,\ \pm\epsilon_3\pm\epsilon_4,\qquad
 \pm(\epsilon_1-\epsilon_2),\\
 \tfrac12\bigl(\pm(\epsilon_1-\epsilon_2)
                         \pm\epsilon_3\pm\epsilon_4\bigr).
 \end{gathered}
\]
Set $u_1=\tfrac12(\epsilon_3+\epsilon_4)$, $u_2=\tfrac12(\epsilon_3-\epsilon_4)$ and $u_3=\tfrac12(\epsilon_1-\epsilon_2)$, which are orthogonal of the same length and span $\gamma^\perp$. The $18$ roots of $\mathfrak c$ are then exactly $\pm2u_i$ (namely $\pm\epsilon_3\pm\epsilon_4$ and $\pm(\epsilon_1-\epsilon_2)$) and $\pm u_i\pm u_j$ for $i\neq j$ (namely $\pm\epsilon_3$, $\pm\epsilon_4$ and $\tfrac12(\pm(\epsilon_1-\epsilon_2)\pm\epsilon_3\pm\epsilon_4)$), which is the root system $C_3$.

Writing $\epsilon_1=\tfrac12\gamma+u_3$ and $\epsilon_2=\tfrac12\gamma-u_3$, the restrictions to $\mathfrak h_\gamma$ of the $14$ roots indexing $R$ are $\pm u_1,\pm u_2,\pm u_3$ and the eight weights $\pm u_1\pm u_2\pm u_3$. The latter form the Weyl orbit of the third fundamental weight $u_1+u_2+u_3$ of $\mathfrak{sp}_6$, whose irreducible representation has dimension $14$. Hence $R$ is this representation, namely $\Lambda^3_0\mathbb C^6=\ker(\Lambda^3\mathbb C^6\to\mathbb C^6)$, the kernel of the contraction by the symplectic form, see \cite[{Section 17.2}]{FultonHarris1991}.

\medskip

The smallest nontrivial irreducible representation $V$ has dimension
$26$ and highest weight $\omega_V=\epsilon_1$. Its nonzero weights are the
$24$ short roots, each with multiplicity one; the zero weight has
multiplicity two.
\[
 V_1=V[\epsilon_1]\oplus V[\epsilon_2]
   \oplus\bigoplus_{s,t\in\{\pm1\}}
       V\!\left[\tfrac12(\epsilon_1+\epsilon_2+s\epsilon_3+t\epsilon_4)\right].
\]
In particular $d=\dim V_1=6$, and $(\dim V_{-1},\dim V_0,\dim V_1)=(6,14,6)$. 

\subsubsection*{\texorpdfstring{$F_4$}{F4}: short root}
Choose $\gamma=\epsilon_1$, so that
$\eta(h)=2(\eta,\epsilon_1)$. Then $R=R_2\oplus R_3$, where
\[
 \begin{aligned}
 R_2&=\bigoplus_{s,t,u\in\{\pm1\}}
        \mathfrak g_{\frac12(\epsilon_1+s\epsilon_2+t\epsilon_3+u\epsilon_4)},
       &\dim R_2&=8,\\
 R_3&=\bigoplus_{i=2}^4
       (\mathfrak g_{\epsilon_1+\epsilon_i}
         \oplus\mathfrak g_{\epsilon_1-\epsilon_i}),
       &\dim R_3&=6.
 \end{aligned}
\]
The centralizer is
\[
 \mathfrak c=\mathfrak h_\gamma\oplus
   \bigoplus_{\substack{2\leq i<j\leq4\\s,t\in\{\pm1\}}}
          \mathfrak g_{s\epsilon_i+t\epsilon_j}.
\]
Its root system is of type $D_3\cong A_3$, and $\dim\mathfrak c=15$.

This triple is not parallelisable. More explicitly, one can check that $[R_3,\theta(R_2)]$ intersects $ \mathfrak c\oplus \langle h\rangle$ trivially as they live in different eigenspaces for $h$, and $[R_3,\theta(R_2)]$ contains the rootspace associated to $\frac{1}{2}(\epsilon_1-\epsilon_2+\epsilon_3+\epsilon_4)$, which is the bracket of the root spaces associated to $\frac{1}{2}(-\epsilon_1+\epsilon_2+\epsilon_3+\epsilon_4)$ and $\epsilon_1-\epsilon_2$.

\subsubsection*{\texorpdfstring{$E_6$}{E6}: root triple}
Let $\epsilon_1,\ldots,\epsilon_6$ satisfy
\[
 \sum_{i=1}^6\epsilon_i=0,\qquad
 (\epsilon_i,\epsilon_j)=\delta_{ij}-\tfrac16.
\]
Let $\epsilon$ be orthogonal to their span, with :
$(\epsilon,\epsilon)=\tfrac12$. Together they span $\mathfrak a^*$;
the $\epsilon_i$ are not linearly independent.
The $72$ roots are :
\[
 \begin{gathered}
 \epsilon_i-\epsilon_j\quad(1\leq i\neq j\leq6),\\
 \epsilon_i+\epsilon_j+\epsilon_k\pm\epsilon
                   \quad(1\leq i<j<k\leq6),\qquad
 \pm2\epsilon.
 \end{gathered}
\]
Choose the simple roots
\[
 \alpha_i=\epsilon_i-\epsilon_{i+1}\quad(1\leq i\leq5),\qquad
 \alpha_6=\epsilon_4+\epsilon_5+\epsilon_6+\epsilon,
\]
with the highest root $\gamma=2\epsilon$.
Then
\[
 R=\bigoplus_{1\leq i<j<k\leq6}
          \mathfrak g_{\epsilon_i+\epsilon_j+\epsilon_k+\epsilon},
 \qquad \dim R=20,
\]
and
\[
 \mathfrak c=\mathfrak h_\gamma\oplus
      \bigoplus_{1\leq i\neq j\leq6}\mathfrak g_{\epsilon_i-\epsilon_j}.
\]
The roots of $\mathfrak c$ are the $30$ roots $\epsilon_i-\epsilon_j$ with $1\leq i\neq j\leq6$, where $\sum_i\epsilon_i=0$ and $(\epsilon_i,\epsilon_j)=\delta_{ij}-\tfrac16$. This is the root system $A_5$ of $\mathfrak{sl}_6(\mathbb C)$, in which the $\epsilon_i$ are the weights of $\mathbb C^6$.

Since $\epsilon$ vanishes on $\mathfrak h_\gamma=\ker(2\epsilon)$, the weights of $R$ are the $20$ distinct weights $\epsilon_i+\epsilon_j+\epsilon_k$ with $i<j<k$. These are exactly the weights of the irreducible representation $\Lambda^3\mathbb C^6$, which has dimension $20$, so $R\cong\Lambda^3\mathbb C^6$.

\medskip

Choose the $27$-dimensional irreducible representation of highest
weight $\omega_V=\epsilon_1+\epsilon$. Its weights all have multiplicity one and are:
\[
 \epsilon_i\pm\epsilon\quad(1\leq i\leq6),\qquad
 -\epsilon_i-\epsilon_j\quad(1\leq i<j\leq6),
\]
Zero is not a weight.
\[
 V_1=\bigoplus_{i=1}^6V[\epsilon_i+\epsilon]\]
 Hence $(\dim V_{-1},\dim V_0,\dim V_1)=(6,15,6)$.

\subsubsection*{\texorpdfstring{$E_7$}{E7}: root triple}
Let $\epsilon_1,\ldots,\epsilon_8$ span $\mathfrak a^*$ and satisfy
\[
 \sum_{i=1}^8\epsilon_i=0,\qquad
 (\epsilon_i,\epsilon_j)=\delta_{ij}-\tfrac18.
\]
The $126$ roots are
\[
 \epsilon_i-\epsilon_j\quad(1\leq i\neq j\leq8),\qquad
 \epsilon_i+\epsilon_j+\epsilon_k+\epsilon_\ell
                      \quad(1\leq i<j<k<\ell\leq8).
\]
The second family is already closed under negation, by taking
complementary four-element subsets.
Choose the simple roots :
\[
 \alpha_i=\epsilon_i-\epsilon_{i+1}\quad(1\leq i\leq6),\qquad
 \alpha_7=\epsilon_5+\epsilon_6+\epsilon_7+\epsilon_8.
\]
We choose the highest root $\gamma=\epsilon_8-\epsilon_7$.

Then $R$ is the sum of the $32$ root spaces indexed by
\[
 \begin{gathered}
 \epsilon_8-\epsilon_i,\quad\epsilon_i-\epsilon_7
                         \qquad(1\leq i\leq6),\\
 -(\epsilon_i+\epsilon_j+\epsilon_k+\epsilon_7)
                         \qquad(1\leq i<j<k\leq6).
 \end{gathered}
\]
The centralizer is the sum of the kernel $\mathfrak h_\gamma\subset \mathfrak{h}$ of $\gamma$ and the roots spaces associated to the $60$ roots :
\[
 \epsilon_i-\epsilon_j\quad(1\leq i\neq j\leq6),\qquad
 \pm(\epsilon_i+\epsilon_j+\epsilon_7+\epsilon_8)
                         \quad(1\leq i<j\leq6).
\]
Let $s=\sum_{k=1}^6\epsilon_k=-(\epsilon_7+\epsilon_8)$ and $u_i=\epsilon_i-\tfrac12 s$ for $1\leq i\leq6$. A direct computation gives $(u_i,u_j)=\delta_{ij}$ and $(u_i,\gamma)=0$. The roots of $\mathfrak c$ are $\epsilon_i-\epsilon_j=u_i-u_j$ and $\pm(\epsilon_i+\epsilon_j+\epsilon_7+\epsilon_8)=\pm(u_i+u_j)$, that is, the $60$ roots $\pm u_i\pm u_j$ of type $D_6$.

Computing the scalar products with the $u_l$, the roots $\epsilon_8-\epsilon_i$, $\epsilon_i-\epsilon_7$ and $-(\epsilon_i+\epsilon_j+\epsilon_k+\epsilon_7)$ restrict to $\tfrac12\sum_l\pm u_l$ with respectively $1$, $5$ and $3$ minus signs. Hence the weights of $R$ are the $32$ vectors $\tfrac12(\pm u_1\pm\cdots\pm u_6)$ with an odd number of minus signs. These are the weights of one of the two half-spin representations of $\mathfrak{so}_{12}$, which is irreducible of dimension $32$, so $R$ is this half-spin representation \cite[{Section 20.1}]{FultonHarris1991}.

\medskip

The smallest nontrivial irreducible representation $V$ has dimension
$56$ and highest weight $\omega_V=\epsilon_1+\epsilon_8$. Its weights are
\[
 \pm(\epsilon_i+\epsilon_j)\qquad(1\leq i<j\leq8),
\]
all with multiplicity one. We obtain
\[
 V_1=\bigoplus_{i=1}^6
        \bigl(V[\epsilon_i+\epsilon_8]\oplus V[-\epsilon_i-\epsilon_7]\bigr),
 \qquad d=12,
\]
and
\[
 (\dim V_{-1},\dim V_0,\dim V_1)=(12,32,12).
\]

\subsubsection*{\texorpdfstring{$E_8$}{E8}: root triple}
Let $\epsilon_1,\ldots,\epsilon_9$ span $\mathfrak a^*$ and satisfy
\[
 \sum_{i=1}^9\epsilon_i=0,\qquad
 (\epsilon_i,\epsilon_j)=\delta_{ij}-\tfrac19.
\]
The $240$ roots are
\[
 \epsilon_i-\epsilon_j\quad(1\leq i\neq j\leq9),\qquad
 \pm(\epsilon_i+\epsilon_j+\epsilon_k)\quad(1\leq i<j<k\leq9).
\]
Choose the simple roots
\[
 \alpha_i=\epsilon_i-\epsilon_{i+1}\quad(1\leq i\leq7),\qquad
 \alpha_8=\epsilon_6+\epsilon_7+\epsilon_8,
\]
and the highest root $\gamma=\epsilon_1-\epsilon_9$.

The space $R$ is the sum of the $56$ root spaces indexed by
\[
 \begin{gathered}
 \epsilon_1-\epsilon_i,\quad\epsilon_i-\epsilon_9
                           \qquad(2\leq i\leq8),\\
 -\epsilon_9-\epsilon_i-\epsilon_j,\quad
 \epsilon_1+\epsilon_i+\epsilon_j
                           \qquad(2\leq i<j\leq8).
 \end{gathered}
\]
The centralizer is the sum of $\mathfrak{h}_\gamma$ and the root spaces associated to the 126 roots:
\[
 \begin{gathered}
 \epsilon_i-\epsilon_j\qquad(2\leq i\neq j\leq8),\\
 \pm(\epsilon_i+\epsilon_j+\epsilon_k)
                           \qquad(2\leq i<j<k\leq8),\\
 \pm(\epsilon_1+\epsilon_i+\epsilon_9)\qquad(2\leq i\leq8).
 \end{gathered}
\]
The listed roots of $\mathfrak c$ form a reduced, simply laced root system of rank $7=\dim\mathfrak h_\gamma$ with $42+70+14=126$ roots. Among simply laced root systems of rank $7$, $A_7$ has $56$ roots and $D_7$ has $84$, so $\mathfrak c$ is of type $E_7$.

The space $R$ has dimension $56$. It has no zero weight, since a root $\eta$ with $\eta|_{\mathfrak h_\gamma}=0$ would be $\tfrac12\gamma$, which is not a root. Hence $R$ has no trivial summand. Since every nontrivial irreducible representation of $\mathfrak e_7$ has dimension at least $56$, $R$ is the $56$-dimensional representation, (which is also the representation $V$ when $G$ is of type $E_7$).

\medskip

The smallest nontrivial irreducible representation is the adjoint
representation, so we don't have a representation $V$ whose $h$-weights are $\lbrace 1,0,-1\rbrace$.

\section{Classification of parallelisable \texorpdfstring{$\mathfrak{sl}_2$}{sl2}-triples}
\label{appendix:Classification}

In this Appendix, we give the list of all parallelisable $\mathfrak{sl}_2$-triples. The key idea is that for each such parallelisable $\mathfrak{sl}_2$-triple one can associate a complex linear Lie algebra involution whose fixed points are exactly $\mathfrak{s}\oplus \mathfrak{c}$, the $\mathfrak{sl}_2$-Lie subalgebra and its centralizer. With this idea the classification becomes quite simple, using the classification of involutions of simple complex Lie algebras.

\subsection{The criterion}

In order to get this characterization of parralelisable triples we will use the following lemma.
\begin{lem}
    For $m\geq 0$, denote by $V_m=\rm{Sym}_m(\C^2)$ the $(m+1)$-th
    dimensional irreducible $\rm{SL}(2)$-module. Decomposing $V_m=\oplus_{k=0}^m W_m^{m+2k}$ into $h$-weight spaces, we have, for any $m,m'\geq 0$, $W_m^{m}\otimes W_{m'}^{-m'}$
    generates the $\mathfrak{sl}_2$-module $V_m\otimes V_{m'}$.
\end{lem}

Recall that a $\mathfrak{sl}_2$-module is a vector space $V$ with an action $\mathfrak{sl}_2\to \mathfrak{gl}(V)$. Given two $\mathfrak{sl}_2$-modules $V,W$, the action on the module $V\otimes W$ is such that for $x\in \mathfrak{sl}_2$, $v\in V$ and $w\in W$ : $x\cdot (v\otimes w)=(x\cdot v)\otimes w+v\otimes(x\cdot w)$. 
\begin{proof}
    Writing $\mathfrak{sl}_2=\langle e,f,h\rangle$, we have the following:
    $$\ad(e)^k(W_m^m)=W_m^{m-2k},\;\ad(e)(W_m^{-m'})=\lbrace 0\rbrace.$$
    Similarly:
    $$\ad(f)^k(W_{m'}^{-m'})=W_{m'}^{-m'+2k},\;\ad(f)(W_m^{m})=\lbrace 0\rbrace.$$
    As such:
    $$\ad(e)^k(W_m^m\otimes W_{m'}^{-m'})=W_m^{m-2k}\otimes W_{m'}^{-m'}.$$
    In particular the $\mathfrak{sl}_2$-module generated by $W_m^m\otimes W_{m'}^{-m'}$ contains
    $W_m\otimes W_{m'}^{-m'}$ and $W_m^m\otimes W_{m'}$.
    Applying $\ad(f)$ to $W_m\otimes W_{m'}^{-m'}$, and projecting
    orthogonally to $W_m\otimes W_{m'}^{-m'}$, this shows that
    $W_m\otimes W_{m'}^{-m'+2}$ is also in the generated module.
    Applying $\ad(f)$ successively and using the same argument, we deduce that the $\mathfrak{sl}_2$-module generated is $W_m\otimes W_{m'}$ as desired.
\end{proof}
\begin{prop}
    Let $\mathfrak s$ be an $\mathfrak{sl}_2$-triple in $\mathfrak g$. the following assertions are equivalent:
    \begin{enumerate}
        \item $\mathfrak s$ is parallelisable.
        \item $[\mathfrak m,\mathfrak m]\subset\mathfrak s\oplus\mathfrak c$, where $\mathfrak m$ is the orthogonal for the Killing form of $\mathfrak s\oplus \mathfrak c$.
        \item $\mathfrak s\oplus \mathfrak c$ is the fixed subalgebra of an involutive automorphism of $\mathfrak g$.
    \end{enumerate}
\end{prop}
\begin{proof}
    Note that by construction, the restricted slodowy slice $R$ is $R=\mathfrak m\cap\rm{ker}(ad_f)$, and $\Theta(\mathfrak m)=\mathfrak m$, so $(2)$ implies that $\mathfrak s$ is parallelisable.

    It is a classical fact that $(2)$ and $(3)$ are equivalent. Indeed on the one hand given $(2)$ one can define the involution to be the $\Id$ on $\mathfrak{s}\oplus \mathfrak{c}$ and $-\Id$ on $\mathfrak{m}$ and obtain $(3)$. On the other hand the $-1$-eigenspace for a Lie algebra involution must be orthogonal to its $1$-eigenspace, proving point $(2)$.

    It remains to show $(2)$ under the assumption that $\mathfrak s$
    is parallelisable. To this effect, decompose $\mathfrak m=\oplus M_a$ into irreducible $\mathfrak s$-modules, and observe that $R$
    is the direct sum of the highest-weight spaces in each of the $M_a$. Similarly, $\theta(R)$ is the sum of the lowest-space weights in $M_a$ if $\theta$ is such that $\theta(e)=-f$ and $\theta(h)=-h$.
    
    Denote by $T:\mathfrak m\otimes\mathfrak m\rightarrow \mathfrak m$ the tensor defined by $T(X,Y)=\rm{pr}_\mathfrak m([X,Y])$,
    the projection on $\mathfrak m$ of the Bracket for the Killing metric. The parallelizability assumption ensures that $T$ vanishes
    on $M_a^+\otimes M_b^-$ where $M_a^+$ is the highest weight space of $M_a$ and $M_b^-$ the lowest space weight of $M_b$.

    Also, the compatibility of the Killing form with the Bracket ensures that $T$ is $\mathfrak s$-equivariant, i.e.:
    \begin{equation*}
        \forall u\in\mathfrak s,X,Y\in \mathfrak m,\,T(\ad(u)\cdot X,Y)+T(X,\ad(u)\cdot Y)=\ad(u)\cdot T(X,Y)\,.
    \end{equation*}
    By the previous lemma, we deduce that $T=0$, i.e. $[\mathfrak m,\mathfrak m]\subset\mathfrak s\oplus \mathfrak c$, as claimed.
\end{proof}

\begin{rem}
    Note that this complex Lie algebra involution with fixed point $\mathfrak{s}\oplus \mathfrak{c}$ corresponds exactly with the real form $G^\sigma$ that we introduce, associated with the homogeneous space $G/G^\sigma$.
\end{rem}

\begin{rem}
    The difference with the real form associated with magical $\mathfrak{sl}_2$-triples is the following: for the involution associated to a parallelisable $\mathfrak{sl}_2$-triple, $R\subset \mathfrak{m}$ but $e\notin \mathfrak{m}$. For the involution associated to a magical $\mathfrak{sl}_2$-triple,$R\subset \mathfrak{m}$ and $e\in \mathfrak{m}$.
\end{rem}

\subsection{The classification}

We now state and prove the classification of parallelisable $\mathfrak{sl}_2$-triples.

\begin{thm}
    If $\mathfrak{s}\subset \mathfrak{g}$ is a parallelisable $\mathfrak{sl}_2$-triple in a simple complex Lie subalgebra, then either:
\begin{enumerate}
    \item the triple is the minimal $\mathfrak{sl}_2$-triple in $\mathfrak{g}$, 
    \item the triple is the principal triple  $\mathfrak{so}(2,1)\subset \mathfrak{sl}(3,\C)$,
    \item the triple is the non-principal magical triple $\mathfrak{so}(2,1)\subset \mathfrak{so}(n,\C)$ for $n\geq 5$,
    \item the triple is the short root $\mathfrak{sl}_2$-triple in the Lie algebra of type $G_2$.
\end{enumerate}
\end{thm}

The triples $(2)$ and $(3)$ are both parallelisable and \emph{magical} in the sense of \cite{BCGGO24}. The only non-minimal and non-magical case is the short root $\mathfrak{sl}_2$-triple in $G_2$. Representations in its Slodowy slice which also lie in its split real form $G_2'$ are studied in \cite{CT23, DE26b, DE26} ($\beta$-bundles).

\begin{rem}
Note that in this classification one must be careful that the minimal $\mathfrak{sl}_2$-triple is defined only up to automorphisms of the Lie algebras, giving in general $3$ minimal $\mathfrak{sl}_2$-triples in $\mathfrak{so}(8,\C)$ up to inner automorphisms. Note also that in this classification one needs to keep in mind the low dimensional isomorphisms, so that in $\mathfrak{sl}(4,\C)\simeq \mathfrak{so}(6,\C)$ there is a magical parallelisable $\mathfrak{sl}_2$-triple.
\end{rem}

In the proof, we use the classification of complex Lie algebra involutions from \cite[{Section 41.4}]{Eti24} for the classical cases and \cite[{Section 42}]{Eti24} for the exceptional cases. Note that complex involutions are in correspondence with real forms: given a complex involution, composing it with a commuting Cartan involution defines a real form.
\begin{proof}

In this proof we consider case by case all the possible fixed subalgebras of complex Lie algebra involutions. We determine all the cases in which this algebra is isomorphic to $\mathfrak{sl}(2,\C)\oplus \mathfrak{c}$ and we identify the $\mathfrak{sl}(2,\C)$ with one element of the classification. Note that we already proved case by case that all the triples in the classification are parallelisable.

\textbf{Type $A_{n-1}$:} Let $\mathfrak{g}=\mathfrak{sl}(n,\C)$, $n\geq 2$. The fixed Lie algebra of any complex involution has one of the following forms:
\begin{itemize}
    \item $\mathfrak{sl}(p,\C)\oplus \mathfrak{gl}(q,\C)$ for some $p+q=n$,
    \item $\mathfrak{so}(n,\C)$,
    \item $\mathfrak{sp}(2k,\C)$ if $n=2k$. 
\end{itemize}

This fixed sub-algebra has the correct shape only if it is $\mathfrak{sl}(2,\C)\oplus \mathfrak{gl}(n-2,\C)$ (case $(1)$), or $\mathfrak{so}(3,\C)\simeq \mathfrak{sl}(2,\C)$ if $n=3$ (case $(2)$) or $\mathfrak{so}(4,\C)$ if $n=4$ (case $(3)$ using $\mathfrak{sl}(4,\C)\simeq \mathfrak{so}(6,\C)$). 

\medskip

\textbf{Type $B_n$:} Let $\mathfrak{g}=\mathfrak{so}(2n+1,\C)$. 
The fixed Lie algebra of any complex involution has the form:
\begin{itemize}
    \item $\mathfrak{so}(2p+1,\C)\oplus \mathfrak{so}(2q,\C)$ for some $p+q=n$.
\end{itemize}

This fixed sub-algebra has the correct shape only if it is $\mathfrak{so}(3,\C)\oplus \mathfrak{so}(2n-2,\C)$ (case $(3)$), or $\mathfrak{so}(4,\C)\oplus \mathfrak{so}(2n-3,\C)$ (case $(1)$). In the last case note that $\mathfrak{so}(4,\C)=\mathfrak{sl}(2,\C)\oplus \mathfrak{sl}(2,\C)$, and these factors can be exchanged by Lie algebra automorphisms. 

\medskip

\textbf{Type $C_n$:} Let $\mathfrak{g}=\mathfrak{sp}(2n,\C)$. The fixed Lie algebra of any complex involution has the form:
\begin{itemize}
    \item $\mathfrak{sp}(2p,\C)\oplus \mathfrak{sp}(2q,\C)$ for some $p+q=n$,
    \item $\mathfrak{gl}_n(\C)$.
\end{itemize}

This fixed sub-algebra has the correct shape only if it is $\mathfrak{sp}(2,\C)\oplus \mathfrak{sp}(2n-2,\C)$ (case $(1)$), or $\mathfrak{gl}_n(\C)$ if $n=2$ (case $(3)$).

\medskip

\textbf{Type $D_n$:} Let $\mathfrak{g}=\mathfrak{so}(2n,\C)$, $n\geq 3$. The fixed Lie algebra of any complex involution has the following form, up to Lie algebra automorphisms (including the triality in the case $n=4$):
\begin{itemize}
    \item $\mathfrak{so}(2p,\C)\oplus \mathfrak{so}(2q,\C)$ for some $p+q=n$,
    \item $\mathfrak{so}(2p+1,\C)\oplus \mathfrak{so}(2q-1,\C)$ for some $p+q=n$,
    \item $\mathfrak{gl}_n(\C)$.
\end{itemize}

This fixed sub-algebra has the correct shape only if it is $\mathfrak{so}(3,\C)\oplus \mathfrak{so}(2n-3,\C)$ (case $(3)$), or $\mathfrak{so}(4,\C)\oplus \mathfrak{so}(2n-4,\C)$ (case $(1)$). 

\medskip

\textbf{Exceptional cases:}
In the following table we write down the type of the fixed subalgebras of complex Lie algebra involutions for exceptional simple Lie algebras (here we exclude the trivial involution $\Id$):
\begin{center}
    \begin{tabular}{|c|c|c|c|c|}
        \hline
        $G_2$ & $F_4$& $E_6$& $E_7$& $E_8$  \\
        \hline
                 $A_1+A_1$& $B_4$&$C_4$&$A_7$&$D_8$\\

         &$C_3+A_1$ &$D_5+\C$&$E_6+\C$&$E_7+A_1$\\

           & &$F_4$&$D_6+A_1$&\\

           & &$A_5+A_1$&&\\
          \hline
    \end{tabular}
\end{center}

In each column only the last case can have the form $\mathfrak{s}\oplus \mathfrak{c}$. Note that in the case $G_2$ one can choose which one of the $A_1$-factor is $\mathfrak{s}$, giving us at most two different $\mathfrak{sl}_2$-triples. Since we exhibited for each of these groups one parallelisable $\mathfrak{sl}_2$-triple and $2$ different ones in the case $G_2$ we have them all.

\end{proof}
\newpage

\bibliographystyle{alpha}
\bibliography{references}
\end{document}